\documentclass{m2an}

\usepackage{amsmath}
\usepackage{amssymb,amsfonts,amsthm,bm}
\usepackage[usenames,dvipsnames]{color}
\usepackage{graphicx,caption,subfig}
\usepackage{xspace}
\usepackage{placeins}
\newtheorem{thm}{Theorem}[section]
\newtheorem{prop}[thm]{Proposition}
\newtheorem{lem}[thm]{Lemma}

\newtheorem{remark}[thm]{Remark}

\newcommand{\myinclude}[2]{{\includegraphics{pics/ellipses_a_N10_Ex3}}}

\newcommand{\sm}{\hspace{-2pt}-\hspace{-2pt}}
\renewcommand{\sp}{\hspace{-2pt}+\hspace{-2pt}}

\newcommand{\cK}{{\mathcal K}}
\newcommand{\cF}{{\mathcal F}}
\newcommand{\jl}{[\![}
\newcommand{\jr}{]\!]}
\newcommand{\jmp}[1]{\jl#1\jr}
\newcommand{\al}{\{\hspace{-3.5pt}\{}
\newcommand{\ar}{\}\hspace{-3.5pt}\}}
\newcommand{\avg}[1]{\al#1\ar}

\newcommand{\cd}{\hspace{-2pt}\cdot\hspace{-2pt}}
\newcommand{\sscal}[2]{(#1,#2)_{\star,\kappa}}

\newcommand{\tnorm}[1]{|\!|\!| #1|\!|\!|}
\newcommand{\snorm}[1]{\|#1\|_{\star,\kappa}}

\newcommand{\uN}{u_N}
\newcommand{\vN}{v_N}

\newcommand{\real}{\mathbb R}
\newcommand{\Hper}{H^1_\#(\Omega)}
\newcommand{\VN}{\mathbb V_N}

\newcommand{\grad}{\nabla}

\newcommand{\pz}{\Pi^\kappa_0}
\newcommand{\pN}{\Pi_N^\kappa}

\newcommand{\plus}{\raisebox{.4\height}{\scalebox{.6}{+}}}
\newcommand{\minus}{\raisebox{.4\height}{\scalebox{.8}{-}}}
\newcommand{\Vp}{V_{\pmb{\plus}\,}}
\newcommand{\Vm}{V_{\pmb{\minus}\,}}
\newcommand{\bub}{g}

\newcommand{\etaR}{\eta_{{\tt R},\kappa}}
\newcommand{\etaF}{\eta_{{\tt F},\kappa}}
\newcommand{\etaJ}{\eta_{{\tt J},\kappa}}
\newcommand{\ccR}{c_{{\tt R},\kappa}}
\newcommand{\ccF}{c_{{\tt F},\kappa}}
\newcommand{\ccJ}{c_{{\tt J},\kappa}}

\newcommand{\bk}{\beta_\kappa}

\newcommand{\gk}{\gamma_\kappa}
\newcommand{\gF}{\gamma_{\tt F}}
\newcommand{\dku}{{\tt d}_\kappa^u}

\newcommand{\dkpu}{{\tt d}_{\kappa'}^u}

\newcommand{\dkN}{{\tt d}_{\kappa}}
\newcommand{\rak}{{\tt a}_\kappa}
\newcommand{\rbk}{{\tt b}_\kappa}
\newcommand{\rbkp}{{\tt b}_{\kappa'}}
\newcommand{\rbok}{{\tt b}_{\omega(\kappa)}}
\newcommand{\rck}{{\tt c}_\kappa}

\makeatletter
\newcommand{\pushright}[1]{\ifmeasuring@#1\else\omit\hfill$\displaystyle#1$\fi\ignorespaces}
\newcommand{\pushleft}[1]{\ifmeasuring@#1\else\omit$\displaystyle#1$\hfill\fi\ignorespaces}
\makeatother

\newcommand{\REV}[1]{{{#1}}}

\begin{document}

\title{A posteriori error estimates for discontinuous Galerkin methods
using non-polynomial basis functions. \\Part I: Second order linear PDE}
\author{Lin Lin}
\address{
Department of Mathematics, University of California Berkeley
and Computational Research Division, Lawrence Berkeley National
Laboratory, Berkeley, CA 94720. Email: linlin@math.berkeley.edu}
\author{Benjamin Stamm}
\address{
Sorbonne Universit\'{e}s, UPMC Univ. Paris 06, UMR 7598, Laboratoire
Jacques-Louis Lions, F-75005 Paris, France. Email: stamm@ann.jussieu.fr
}



\date{}

\begin{abstract}

We present the first systematic work for deriving a posteriori error estimates for general non-polynomial basis functions in an interior penalty discontinuous Galerkin (DG) formulation for solving second order linear PDEs.  Our residual type upper and lower bound error estimates measure the error in the energy norm. The main merit of our method is that the method is parameter-free, in the sense that all but one solution-dependent constants appearing in the upper and lower bound estimates are explicitly \REV{computable by solving local eigenvalue problems}, and the only non-computable constant can be reasonably approximated by a computable one without affecting the overall effectiveness of the estimates in practice.  As a side product of our formulation, the penalty parameter in the interior penalty formulation can be automatically determined as well.  We develop an efficient numerical procedure to compute the error estimators.  Numerical results for a variety of problems in 1D and 2D demonstrate that both the upper bound and lower bound are effective.
\end{abstract}

\maketitle

%

\section{Introduction}

Let $\Omega$ be a bounded domain. We consider the development of a posteriori error
estimates for the following second order linear PDE
\begin{equation}
  \label{eqn:problem}
	-\Delta u + V u = f,\qquad\mbox{in }\Omega, 
\end{equation}
using the discontinuous Galerkin (DG) formulation with general
non-polynomial basis sets. 

Such equation arises in many scientific and engineering problems
such as in electromagnetism, geophysics, quantum physics, to name a few.
In order to solve Eq.~\eqref{eqn:problem} in practice, it is desirable
to reduce the number of degrees of freedom for discretizing
Eq.~\eqref{eqn:problem} to have a smaller algebraic problem to solve.
While standard polynomial basis functions can approach a complete basis
set and is versatile enough to represent almost any function of
interest, the resulting number of degrees of freedom is usually large even when
high order polynomials are used.  Non-polynomial basis functions are
therefore often employed to reduce the number of degrees of freedom, and
are widely used to solve Eq.~\eqref{eqn:problem} and other equations,
including the planewave basis set for solving Helmholtz
equation~\cite{HiptmairMoiolaPerugia2011,TezaurFarhat2006}, the
heterogeneous multiscale method (HMM)~\cite{EEngquist2003} and the
multiscale finite element method~\cite{HouWu1997} for solving multiscale
elliptic equations, and the various non-polynomial basis set used in
quantum chemistry such as the Gaussian basis
set~\cite{FrischPopleBinkley1984}, atomic orbital basis
set~\cite{Junquera:01}, and adaptive local basis
set~\cite{LinLuYingE2012} etc.  

Besides solving the equation, it is also often desirable to assess the
accuracy of the numerical solution via a posteriori error estimates and to design approximation spaces that result in a uniform distribution of the error in space to achieve best accuracy for a given number of degrees of freedom. 
\REV{
In this paper we focus on the a posteriori error estimates in the
interior penalty DG formulation~\cite{Wheeler78,Arnold1982,Nitsche71,BabuskaZlamal:73,DD76,Baker77}. 
}

The DG formulation has the advantage that it formally relaxes the continuity constraint of
basis functions at the inter-element boundary, and is therefore
particularly suitable for incorporating general basis functions, which
are difficult to match at the inter-element boundary.


\subsection{Previous work}

Compared to the many existing works on a posteriori error estimates
using polynomial basis functions in the DG
formulation~\cite{HoustonSchotzauWihler2007,KarakashianPascal2003,SchotzauZhu2009},
it is much more difficult to develop systematic a posteriori error
estimates for general non-polynomial basis functions. One of the
important reasons is that approximation and scaling properties of the function space
spanned by non-polynomial basis functions, which are key to a posteriori
error estimates, are generally difficult to deduce. For instance, Amara et
al~\cite{AmaraDjellouliFarhat2009} developed the upper bound error
estimates for the Helmholtz equation in planewave basis enriched DG
method, and the error is measured in the $L^{2}$--norm. 
Kaye et al~\cite{KayeLinYang2014} developed the upper bound error
estimates for solving linear eigenvalue problems using
non-polynomial basis functions in a DG framework, which generalizes the
work of Giani et al~\cite{GianiHall2012} for polynomial basis functions.
However, the assumption of approximation properties on the function
space is in general difficult to verify.
Though not in the DG framework, Ohlberger~\cite{Ohlberger2005} developed the a
posteriori error estimates for the HMM method for elliptic
homogenization problems.

The difficulties of a posteriori error analysis for general
non-polynomial basis functions are largely due to the lack of credible
methods for measuring the ratio of the error using different norms, defined
in proper function spaces. For instance, approximately speaking, in a
residual based error estimator, the constants associated with the
residual requires the estimation of ratio of the error measured using
$L^{2}$--norm and the $H^{1}$--norm. 
The scaling properties of such constants with respect to the increase of
the number of basis functions on a particular element can
be rather intrigue for non-polynomial basis functions. The estimation of
such constants is already complicated for polynomial basis functions or
planewave basis functions, not to mention the case when the
non-polynomial basis functions come from numerical solution without a
analytic recipe, or even worse, the basis functions do not in practice
form a complete basis set with only saturating accuracy.

\subsection{Contribution}

To the extent of our knowledge, this is the first systematic work for
deriving a posteriori error estimates for general non-polynomial basis
functions in a DG framework. Our upper and lower bound error estimates are
residual type estimators for the error in the energy norm.
In our formulation, all but one basis-dependent constants appearing in
the upper and lower bound estimates are explicitly computable
\REV{by solving local eigenvalue
problems}. For
solution with sufficient regularity (for instance $u\in H^{2}(\Omega)$),
the only non-computable constant can be reasonably approximated by a
computable one without affecting the overall effectiveness of the
estimates. \REV{While the requirement of $H^{2}(\Omega)$
regularity appears to be a formal drawback in the context of a
posteriori error estimates, the main goal of this work is to develop a
posteriori error estimates for general basis sets rather than for
$h$-refinement, and the difficulty of general basis sets holds even if
the solution has $C^{\infty}(\Omega)$ regularity. Therefore we think our method can
have important practical values. }
As a side product of our results, the penalty parameter in the interior
penalty formulation is also automatically computed, and the computed
constants guarantees that the coercivity of the resulting DG bilinear
form for Poisson's equation.

We develop an efficient numerical procedure to compute these
constants.  Both the formulation and the practical implementation of our
method are independent of how the basis functions are generated.
Although the numerical procedure is developed for general non-polynomial
basis functions, we find that the procedure, when applied to standard
polynomial basis functions, generates constants are even more accurate
than the analytical asymptotic result.  Numerical results for a variety
of problems in 1D and 2D indicate that both the upper bound and lower
bound are sharp, and the effectiveness of the estimators holds even at
the level of each element.




\subsection{Outline}

The rest of the paper is organized as follows. After an introduction to some technical
results in section \ref{sec:Prelim}, we start with the derivation of the
upper bound a posteriori error estimates for the Poisson's equation in
section \ref{sec:Poisson}, without the potential term $V$.  We then
generalize the derivation of the upper bound error estimates to
indefinite problems with the potential term, as well as the lower bound
error estimates in section~\ref{sec:Indefinite}. We elaborate in section
\ref{sec:computation} on the numerical methods for computing the
constants appearing in the upper and lower bound estimates needed in our
analysis.  Finally, we present numerical results in
section~\ref{sec:numerical}, before we conclude in
section~\ref{sec:Conclusion} followed by an appendix.

\section{Preliminary results}
\label{sec:Prelim}
\subsection{Mesh, broken spaces, jump and average operators}
Let $\Omega=(0,1)^d$, $d=1,2,3$ and let $\cK$ be a regular partition of
$\Omega$ into elements $\kappa\in\cK$.  That is, we assume that
\REV{the interior of $\overline{\kappa}\cap\overline{\kappa}'$, for any $\kappa,\kappa'\in \cK$, is either an
element of $\cK$,} a common face, edge, vertex of the partition or the
empty set.
For simplicity, we identify the boundary of $\Omega$ in a periodical
manner. That means, that we also assume the partition to be regular
across the boundary $\partial \Omega$.  
\REV{We remark that although the assumption of a rectangular
domain with periodic boundary condition appears to be 
restrictive, such setup already directly finds its application in
important areas such as quantum chemistry and materials science.}
However, the analysis below is not restricted to equations with periodic
boundary condition.  Other boundary conditions, such as Dirichlet or
Neumann boundary conditions can be employed as well with minor
modification.  \REV{Generalization to non-rectangular domain does not
introduce conceptual difficulties either, but may lead to changes in numerical
schemes for estimating relevant constants in
section~\ref{sec:computation}, if the tensorial structure of the grid
points is not preserved.}

Let $N=(N_\kappa)_{\kappa\in\cK}$ denote the vector of the local number of degrees of freedom $N_\kappa$ on each element $\kappa\in\cK$.
Let $\VN=\bigoplus_{\kappa\in\cK}\VN(\kappa)$ by any piecewise discontinuous approximation space on a partition $\cK$ of the domain $\Omega$. 
\REV{
It is important to highlight that 
little is known about the a priori information of $\VN$
except that we assume that each $\VN(\kappa)$ contains
constant functions and that $\mathbb{V}_{N}(\kappa)\subset
H^{\frac32}(\kappa)$, so that the traces of $\grad v_N$ on the boundary
$\partial \kappa$ are well-defined for all $v_N\in\VN(\kappa)$, for all
$\kappa\in \cK$.
}
We denote by $H^s(\kappa)$ the standard Sobolev space of $L^2(\kappa)$-functions such that all partial derivatives of order $s\in\mathbb N$ or less lie as well in $L^2(\kappa)$.
By $H^s(\cK)$, we denote the set of piecewise $H^s$-functions defined by
\[
	H^s(\cK) = \left\{ v\in L^2(\Omega) \,\middle|\, v|_\kappa \in H^s(\kappa), \forall \kappa \in \cK \right\},
\]
also referred to as the broken Sobolev space.
We denote by $\Hper$ the space of periodic $H^1$-functions on $\Omega$.
We further define the element-wise resp. face-wise scalar-products and norms as
\[
	(v,w)_\cK = \sum_{\kappa\in\cK} (v,w)_\kappa
	\qquad\mbox{and}\qquad
	\|v\|_\cK = (v,v)_\cK^\frac12.
\]
The $L^2$-norm on $\kappa$ and $\Omega$ are denoted by $\|\cdot\|_\kappa$ and $\|\cdot\|_\Omega$, respectively.

The jump and average operators on a face \REV{$\overline{F}=\overline{\kappa}\cap\overline{\kappa}'$} are defined in a standard manner by
\begin{align*}
			&\avg{v} = \tfrac12(v|_{\kappa} + v|_{\kappa'}),
			\qquad&&\mbox{and}\qquad\qquad
			\jmp{v} = v|_{\kappa} n_\kappa + v|_{\kappa'} n_{\kappa'},\\
			&\avg{\grad v} = \tfrac12(\grad v|_{\kappa} + \grad v|_{\kappa'}),
			\qquad&&\mbox{and}\qquad\qquad
			\jmp{\grad v} = \grad v|_{\kappa} n_\kappa +\grad v|_{\kappa'} n_{\kappa'},
\end{align*}
where $n_\kappa$ denotes the exterior unit normal of the element $\kappa$. 
Further we state some standard results


\REV{
Finally we recall the standard result of piecewise integration by parts formula that will be employed several times in the upcoming analysis.
\begin{lem}
	\label{lem:IntByParts}
	Let $v,w\in H^{2}(\cK)$. Then, there holds
	\[
		\sum_{\kappa\in\cK} \Big[
			(\Delta v,w)_{\kappa} + (\grad v,\grad w)_{\kappa}
		\Big]
		= 
		\tfrac12
		\sum_{\kappa\in\cK} \Big[
			(\jmp{\grad v},w)_{\partial\kappa} + (\grad v,\jmp{w})_{\partial\kappa}
		\Big].
	\]
\end{lem}
}

\subsection{Projections}
For any element $\kappa\in\cK$, let us denote by $\pz:L^2(\kappa)\to \real$ the $L^2(\kappa)$-projection onto constant functions defined by
\[
	( \pz v, w)_\kappa = (v,w)_\kappa,\qquad \forall w\in\real,
\]
that is explicitly given by $\pz v =\frac{1}{|\kappa|} \int_\kappa v \,dx$.
On $H^1(\kappa)$ we define the following scalar product and norm
\begin{align}
	\label{eq:ScalProd}
	\sscal{v}{w}
	&= (\pz v,\pz w)_\kappa + (\grad v, \grad w)_\kappa,
	\\
	\nonumber
	\snorm{v}
	&= \sscal{v}{v}^\frac12,
\end{align}
for all $v,w\in H^1(\kappa)$
and the corresponding projection $\pN:H^1(\kappa)\to \VN(\kappa)$  by 
\begin{equation}
	\label{eq:DefProj}
	\sscal{\pN v}{w_N}
	=
	\sscal{v}{w_N}
	\qquad \forall w_N\in\VN(\kappa).
\end{equation}
Then, it is easy to see that this projection satisfies the following properties
\begin{align}
	\nonumber
	( v - \pN v, c)_\kappa &= 0, &&\forall c\in\real,\forall v\in H^1(\kappa),
\end{align}
or equivalently expressed as $\pz (v-\pN v)=0$.
This implies that
\begin{align}
	\label{eq:ProjOrtho}
	(\grad (v-\pN v),\grad w_N)_\kappa &= 0, && \forall w_N\in\VN(\kappa),\forall v\in H^1(\kappa),\\
	\label{eq:ProjStability}
	\|\grad (v-\pN v)\|_\kappa &\le \|\grad v\|_\kappa, && \forall v\in H^1(\kappa),\\
	\nonumber
	\snorm{v-\pN v} &\le \snorm{v}, && \forall v\in H^1(\kappa).
\end{align}

\subsection{Local scaling constants}
In this section, we are going to define some local constants that will be used in the upcoming a posteriori error analysis.
We start with defining the local trace inverse inequality constant $\dkN$ for each $\kappa\in\cK$ defined by
\[
	\dkN \equiv \sup_{v_N\in\VN(\kappa)} \frac{\|\grad v_N\cd n_\kappa\|_{\partial\kappa}}{\snorm{v_N}}>0.
\]
Further, let
\[
	\rak 
	\equiv \sup_{\substack{v\in H^1(\kappa),\\ v\perp\VN(\kappa)}} \frac{\|v\|_{\kappa}}{\snorm{v}}
	\qquad\mbox{and}\qquad
	\rbk \equiv \sup_{\substack{v\in H^1(\kappa),\\ v\perp\VN(\kappa)}} \frac{\|v\|_{\partial\kappa}}{\snorm{v}},
\]
where $\perp$ is in the sense of the scalar product $\sscal{\cdot}{\cdot}$ defined by \eqref{eq:ScalProd}.
\begin{remark}[The computation of the constants $\rak$, $\rbk$ and $\dkN$]
	We provide more details in Section \ref{sec:computation} on how these local constants can be approximated by solving local eigenvalue problems.
\end{remark}

\begin{lem}
	\label{lem:LemA}
	Let $\kappa\in\cK$, $v\in H^1(\kappa)$.
	Then, there hods that
	\begin{align*}
			\|v-\Pi_N^\kappa v\|_\kappa 
			&\le 
			\rak\, \|\grad v \|_\kappa, \\
			\|v-\Pi_N^\kappa v\|_{\partial\kappa} 
			&\le 
			\rbk\, \|\grad v \|_\kappa.
	\end{align*}
\end{lem}
\begin{proof}
The proof consists of simply combining the definition of $\rak$ resp. $\rbk$ and the stability of the projection $\Pi_N^\kappa$ described in \eqref{eq:ProjStability}
\[
			\|v-\Pi_N^\kappa v\|_\kappa 
			\le 
			\rak\, \snorm{v-\Pi_N^\kappa v}
			=
			\rak\, \|\grad (v-\Pi_N^\kappa v) \|_\kappa
			\le 
			\rak\, \|\grad v\|_\kappa,
\]
since $\pz (v-\pN v)=0$. 
The proof for the second inequality is almost identical.
\end{proof}

%

\section{Poisson's equation}
\label{sec:Poisson}
As has been motivated in the introduction we start with a simple model problem that however reflects the difficulties associated to the discontinuous Galerkin method using non-polynomial functions.
The problem then reads: find $u\in \Hper\cap H^2(\cK)$ such that
\begin{equation}
	\label{eq:Poisson}
	-\Delta u = f,\quad\mbox{in } \Omega,
\end{equation}
for some $f\in L^2(\Omega)$.

\REV{
Given a piecewise constant and positive penalty function $\gamma$ such that $\gamma|_\kappa=\gk\in\real^+$ for all $\kappa\in\cK$, the discontinuous bilinear form is defined by
\begin{align*}
	a(w,v)
	&= \sum_{\kappa\in\cK}		
		\Big[ (\grad w, \grad v)_\kappa 
		- \tfrac12(\grad w,\jmp{v})_{\partial\kappa} 
		-  \tfrac{\theta}{2} (\jmp{w},\grad v)_{\partial\kappa} 
		+ \tfrac{\gk}{2} ( \jmp{w},\jmp{v})_{\partial\kappa}
		\Big],
\end{align*}
for any $w,v\in H^2(\cK)$ and for $\theta\in\real$.
Note that this is equivalent to the somewhat more standard notation
\[
	a(w,v)
	= (\grad w, \grad v)_\cK - (\avg{\grad w},\jmp{v})_\cF - \theta (\jmp{w},\avg{\grad v})_\cF + (\gF \jmp{w},\jmp{v})_\cF,
\]
with $\gF = \avg{\gamma}$ and where $(\cdot,\cdot)_\cF$ denotes the face-wise $L^2$-inner product over all faces of the mesh.}
Note that the choice of $\theta=1,-1$ corresponds to the symmetric and
non-symmetric interior penalty discontinuous Galerkin (SIPG
\REV{\cite{Wheeler78,Arnold1982}} or NIPG \REV{\cite{BaumannOden99}}) method, respectively.
The former case results in a symmetric bilinear form.

Then, the {\bf discontinuous Galerkin} approximation is defined by: Find $\uN\in\VN$ such that
\begin{equation}
	\label{eq:DG_Poisson}
	a(\uN,\vN) = (f,\vN)_\Omega,
	\qquad\forall \vN\in\VN.
\end{equation}

In this context we define the following {\bf broken energy norm} by
\REV{
\begin{equation}
	\label{eq:DGnorm}
	\tnorm{v}^2
	= \sum_{\kappa\in\cK}
		\Big[ \| \grad v \|_\kappa^2 + \tfrac{\gk}{2} \| \jmp{v} \|_{\partial\kappa}^2\Big],\quad \forall v\in H^1(\cK).
\end{equation}
}
Observe that 
\[
	 \tnorm{v}^2 = \sum_{\kappa\in\cK} \tnorm{v}^2_\kappa
	\qquad\mbox{with}\qquad 
	 \tnorm{v}^2_\kappa =\| \grad v \|_\kappa^2 + \tfrac{\gk}{2} \| \jmp{v} \|_{\partial\kappa}^2,
\]
and that this is indeed a norm as $\gamma>0$.
\REV{
As usual, the penalty parameter $\gamma$ needs to be chosen carefully to ensure coercivity.
}
Even when polynomial basis functions are used, the choice of an optimal $\gamma$ is not completely trivial and related discussions can be found in \cite{Epshteyn:2007vz,Ainsworth:2012kv}. 
The scaling in the element sizes and the polynomial order is however known \cite{Houston:2002uo,StammWihler2010}. 
The involved constants are resulting from applying trace and inverse inequalities, but no inverse inequality is known if general non-polynomial basis functions are employed.
To have a precise idea of the values of the combined trace and inverse inequalities for the generic non-polynomial basis functions spanning $\VN$, we propose here to use the local constants that were defined in Section \ref{sec:Prelim}.
In consequence, we can give a precise value of the piecewise constant function $\gamma$ that is needed to ensure coercivity of the bilinear form $a(\cdot,\cdot)$. 
This is stated in the following lemma.
\begin{lem}
	\REV{
	Under the assumption that $\tnorm{\cdot}$ is a norm (which is assumed here since $\gamma>0$) and if additionally $\gk \ge \frac12 \,(1+\theta)^2 \,(\dkN)^2$ for each $\kappa\in\cK$,
	} 
	then the bilinear form is coercive on $\VN$, i.e., there holds
	\[
		\tfrac12 \tnorm{v_N}^2 \le a(v_N,v_N), \qquad \forall v_N\in\VN.
	\]
\end{lem}
\begin{proof}
Since for any $v_N\in\VN$ we have $\grad v_N=\grad (v_N-\pz v_N)$ and $\snorm{v_N-\pz v_N}=\|\grad v_N\|_\kappa$
we can develop 
\begin{align*}
	a(v_N,v_N)
	&= \tfrac12\sum_{\kappa\in\cK}
		\Big[
			2 \, \|\grad v_N\|^2_\kappa 
			- (1+\theta)(\grad (v_N-\pz v_N),\jmp{v_N})_{\partial\kappa} 
			+ \gk \|\jmp{v_N}\|^2_{\partial\kappa}
		\Big]	\\
	&\ge 
		\tfrac12 \sum_{\kappa\in\cK}
		\Big[
			2 \, \|\grad v_N\|^2_\kappa 
			- (1+\theta) \, \dkN \snorm{v_N-\pz v_N}\|\jmp{v_N}\|_{\partial\kappa} 
			+ \gk \|\jmp{v_N}\|^2_{\partial\kappa}
		\Big]	\\
	&\ge 
		\tfrac12\sum_{\kappa\in\cK}
		\Big[
			\|\grad v_N\|^2_\kappa 
			+ \Big(\gamma_\kappa - \tfrac14(1+\theta)^2\, (\dkN)^2\Big)  \|\jmp{v_N}\|^2_{\partial\kappa}
		\Big]	
\end{align*}
and obtain
\[
	\tfrac12 \tnorm{v_N}^2 \le a(v_N,v_N)
\]
for the particular choice $\gamma_\kappa \ge \tfrac12(1+\theta)^2 \,(\dkN)^2$. Note however that for the particular choice of $\theta=-1$, $\gamma_\kappa$ stills needs to be positive in order that $\tnorm{\cdot}$ is indeed a norm.
\end{proof}


\subsection{Error representation}
Define the scaled error function $\varphi = \frac{u-\uN}{\tnorm{u-\uN}}$ and develop
\REV{
\begin{align*}
	\tnorm{u-\uN}
	&= \sum_{\kappa\in\cK}  
			\Big[ 
				( \grad(u-\uN),\grad\varphi)_\kappa
				+ \tfrac{\gk}{2} ( \jmp{u-\uN},\jmp{\varphi})_{\partial\kappa}
			\Big]\\
	&= a(u-\uN,\varphi) + \tfrac{1+\theta}{2} \sum_{\kappa\in\cK}    (\grad \varphi,\jmp{u-\uN})_{\partial\kappa}.
\end{align*}
}
We prefer to work with the scaled error function $\varphi$ for sake of a simple presentation of the upcoming error analysis.
Observe that due to the regularity of $u\in\Hper$, which implies $\jmp{u}=0$, and since $u$ is indeed the solution of \eqref{eq:Poisson}, there holds
\REV{
\[
	a(u,\varphi) 
	= 
	\sum_{\kappa\in\cK}		
		\Big[ (\grad u, \grad \varphi)_\kappa 
		- \tfrac12(\grad u,\jmp{\varphi})_{\partial\kappa} 
		\Big]
	= - (\Delta u, \varphi)_\Omega = (f,\varphi)_\Omega.
\]
}
On the other hand, since $u_N\in\VN$ is the DG-solution solution of \eqref{eq:DG_Poisson}, we obtain
\[
	-a(\uN,\varphi) 
	= -a(\uN,\varphi-\varphi_N) - (f,\varphi_N)_\Omega,
\]
for any $\varphi_N\in\VN$.
Thus, using the integration by parts and Lemma \ref{lem:IntByParts}, we can develop
\REV{
\begin{align*}
	-a(\uN,\varphi)
	&=  
		\sum_{\kappa\in\cK}		
		\Big[ 
		- (\grad \uN, \grad (\varphi-\varphi_N))_\kappa 
		+ \tfrac12(\grad \uN,\jmp{\varphi-\varphi_N})_{\partial\kappa} 
		+  \tfrac{\theta}{2} (\jmp{\uN},\grad (\varphi-\varphi_N))_{\partial\kappa} 
		\\
		&\quad \qquad\quad
		- \tfrac{\gk}{2} ( \jmp{\uN},\jmp{\varphi-\varphi_N})_{\partial\kappa}
		\Big]
		- (f,\varphi_N)_\Omega\\
	&= 
	    \sum_{\kappa\in\cK}		
		\Big[
		 (\Delta \uN, \varphi-\varphi_N)_\kappa 
		- \tfrac12(\jmp{\grad \uN},\varphi-\varphi_N)_{\partial\kappa} 
		+  \tfrac{\theta}{2} (\jmp{\uN},\grad (\varphi-\varphi_N))_{\partial\kappa} 
		\\
		&\quad \qquad\quad
		- \tfrac{\gk}{2} ( \jmp{\uN},\jmp{\varphi-\varphi_N})_{\partial\kappa}
		\Big]
		- (f,\varphi_N)_\Omega,
\end{align*}
and obtain the {\bf error representation equation}
\begin{align}
	\tnorm{u-\uN}
	&= 
	\label{eq:ErrRepPoisson}	
	\sum_{\kappa\in\cK}
	\Big[
	(f+\Delta \uN, \varphi\sm\varphi_N)_\kappa
	- \tfrac12(\jmp{\grad \uN},\varphi\sm\varphi_N)_{\partial\kappa} 
	\\
	\nonumber
	&\qquad\qquad
	- \gk(\jmp{\uN},(\varphi\sm\varphi_N)n_\kappa)_{\partial\kappa} 
	- \tfrac12(\jmp{\uN},\grad \varphi\sp\theta\grad \varphi_N)_{\partial\kappa} 
	\Big].
\end{align}
}

\subsection{A posteriori error estimation}
After recalling that we assumed that $u\in H^2(\kappa)$, we start by introducing the constant $\dku(u_N)$ defined by
\begin{align*}
	\dku(u_N) = \frac{\|\grad(u-u_N)\cd n_\kappa\|_{\partial\kappa}}{\|\grad(u-u_N)\|_{\kappa}},
\end{align*}
and define the constant $\rck$ by
\[
	\rck = \dku(u_N)+\dkN|\theta|.
\]
We note that in practice, the constant $\dku(u_N)$ can not be evaluated since $u$ is unknown.
See in the upcoming numerical examples how we deal with this term.
\begin{remark}
Observe that $\dku(u_N)$ is bounded by the constant 
\begin{align*}
	\sup_{\substack{v_N\in\VN(\kappa)}} \frac{\|\grad(u-v_N)\cd n_\kappa\|_{\partial\kappa}}{\|\grad(u-v_N)\|_{\kappa}}<\infty,
\end{align*}
which, in turn, is independent of the approximation $u_N$ (but still depends on the exact solution $u$ and the approximation space $\VN$).
\end{remark}

Define the following estimators
\begin{align}
	\label{eq:Estim1}
	\etaR &\equiv \rak\|f+\Delta \uN\|_\kappa,\\
	\label{eq:Estim2}
	\etaF &\equiv \tfrac{\rbk}{2} \|\jmp{\grad \uN}\|_{\partial\kappa},\\
	\label{eq:Estim3}
	\etaJ &\equiv (\rbk\,\gk+\tfrac{\rck}{2})\|\jmp{\uN}\|_{\partial \kappa},
\end{align}
in order to state the first Theorem.

\begin{thm}
	\label{thm:PoissonApost}
	Let $u\in \Hper\cap H^2(\cK)$ be the solution of \eqref{eq:Poisson} and $u_N\in\VN$ the DG-approximation defined by \eqref{eq:DG_Poisson}. Then, we have the following a posteriori upper bound
	\[
	\tnorm{u-\uN}
  \le 		
  \left(
			\sum_{\kappa\in\cK}
			\Big[ \etaR+\etaF+\etaJ \Big]^2
	\right)^\frac12.
	\]
\end{thm}
\begin{proof}

Using the triangle inequality, observe that
\[
	\|(\grad \varphi+\theta\grad \varphi_N)\cd n_\kappa\|_{\partial\kappa}
	\le
	\|\grad \varphi\cd n_\kappa\|_{\partial\kappa}
	+|\theta|\|\grad \varphi_N\cd n_\kappa\|_{\partial\kappa}
	\le 
	\dku(u_N)\|\grad \varphi \|_\kappa
	+\dkN|\theta|\|\grad \varphi_N\|_{\kappa}.
\]
So far, the results where valid for any arbitrary discrete function $\varphi_N\in\VN$.
In this proof we consider the particular choice $\varphi_N|_\kappa =\Pi_N^\kappa \varphi$ so that we can easily state
\[
	\|\grad \varphi_N\|_{\kappa} 
	\le \|\grad \varphi\|_{\kappa}
\]
by splitting $\varphi = \Pi_N^\kappa \varphi + (\varphi-\Pi_N^\kappa \varphi)$ and using the orthogonality relation \eqref{eq:ProjOrtho}.
Then, there holds
\begin{equation}
	\label{eq:DifficultTerm}
	\|(\grad \varphi+\theta\grad \varphi_N)\cd n_\kappa\|_{\partial\kappa}
	\le 
	\underbrace{\left(\dku(u_N)+\dkN|\theta|\right)}_{=\rck} \|\grad \varphi\|_\kappa
	\le 
	\rck \|\grad \varphi\|_\kappa,
\end{equation}
by applying a simple triangle inequality.

\REV{
If we apply the Cauchy-Schwarz inequality to the error representation
formula \eqref{eq:ErrRepPoisson}, in combination with Lemma
\ref{lem:LemA}, Eq.~\eqref{eq:DifficultTerm} and another
Cauchy-Schwarz inequality, we have (recall that $\varphi =
\frac{u-\uN}{\tnorm{u-\uN}}$)}
\begin{align*}
	\tnorm{u-\uN}
	&\le
	\sum_{\kappa\in\cK}
	\Big[
		\|f+\Delta \uN\|_\kappa \, \|\varphi-\varphi_N\|_\kappa
		+
		\tfrac12\|\jmp{\grad \uN}\|_{\partial\kappa} \, \|\varphi-\varphi_N\|_{\partial\kappa} 
		+ 
		\gk \|\jmp{\uN}\|_{\partial \kappa} \, \| \varphi-\varphi_N \|_{\partial\kappa} 		\\
		&\qquad\qquad
		+\tfrac12 \|\jmp{\uN}\|_{\partial \kappa} \, \|(\grad \varphi+\theta\grad \varphi_N)\cd n_\kappa\|_{\partial\kappa}
	\Big]
	\\
	&\le
	\sum_{\kappa\in\cK}
	\Big[
		\rak \, \|f+\Delta \uN\|_\kappa
		+
		\tfrac{\rbk}{2}\|\jmp{\grad \uN}\|_{\partial\kappa}
		+ 
		(\gk \rbk+\tfrac{\rck}{2}) \|\jmp{\uN}\|_{\partial \kappa}	
	\Big] 
	\, \|\grad \varphi \|_\kappa
	\\
	&\le
	\left(
	\sum_{\kappa\in\cK}
	\Big[
		\rak \, \|f+\Delta \uN\|_\kappa
		+
		\tfrac{\rbk}{2}\|\jmp{\grad \uN}\|_{\partial\kappa}
		+ 
		(\gk \rbk+\tfrac{\rck}{2}) \|\jmp{\uN}\|_{\partial \kappa}	
	\Big]^2\right)^{\frac12}\\
	&=
	  \left(
			\sum_{\kappa\in\cK}
			\Big[ \etaR+\etaF+\etaJ \Big]^2
	\right)^\frac12.
\end{align*}
\end{proof}

\section{Second order indefinite problems}
\label{sec:Indefinite}
In this section we consider the more general indefinite equation: find $u\in \Hper$ such that
\begin{equation}
	\label{eq:Indef}
	-\Delta u + V u = f,\qquad\mbox{in }\Omega, 
\end{equation}
for some $f\in L^2(\Omega)$ and where we only assume that $V\in L^\infty(\Omega)$ is bounded and that the operator $-\Delta +V$ has no zero eigenvalue. 
For the particular choice of $V=-k^2\in \real$ this framework includes the Helmholtz equation.
The DG-bilinear form is provided by
\REV{
\begin{align*}
	a(w,v)
	&= \sum_{\kappa\in\cK}		
	\Big[
		(\grad w, \grad v)_\kappa 
		+ (Vw,v)_\kappa
		- \tfrac12(\grad w,\jmp{v})_{\partial\kappa}
		- \tfrac{\theta}{2} (\jmp{w},\grad v)_{\partial\kappa} 
		+ \tfrac{\gk}{2}( \jmp{w},\jmp{v})_{\partial\kappa}
	\Big],
\end{align*}
}
such that the DG-approximation is defined by: Find $\uN\in\VN$ such that
\begin{equation}
	\label{eq:DG_Indef}
	a(\uN,\vN) = (f,\vN)_\Omega,
	\qquad\forall \vN\in\VN,
\end{equation}
and we keep the definition of the broken energy norm of \eqref{eq:DGnorm}.
Of course the choice $\gk = 2\,(1+|\theta|)^2 \,(\dkN)^2$ does not imply coercivity of the bilinear form in this setting any more. 
We assume that $\gk$ has been chosen by the user to insure that DG-problem has a unique solution and focus on how to quantify the error a posteriori. 
Observe that whenever the DG-problem is not uniquely solvable, the solver of the numerical system typically reveals the lack of well-posedness. 
The following analysis requires that the DG-solution satisfies \eqref{eq:DG_Indef}.

\subsection{Computable upper bounds}

We first introduce a modified norm. For this consider $\Vp$ and $\Vm$ defined by $\Vp=\max(V,0)\ge 0$ and $\Vm = \max(-V,0)\ge 0$ so that $V=\Vp-\Vm$ and $|V| = \Vp+\Vm$.
Then, define
\[
	 \tnorm{v}^2 
	 = \sum_{\kappa\in\cK} \tnorm{v}^2_\kappa
	\qquad\mbox{with}\qquad 
	 \tnorm{v}^2_\kappa 
	 =\| \grad v \|_\kappa^2 
 			+\| \Vp^\frac12  v \|_\kappa^2 
			+ \tfrac{\gk}{2} \| \jmp{v} \|_{\partial\kappa}^2,
	\qquad \forall v\in H^1(\cK).
\]

Applying similar arguments as in Section \ref{sec:Poisson} the following error representation can be developed
\begin{align}
	\tnorm{u-\uN}
	&= 
	\label{eq:ErrRepIndef}	
	\sum_{\kappa\in\cK}
	\Big[
	(f + \Delta \uN - V\uN, \varphi\sm\varphi_N)_\kappa
  +(\Vm(u-\uN),\varphi)_\kappa \Big]
	\\
	\nonumber
	&\quad
	-
	\tfrac12
	\sum_{\kappa\in\cK}
	\Big[
	(\jmp{\grad \uN},\varphi\sm\varphi_N)_{\partial\kappa} 
	+ \gk(\jmp{\uN},\jmp{\varphi\sm\varphi_N})_{\partial\kappa} 
	+ (\jmp{\uN},\grad \varphi\sp\theta\grad \varphi_N)_{\partial\kappa} 
	\Big].
\end{align}
Redefining the residual as 
\begin{equation}
	\label{eq:Redef_etaR}
	\etaR \equiv \rak\|f+\Delta \uN - V \uN\|_\kappa,
\end{equation}  
the following bound can be developed.
\begin{thm}
	\label{thm:IndefApost}
	Let $u\in \Hper\cap H^2(\cK)$ be the solution of \eqref{eq:Indef} and $u_N\in\VN$ the DG-approximation defined by \eqref{eq:DG_Indef}. Then, we have the following a posteriori upper bound
\[
	\tnorm{u-\uN}
 	\le 		
  	\left(
			\sum_{\kappa\in\cK}
			\Big[\etaR+\etaF+\etaJ\Big]^2
	\right)^\frac12
	+ \frac{\|\Vm^\frac12 (u-\uN)\|_\cK^2}{\tnorm{u-\uN}},
	\]
	where $\etaR$ is defined by \eqref{eq:Redef_etaR} and $\etaF,\etaJ$ are defined by \eqref{eq:Estim2}--\eqref{eq:Estim3}.
\end{thm}
\begin{proof}
	This estimate can be obtained by applying the Cauchy-Schwarz inequality to the error representation equation \eqref{eq:ErrRepIndef} similar as in the proof of Theorem \ref{thm:PoissonApost}. 
	Only the additional term 
	\[
		(\Vm(u-\uN),\varphi)_\cK
		= \frac{\|\Vm^\frac12 (u-\uN)\|_\cK^2}{\tnorm{u-\uN}}
	\]
	is not estimated.
\end{proof}

\begin{remark}
  For $\Vm\in L^{\infty}$, the term  $\frac{\|\Vm^\frac12
  (u-\uN)\|_\cK^2}{\tnorm{u-\uN}}$ is small compared to
  the upper bound estimator in the limit of complete basis set.
  On the other hand, when only a small number of basis functions are
  used, this term can become large, and the upper bound error estimator
  can underestimate the true error in energy norm.
\end{remark}

\subsection{Computable lower bounds}
The goal of this section is to derive computable lower bounds of the approximation error. We note that the following theory applies also to the Poisson problem with $V=0$.

Observe that 
\[
	\etaJ
	= \left(\rbk\,\gk+\tfrac{\rck}{2}\right)
	\|\jmp{\uN}\|_{\partial\kappa}
	\le  
	\sqrt{\tfrac{2}{\gk}}	\left(\rbk\,\gk+\tfrac{\rck}{2}\right) \,\tnorm{u-u_N}_\kappa.
\]
Second, for any face $F$ of $\partial\kappa$, denote by $\kappa'$ the adjacent element such that $F=\partial \kappa \cap \partial\kappa'$ such that there holds
\begin{align}
	\label{eq:loceq}
	\etaF^2
	&= \tfrac{\rbk^2}{4} \|\jmp{\grad \uN}\|_{\partial\kappa}^2
	= \tfrac{\rbk^2}{4} \|\jmp{\grad (u-\uN)}\|_{\partial\kappa}^2
	\le \tfrac{\rbk^2}{2} 
		\sum_{F\in\partial\kappa}		
		\Big(
				\|\grad (u-\uN)|_{\kappa}\cdot n_\kappa\|_{F}^2
			+ \|\grad (u-\uN)|_{\kappa'}\cdot n_{\kappa'}\|_{F}^2
		\Big).
\end{align}
Recall that $\omega(\kappa)$ is the patch consisting of $\kappa$ and its adjacent elements sharing one face, then
\[
	\etaF^2
	\le \tfrac{\rbk^2}{2} 
			\sum_{\kappa'\in\omega(\kappa)}  
				\big(
					\dkpu(u_N)
					\|\grad (u-\uN)\|_{\kappa'}
				\big)^2
    \le 
    	\tfrac{\rbk^2}{2} 
    	\Big(
		\max_{\kappa'\in\omega(\kappa)}\dkpu(u_N)
		\Big)^2
			\sum_{\kappa'\in\omega(\kappa)}  
				\|\grad (u-\uN)\|_{\kappa'}^2
\]
Further, let \REV{$\bub_\kappa$} be a smooth non-negative bubble
function with \REV{$\sup_{x\in\kappa}\bub_\kappa(x)=1$} and local support, i.e. \REV{$\mbox{supp}(\bub_\kappa)\subset \kappa$}, which in turn implies that \REV{$\bub_\kappa|_{\partial\kappa} = 0$}.
Finally, let us denote the residual by $R=f+\Delta u_N-Vu_N$ and define
\REV{
\[
	\sigma_\kappa = \rak\frac{\|R\|_\kappa}{\|\bub_\kappa^\frac12R\|_\kappa^2}.
\]
Denote by $\varphi_\kappa\in H^1_0(\kappa)$ the solution to 
\[
	-\Delta\varphi_\kappa = V\bub_\kappa R,\qquad\mbox{on } \kappa,
\]
so that 
\begin{align*}
	\etaR 
	&= \rak\|R\|_\kappa 
	= \sigma_\kappa\|\bub_\kappa^\frac12R\|_\kappa^2
	=  \sigma_\kappa
		\int_\kappa \bub_\kappa \Big[-\Delta(u-u_N)+ V(u-\uN)\,\Big]\, R
	\\
		&=-  \sigma_\kappa
		\int_\kappa \Big[\Delta(u-u_N)\,\bub_\kappa \,R - \Delta\varphi_\kappa(u-u_N) \,\Big]
	\\
	&= \sigma_\kappa
			\int_\kappa 
			\Big[\nabla(u-u_N) \cdot \nabla (\bub_\kappa \, R)
			-
			\nabla(u-\uN)\cdot \nabla \varphi_\kappa \Big]
	\\ 
	&\le 
	 \sigma_\kappa
			\|\nabla(u-u_N)\|_\kappa
			 \|\nabla (\bub_\kappa \, R - \varphi_\kappa)\|_\kappa,
\end{align*}
and in consequence
\[
	\frac{\etaR}{\tnorm{u-u_N}_\kappa}
	\le 
	 \sigma_\kappa
			 \|\nabla (\bub_\kappa \, R - \varphi_\kappa)\|_\kappa.
\]}
The results above indicate that
\begin{equation}
	\label{eq:locestim}
  \tnorm{u-u_{N}}_{\kappa} 
  \ge 
  \max\left\{
  	\frac{\etaR}{\ccR},\frac{\etaJ}{\ccJ}
  \right\},
  \qquad
 \tnorm{u-u_{N}}_{\omega(\kappa)} \ge 
  \frac{\etaF}{\ccF},
\end{equation}
where,  denoting by $|\omega(\kappa)|$ the cardinality of the set $\omega(\kappa)$, we use the definitions
\[
	 \tnorm{v}_{\omega(\kappa)}^2 
	 = \frac{1}{|\omega(\kappa)|}\sum_{\kappa'\in\omega(\kappa)}  
				\|\grad v\|_{\kappa'}^2
		+ \tfrac{\gk}{2} \| \jmp{v}\|^2_{\partial\kappa},
\]
and
\begin{align*}
    \ccR &= \rak \frac{\|R\|_{\kappa} \|\grad (b_\kappa
    \, R - \varphi_{\kappa})\|_\kappa}{\| b_\kappa^{1/2} \, R \|^2_{\kappa}},
    \\
   	\ccF &= \rbk \sqrt{\tfrac{|\omega(\kappa)|}{2}}
   	\max_{\kappa'\in\omega(\kappa)}\dkpu(u_N),
   	\\
    \ccJ &= \sqrt{\tfrac{2}{\gk}}	\left(\rbk\,\gk+\tfrac{\rck}{2}\right).
\end{align*}
%
We summarize the results in the following proposition.
\begin{prop}[Local lower bound]
	Let $u\in \Hper\cap H^2(\cK)$ be the solution of \eqref{eq:Indef} and $u_N\in\VN$ the DG-approximation defined by \eqref{eq:DG_Indef}. 
	Then, the quantity
\[
	\xi_{\kappa}
	=	
	\max\left\{ 
	\frac{\etaR}{\ccR},
	\frac{\etaF}{\ccF},
	\frac{\etaJ}{\ccJ}
  \right\},
\]
is a local lower bound of the local error 
	\[
		\max\left\{\tnorm{u-u_{N}}_{\kappa},\tnorm{u-u_{N}}_{\omega(\kappa)}\right\}.
	\]
\end{prop}

\begin{remark}
Since in practice, the nominator as well as the denominator of any of those fractions might become very small, these ratios are not numerically stable. It turns out that 
\[
	\xi_{\kappa} = \frac{\etaR+\etaF+\etaJ}{\ccR+\ccF+\ccJ}
\]
is numerically more robust and still meaningful as it replaces the maximum by the average.
\end{remark}

On a global level, the following result holds.
\begin{prop}[Global lower bound]
	\label{prop:GlobLowerBound}
	Let $u\in \Hper\cap H^2(\cK)$ be the solution of \eqref{eq:Indef} and $u_N\in\VN$ the DG-approximation defined by \eqref{eq:DG_Indef}. 
	Then, there holds that
	\[
		\xi 
		= \frac{
				\left(
			\sum_{\kappa\in\cK} 
			\Big[ \etaR + \etaF + \etaJ \Big]^2
		\right)^\frac12
		}{
		\sqrt{3}
		\max_{\kappa\in\cK}
			\left(
			\ccR^2 
			+
			{\rbok^2} \dku(\uN)^2
			+ 
			\ccJ^2
		\right)^\frac12}
		\le 
		\tnorm{u- \uN},
	\]
	where 
	\[
		\rbok^2 
		= \max_{F\in\partial\kappa}\avg{\rbk^2}|_F
		= \max_{F\in\partial\kappa}
			\left(\tfrac{\rbk^2}{2} + \tfrac{\rbkp^2}{2} \middle)\right|_F.
	\]
\end{prop}
\begin{proof}
Observe that by \eqref{eq:loceq} there holds
\begin{align*}
	\sum_{\kappa\in\cK} \etaF^2
	&\le 
	\sum_{\kappa\in\cK} \tfrac{\rbk^2}{2} 
		\sum_{F\in\partial\kappa}
		\|\grad(u- \uN)|_{\kappa}\cdot n_\kappa\|_{F}^2
	+
	\sum_{\kappa\in\cK} \tfrac{\rbk^2}{2} 
		\sum_{F\in\partial\kappa}
		\|\grad(u- \uN)|_{\kappa'}\cdot n_{\kappa'}\|_{F}^2
	\\
	& =  
	\sum_{\kappa\in\cK} \tfrac{\rbk^2}{2} 
		\sum_{F\in\partial\kappa}
		\|\grad(u- \uN)|_{\kappa}\cdot n_\kappa\|_{F}^2
	+
	\sum_{\kappa\in\cK} 
		\sum_{F\in\partial\kappa}\tfrac{\rbkp^2}{2} 
		\|\grad(u- \uN)|_{\kappa}\cdot n_{\kappa}\|_{F}^2
	\\
	& =  
	\sum_{\kappa\in\cK} 
		\sum_{F\in\partial\kappa}
		\left(\tfrac{\rbk^2}{2} + \tfrac{\rbkp^2}{2} \right)
		\|\grad(u- \uN)|_{\kappa}\cdot n_\kappa\|_{F}^2
	 =  
	\sum_{\kappa\in\cK} 
		\sum_{F\in\partial\kappa}
		\avg{\rbk^2}
		\|\grad(u- \uN)|_{\kappa}\cdot n_\kappa\|_{F}^2
	\\
	& 	 \le
	\sum_{\kappa\in\cK} {\rbok^2}
		\|\grad(u- \uN)|_{\kappa}\cdot n_\kappa\|_{\partial\kappa}^2
	\le
	\sum_{\kappa\in\cK} {\rbok^2} \dku(\uN)^2
		\|\grad(u- \uN)\|_{\kappa}^2.
\end{align*}
Then, using the other local estimates for $\etaR$ and $\etaJ$ given by \eqref{eq:locestim} yields
\begin{align*}
	\sum_{\kappa\in\cK} 
	\Big[ \etaR + \etaF + \etaJ \Big]^2
	&\le 
	3\sum_{\kappa\in\cK} 
	\left(
	\etaR^2 + \etaF^2 + \etaJ^2
	\right)
	\le 3
	\sum_{\kappa\in\cK} 
	\left(
		\ccR^2 
		+
		{\rbok^2} \dku(\uN)^2
		+ 
		\ccJ^2
	\right)
	\tnorm{u- \uN}_\kappa^2\\
	&
	\le 3
	\max_{\kappa\in\cK}
		\left(
		\ccR^2 
		+
		{\rbok^2} \dku(\uN)^2
		+ 
		\ccJ^2
	\right)
	\tnorm{u- \uN}^2.
\end{align*}
\end{proof}

%

\newcommand{\wt}[1]{\widetilde{#1}}
\newcommand{\bvec}[1]{\mathbf{#1}}
\newcommand{\mI}{\mathcal{I}}
\newcommand{\mJ}{\mathcal{J}}

\section{Practical strategies for estimating the
constants}\label{sec:computation}

In this section we discuss how to compute the constants $\dkN,\rak,\rbk$ as
defined in Section~\ref{sec:Prelim} in the a posteriori error
estimator for general non-polynomial basis functions in the discontinuous Galerkin
framework. The basic strategy is to discretize the infinite dimensional
representative space $H^{1}(\kappa)$ using a finite dimensional space such as high order
polynomials, and to replace the various inner products defined in
Section~\ref{sec:Prelim} by discrete bilinear forms using Gauss
quadrature.  With the help of these bilinear forms, $\dkN,\rak,\rbk$ can be
estimated by solving an eigenvalue problem, locally on each element
$\kappa$.

\subsection{Finite dimensional discretization}

For simplicity let $\kappa=[0,h]^d,d=1,2,3$ and all quantities be real.  We start the discussion with
$d=1$, i.e. $\kappa=[0,h]$.  All  numerical quadrature are to be performed
using the Legendre-Gauss-Lobatto (LGL) quadrature with $N_{g}$ points.
The LGL grid points are denoted by $\{y_{j}\}_{j=1}^{N_{g}}$, and the
corresponding LGL weights by $\{\omega_{j}\}_{j=1}^{N_{g}}$.
The Lobatto quadrature implies that 
\[
y_{1}=0, \quad y_{N_{g}}=h,
\] 
which facilitates the description of the boundary
integrals as in the estimate of $\dkN$ and $\rbk$.  The LGL grid points
$\{y_{j}\}_{j=1}^{N_{g}}$ correspond to a unique set of Lagrange
polynomials of degree $(N_{g}-1)$, denoted by
$\{p_{j}(x)\}_{j=1}^{N_{g}}$, and satisfy
\[
p_{j}(y_{i}) = \delta_{ij}, \quad 1\le i,j\le N_{g},
\]
where $\delta_{ij}$ is the Kronecker $\delta$ function.  
We can then approximate $v\in H^1(\kappa)$ using the linear combination of Lagrange polynomials
as
\[
  v(x) \approx \sum_{j=1}^{N_{g}} v_j \, p_{j}(x).
\]
\REV{The} sequence of  spaces $\mathbb P_{N_g}$ of polynomials of degree $N_g$ being dense in  $H^1(\kappa)$ implies that, for any $v\in H^1(\kappa)$ and any $\varepsilon>0$, there exists $N_g$ and $v^1_{N_g},v_{N_g}^2\in\mathbb P_{N_g}$ such that
\[
	\frac{ \| v - v^1_{N_g} \|_{\kappa} }{\| v \|_{\kappa}} \le \varepsilon
	\qquad\mbox{and similarly}\qquad
	\frac{ \snorm{v - v^2_{N_g} } }{\snorm{v}} \le \varepsilon,
\]
if choosing $N_g$ large enough.
That is, elements in $H^1(\kappa)$ can be approximated, in the sense of $L^2$ and $H^1$ with any desired accuracy by elements of $\mathbb P_{N_g}$.
This motivates us to work in $\mathbb P_{N_g}$ instead of $H^1(\kappa)$ for $N_g$ large enough.
We assume that $N_{g}$ is large enough so that the above approximation error in the local $L^2$ and $H^1$-norms are very small.
Further, for functions $u,v\in\mathbb P_{N_g}$, the LGL quadrature for computing the
inner product $(u,v)_{\kappa}$ converges rapidly with respect to the
increase of $N_{g}$.  

We denote by
$v=(v_{1},\ldots,v_{N_{g}})^T$ the column vector corresponding to the
coefficients of $v\in\mathbb P_{N_g}$, and denote by 
$Y=(y_{1},\ldots,y_{N_{g}})^{T}$, $w=(\omega_{1},\ldots,\omega_{N_{g}})^{T}$ the column vector
corresponding to the LGL grid points and weights, respectively.  With a slight abuse of notation we can compute the inner product
using linear algebra notation as
\begin{equation}
  (u,v)_{\kappa} = \sum_{j=1}^{N_{g}} u_{j} \omega_{j} v_{j} \equiv u^T W v,
  \label{eqn:innerproduct1D}
\end{equation}
where $W=\text{diag}[w]$ is a diagonal matrix with the entries of vector $w$ on the
diagonal entries.  

The Lagrange polynomials also induce a \textit{differentiation matrix}
$D$ of size $N_{g}\times N_{g}$, defined as
\begin{equation}
  D_{ij}=p'_{j}(y_{i}), \quad 1\le i,j \le N_{g}.
  \label{eqn:DifferentiationMatrix1D}
\end{equation}
Taking the derivative of a polynomial yields
\[
  v'(x) = \sum_{j=1}^{N_{g}} p'_{j}(x) v_{j}. 
\]
Let $v'=(v'(y_{1}),\ldots,v'(y_{N_{g}}))^{T}$ be the column vector of
the derivative quantity $v'(x)$ on the LGL grid points, then
\begin{equation}
  v' = D v.
  \label{eqn:Differentiation1Dapply}
\end{equation}
Eq.~\eqref{eqn:Differentiation1Dapply} shows that the differentiation
matrix maps the values of a function to the values of its derivative on
the LGL grid points. Using the differentiation matrix, inner products of
the form $(u',v')_{\kappa}$ can be expressed in linear algebra notation as
\begin{equation}
  (u',v')_{\kappa} = (Du)^{T} W (Dv) = u^{T} (D^{T}WD) v. 
  \label{eqn:Gradient1Dinnerproduct}
\end{equation}
In order to compute the inner product $(u,v)_{\star,\kappa}$ we also need to
compute $(\pz u, \pz v)_{\kappa}$.  Note
that 
\[
\pz v = \frac{1}{|\kappa|} (1,v)_{\kappa} = \frac{1}{|\kappa|} w^{T} v,
\]
with $|\kappa|=h$.  Then
\[
(\pz u, \pz v)_{\kappa} = \frac{1}{|\kappa|^2} u^{T}w w^{T} v |\kappa|
= u^{T} \left(w \frac{1}{|\kappa|} w^{T}\right) v.
\]
Therefore the inner product $(u,v)_{\star,\kappa}$ can be computed as
\begin{equation}
  (u,v)_{\star,\kappa} = u^{T} \left(D^{T}WD + w \frac{1}{|\kappa|}
  w^{T}\right) v.
  \label{eqn:star1Dinnerproduct}
\end{equation}

We also need to compute inner products on the boundary $\partial\kappa$.  In
1D, $v\vert_{\partial\kappa}(x)$ is completely described by two points
$v(0)$ and $v(h)$, which is given by the discretization on the LGL grid
points as $v_{1}$ and $v_{N_{g}}$.  Define the weight vector at
$0$-dimension as $\wt{w}=(1,0,\ldots,0,1)^{T}$, and $\wt{W}=\text{diag}[\wt{w}]$, then
the inner product on the boundary can be expressed as
\begin{equation}
  (u,v)_{\partial\kappa} = u_{1} v_{1} + u_{N_{g}} v_{N_{g}} \equiv u^T
  \wt{W} v.
  \label{eqn:innerproductBoundary1D}
\end{equation}
Similarly
\begin{equation}
  (u',v')_{\partial\kappa} = u'_{1} v'_{1} + u'_{N_{g}} v'_{N_{g}} \equiv
  u^{T} D^{T} \wt{W} D v.
  \label{eqn:innerproductGradientBoundary1D}
\end{equation}
The inner
products~\eqref{eqn:innerproduct1D},~\eqref{eqn:Gradient1Dinnerproduct},~\eqref{eqn:star1Dinnerproduct}
and~\eqref{eqn:innerproductGradientBoundary1D} are sufficient for
estimating $\dkN,\rak,\rbk$ for $d=1$.

Now we generalize all the definition above to $d>1$.  Though in practice
we only consider $d=2,3$, the formalism developed here holds for any dimension.
For any $x\in \kappa=[0,h]^{d}$, we denote by $x=(x^{(1)},\ldots,x^{(d)})^{T}$, with
$x^{(l)}$ being the component of $x$ along the $l$-th dimension. Then
the set of $N_{g}^{d}$ LGL grid points in the dimension $d$ is given by
\begin{equation}
  Y^{[d]} = \{y_{j_{1},\ldots,j_{d}}\equiv (y_{j_1},\ldots,y_{j_d})^{T}| 1\le j_1,\ldots,j_{d}\le N_{g}\}.
  \label{eqn:gridDD}
\end{equation}

We define the tensor product of $d$ matrices  $A^{(1)},\ldots,A^{(d)}$ 
of size $N_{g}\times N_{g}$ as
\begin{equation}
  A_{i_{1}j_{1},\ldots,i_{d}j_{d}} = \prod_{l=1}^{d}
  A^{(l)}_{i_{l}j_{l}}, \quad
  1\le i_{1},j_{1},\ldots,i_{d},j_{d}\le N_{g},
  \label{eqn:matTensor}
\end{equation} 
which can be written in a compact form as
\begin{equation}
 A \equiv \bigotimes_{l=1}^{d} A^{(l)}. 
  \label{eqn:tensorCompact1}
\end{equation}
From the computational point of view, it is more convenient to rewrite
the tensor product $A$ as a matrix by stacking the $i_{1},\ldots,i_{d}$ and
$j_{1},\ldots,j_{d}$ indices, respectively.  In other words, we can
view $A$ as a large matrix  of size
$N_{g}^{d}\times N_{g}^{d}$, and each matrix element
$A_{i_{1}j_{1},\ldots,i_{d}j_{d}}$ corresponds to a matrix element
$A_{\mI\mJ}$, with the index
\[
\mI=1+\sum_{l=1}^{d} (i_{l}-1) N_{g}^{(l-1)},\quad
\mJ=1+\sum_{l=1}^{d} (j_{l}-1) N_{g}^{(l-1)}.
\]
Note that when $d=2$, the stacked representation of the tensor product
of $A^{(1)}$ and $A^{(2)}$ is the Kronecker product of $A^{(2)}$ and
$A^{(1)}$.  

We also define a special case for the tensor product of
$d$ vectors $v^{(1)},\ldots,v^{(d)}$ of
size $N_{g}$.  By viewing each $v^{(l)}$ as a matrix of size
$N_{g}\times 1$, we have
\begin{equation}
  v_{j_{1},\ldots,j_{d}} = \prod_{l=1}^{d} v^{(l)}_{j_{l}},\quad
  1\le j_{1},\ldots,j_{d}\le N_{g}.
  \label{eqn:vecTensor}
\end{equation}
Eq.~\eqref{eqn:vecTensor} can be written in a
compact form as
\begin{equation}
  v\equiv \bigotimes_{l=1}^{d} v^{(l)}.
  \label{eqn:tensorCompact2}
\end{equation}
By stacking the indices $j_{1},\ldots,j_{d}$ together, we can view
$v$ as a vector of size $N_{g}^{d}$, and each element
$v_{j_{1},\ldots,j_{d}}$ corresponds to an element $v_{\mJ}$ with
$\mJ=1+\sum_{l=1}^{d} (j_{l}-1) N_{g}^{(l-1)}$.
Using the notation of tensor product, the set of LGL weights is
described by a vector
\begin{equation}
  w^{[d]} = \bigotimes_{l=1}^{d} w.
  \label{eqn:weightDD}
\end{equation}

Similar to the 1D case, each LGL grid point
$y_{j_{1},\ldots,j_{d}}$ uniquely corresponds to a Lagrange polynomial
\[
p_{j_{1},\ldots,j_{d}}(x) = \prod_{l=1}^{d}
p_{j_{l}}(x^{(l)}).
\]
It can be readily seen that
\[
p_{j_{1},\ldots,j_{d}}( y_{i_{1},\ldots,i_{d}} ) =
\prod_{l=1}^{d} \delta_{i_{l}j_{l}}.
\]
As in the 1D case, a polynomial $u(x)$ defined on $\kappa$ can be expressed using the Lagrange polynomials
as
\begin{equation}
  u(x) = \sum_{1\le j_{1},\ldots,j_{d}\le N_{g}} p_{j_{1},\ldots,j_{d}}(x) 
  u(y_{j_{1},\ldots,j_{d}}) 
  \equiv \sum_{1\le j_{1},\ldots,j_{d}\le N_{g}} p_{j_{1},\ldots,j_{d}}(x) 
  u_{j_{1},\ldots,j_{d}}. 
  \label{eqn:vectorDiscretizeDD}
\end{equation}
Denote by $W^{[d]}=\text{diag}[w^{[d]}]$ as a matrix of size
$N_{g}^{d}\times N_{g}^{d}$, the inner product $(u,v)_{\kappa}$ can be
written as
\begin{equation}
  (u,v)_{\kappa} = \sum_{1\le j_{1},\ldots,j_{d}\le N_{g}}
  u_{j_{1},\ldots,j_{d}} v_{j_{1},\ldots,j_{d}} w^{[d]}_{j_{1},\ldots,j_{d}}
  = u^{T} W^{[d]} v.
  \label{eqn:innerproductDD}
\end{equation}

The Lagrange polynomials $p_{j_{1},\ldots,j_{d}}(x)$ can be used to
define $d$ differentiation matrices, defined as
\begin{equation}
  D_{l}^{[d]} = \left(\bigotimes_{i=1}^{l-1} I\right)\bigotimes D
  \bigotimes \left(\bigotimes_{i=l+1}^{d} I\right).
  \label{eqn:DifferentiationMatrixDD}
\end{equation}
Here $I$ is an $N_{g}\times N_{g}$ identity matrix.
$D_{l}^{[d]}$ can be understood as the discretized differential operator
$\partial_{l},1\le l \le d$.  Similar to
Eq.~\eqref{eqn:Differentiation1Dapply}, we denote by
$\partial_{l} v$ a column vector with its entries defined as below
\[
(\partial_{l} v)_{j_{1},\ldots,j_{d}} = (\partial_{l}
v)(y_{i_{1},\ldots,i_{d}}),
\]
then $\partial_{l} v$ can be expressed in the linear algebra notation
as
\begin{equation}
  \partial_{l} v = D_{l}^{[d]} v.
  \label{eqn:DifferentiationDDapply}
\end{equation}
Therefore the inner product $(\nabla u,\nabla v)_{\kappa}$ can be computed
as
\begin{equation}
  (\nabla u,\nabla v)_{\kappa} = u^{T} \left(\sum_{l=1}^{d} 
  (D_{l}^{[d]})^{T} W^{[d]} D_{l}^{[d]}\right) v.
  \label{eqn:GradientDDinnerproduct}
\end{equation}
The inner product $(u,v)_{\star,\kappa}$ can be evaluated similar to
Eq.~\eqref{eqn:star1Dinnerproduct} as
\begin{equation}
  (u,v)_{\star,\kappa} = u^{T} \left( \sum_{l=1}^{d} 
  (D_{l}^{[d]})^{T} W^{[d]} D_{l}^{[d]} + w^{[d]} \frac{1}{|\kappa|}
  (w^{[d]})^{T}\right) v,
  \label{eqn:starDDinnerproduct}
\end{equation}
with $|\kappa|=h^{d}$.

In order to evaluate the inner product on the boundary $\partial \kappa$,
we define $d$ weight vectors corresponding to the $(d-1)$ dimensional
surface for each dimension $l$ ($l=1,\ldots,d$), denoted by
$\wt{w}_{l}^{[d]}$ with the
expression
\begin{equation}
  \left(\wt{w}_{l}^{[d]}\right)_{j_{1},\ldots,j_{d}} 
  = \begin{cases}
    w^{[d-1]}_{j_{1},\ldots,j_{l-1},j_{l+1},\ldots,j_{d}}, &
    j_{l}=1 \quad \text{or} \quad j_{l}=N_{g},\\
    0, & 1<j_{l}<N_{g}.
  \end{cases}
  \label{}
\end{equation}
Define $\wt{W}_{l}^{[d]}=\text{diag}\left[\wt{w}_{l}^{[d]}\right]$, then
the inner product on the boundary can be expressed as
\begin{equation}
  (u,v)_{\partial \kappa} = u^{T} \left(\sum_{l=1}^{d}
  \wt{W}_{l}^{[d]}\right) v,
  \label{eqn:innerproductBoundaryDD}
\end{equation}
and
\begin{equation}
  (\nabla u \cdot n_{\kappa},\nabla v \cdot n_{\kappa})_{\partial \kappa} = 
  u^{T} \left(\sum_{l=1}^{d}
  (D_{l}^{[d]})^{T} \wt{W}_{l}^{[d]} D_{l}^{[d]}\right) v.
  \label{eqn:innerproductGradientBoundaryDD}
\end{equation}

Now we are ready to use the finite dimensional representation of the
inner products to evaluate the constants $\dkN,\rak,\rbk$.  

\subsection{Estimation of $\dkN$}

Recall that
\[
	(\dkN)^2 = \sup_{v_N\in\VN(\kappa)} \frac{\|\grad v_N\cd
  n_\kappa\|^2_{\partial\kappa}}{\snorm{v_N}^2}
  \equiv \sup_{v_N\in\VN(\kappa)} 
  \frac{(\grad v_N\cd n_\kappa, \grad v_N\cd n_\kappa)_{\partial\kappa}}{(v_{N},v_{N})_{\star,\kappa}}.
\]
Using Eq.~\eqref{eqn:innerproductGradientBoundaryDD} and
Eq.~\eqref{eqn:starDDinnerproduct}, we have
\begin{equation}
  (\dkN)^2 = \sup_{v_N\in\VN(\kappa)} \frac{v_{N}^{T} M_{\delta}
  v_{N}}{v_{N}^{T} K v_{N}}.
  \label{eqn:dkN1}
\end{equation}
Here
\begin{equation}
  M_{\delta} = \sum_{l=1}^{d} 
  (D_{l}^{[d]})^{T} \wt{W}_{l}^{[d]} D_{l}^{[d]}, 
  \label{eqn:MdeltaMat}
\end{equation}
\begin{equation}
  K = \sum_{l=1}^{d} (D_{l}^{[d]})^{T} W^{[d]} D_{l}^{[d]} +
  w^{[d]} \frac{1}{|\kappa|} (w^{[d]})^{T}.
  \label{eqn:KMat}
\end{equation}
Let $\{\varphi_{1}(x),\ldots,\varphi_{N}(x)\}$ be a set of basis
functions of the finite dimensional space $\VN(\kappa)$. We
denote by $\varphi_{i} (i=1,\ldots,N)$ the column vector corresponding
to the values of $\varphi_{i}(x)$ evaluated at the LGL grid points, and
denote
by 
\begin{equation}
  \Phi=[\varphi_{1},\ldots,\varphi_{N}],
  \label{eqn:PhiMat}
\end{equation}
the collection of all column
vectors which is an $N_{g}^{d}\times N$ matrix.  Then for any
vector $v(x)\in \VN(\kappa)$, the corresponding column vector $v$ can be
represented as
\[
  v = \Phi c,
\]
where $c$ is a coefficient vector of size $N$. Then
Eq.~\eqref{eqn:dkN1} can be rewritten as
\begin{equation}
  (\dkN)^2 = \sup_{c\in\mathbb{R}^{N}} \frac{c^{T} (\Phi^{T} M_{\delta}
  \Phi) c}{c^{T} (\Phi^T K_{\delta} \Phi) c}.
  \label{eqn:dkN2}
\end{equation}
Eq.~\eqref{eqn:dkN2} can be solved as an eigenvalue problem, 
\begin{equation}
  \Phi^{T} M_{\delta} \Phi c = \lambda \Phi^T K_{\delta} \Phi c,
  \label{eqn:dkN3}
\end{equation}
and $(\dkN)^2$ is equal to the largest eigenvalue $\lambda$.  Since the
size of the matrix $\Phi^{T} M_{\delta} \Phi$ is $N\times N$ and
$N$ is relatively small, Eq.~\eqref{eqn:dkN3} can be solved as a
generalized eigenvalue problem using dense linear algebra.

\subsection{Estimation of $\rak,\rbk$}

Recall that 
\[
	\rak^2 = \sup_{\substack{v\in H^1(\kappa),\\ v\perp\VN(\kappa)}}
  \frac{\|v\|^2_{\kappa}}{\snorm{v}^2} = 
  \sup_{\substack{v\in H^1(\kappa),\\ v\perp\VN(\kappa)}}
  \frac{(v,v)_{\kappa}}{(v,v)_{\star,\kappa}},
\]
then using Eq.~\eqref{eqn:innerproductDD},
Eq.~\eqref{eqn:starDDinnerproduct} and the density arguments shown above, it can be shown that
\begin{equation}
  \sup_{\substack{v\in \mathbb P_{N_g},\\ v\perp\VN(\kappa)}}
  \frac{v^{T} M_{{\tt a}} v}{v^{T} K v}
  \quad
 \overset{N_g\to \infty}{\longrightarrow} 
 \quad
 \rak^2,
  \label{eqn:ak1}
\end{equation}
where $M_{{\tt a}} = W^{[d]}$, and $K$ is given in Eq.~\eqref{eqn:KMat}. We
can express the orthogonality condition $v\perp\VN(\kappa)$ in terms of
a projection operator $Q=\mathbb{I}-\pN$ so that for any $v\in H^1(\kappa)$,
$Qv\perp\VN(\kappa)$, where $\mathbb{I}$ is the identity operator.
Denoting by $\Phi$ as in Eq.~\eqref{eqn:PhiMat} the collection of
spanning vectors of the space $\VN(\kappa)$, then using the Lagrange
polynomials corresponding to the LGL grid points as a basis, the
projection operator $\pN$ can be expressed as an $N_{g}^{d}\times
N_{g}^{d}$ matrix 
\begin{equation}
  \pN = \Phi (\Phi^{T} K \Phi)^{-1} \Phi^{T} K \equiv \Phi
  \Psi^{T}.
  \label{eqn:pNMat}
\end{equation}
where $\Psi = K \Phi (\Phi^{T} K \Phi)^{-1}$.
Therefore the $\pN$ is a low rank matrix with rank $N$. The projection
operator $Q$ and its adjoint operator $Q^{T}$ expressed in the basis of Lagrange polynomials become
\begin{equation}
  Q = I - \Phi \Psi^{T}, \quad Q^{T} = I - \Psi \Phi^{T}.
  \label{eqn:QMat}
\end{equation}
Using Eq.~\eqref{eqn:QMat}, the computation of $\rak$ can be simplified
as
\begin{equation}
  \rak^2 \REV{\approx} \sup_{v\in\mathbb{R}^{N_{g}^{d}}}
 \frac{v^{T} Q^{T} M_{{\tt a}} Q v}{v^{T} Q^{T} K Q v}.
  \label{eqn:ak2}
\end{equation}
In other words, $\rak^2$ corresponds to the largest eigenvalue of the
generalized eigenvalue problem
\begin{equation}
  Q^{T} M_{{\tt a}} Q v = \lambda  Q^{T} K Q v.
  \label{eqn:ak3}
\end{equation}

From a computational point of view, there are two major differences
between Eq.~\eqref{eqn:dkN3} and ~\eqref{eqn:ak3}.  First, the dimension
of the matrices in Eq.~\eqref{eqn:dkN3} is $N\times N$, and the
dimension of the matrices in Eq.~\eqref{eqn:ak3} is $N_{g}^{d}\times
N_{g}^{d}$.  For 3D simulation, if $N_{g}=30$ then $N_{g}^{d}=27000$,
and the corresponding eigenvalue problem is very costly to solve if
$Q^{T} M_{{\tt a}} Q$ and $Q^{T} K Q$ are treated as dense matrices.
Second, the matrix $\Phi^{T}K\Phi$ in Eq.~\eqref{eqn:dkN3} is a positive
definite matrix since $K$ is positive definite, and the
problem~\eqref{eqn:dkN3} can be solved directly as a dense generalized
eigenvalue problem.  On the other hand, $Q^{T} K Q$ is a rank deficient
matrix with the rank of its kernel being $N$.  Therefore it can 
potentially cause a large numerical error if Eq.~\eqref{eqn:ak3} is solved
directly as a dense generalized eigenvalue problem.

In order to overcome the two difficulties mentioned above, we note that for
any vector $v$, the computational cost for the matrix vector
multiplication $Qv, Q^{T}v, M_{{\tt a}}v, K v$ is only proportional to
$N_{g}^{d}$ thanks to the low rank representation of the
operators.  Therefore Eq.~\eqref{eqn:ak3} can be solved
using iterative methods.  Another advantage of using iterative methods is that since
we only need the largest eigenvalue corresponding to
Eq.~\eqref{eqn:ak3}, at the $k$-th step of the CG iteration we only need
to keep three vectors: the current approximation of eigenvector
$v^{(k)}$, the conjugate direction $p^{(k)}$ and the residual
$r^{(k)}$.  Even though the matrix $Q^{T} K Q$ is singular, the
projection onto the $3$ dimensional subspace $[v^{(k)}, p^{(k)},
r^{(k)}]$ is usually well conditioned.  In practice we use the Locally
Optimal Block Preconditioned Conjugate Gradient (LOBPCG)
method~\cite{Knyazev2001} (with block size equal to $1$)
for evaluating the largest eigenvalue for Eq.~\eqref{eqn:ak3}. It
should be noted that since there is no apparent preconditioner that can
be applied efficiently to solve Eq.~\eqref{eqn:ak3}, the convergence of
the largest eigenvalue may be slow.  However, we should keep in mind that the estimation of
$\rak,\rbk$ is only used in the a posteriori error estimator, and only
low accuracy is needed.  In fact $\rak,\rbk$ is already very accurate in
the sense of the preconstant in the estimator even if the relative error
is $10\%$.  Therefore the slow convergence of the conjugate gradient
method is compensated by the low accuracy required in the computation of
the constants.  

The constant $\rbk$ can be estimated similarly to $\rak$. Recall that
\[
	\rbk^2 = \sup_{\substack{v\in H^1(\kappa),\\ v\perp\VN(\kappa)}}
  \frac{\|v\|^2_{\partial\kappa}}{\snorm{v}^2} = 
  \sup_{\substack{v\in H^1(\kappa),\\ v\perp\VN(\kappa)}}
  \frac{(v,v)_{\partial\kappa}}{(v,v)_{\star,\kappa}},
\]
and using the same projection operator $Q$, $\rbk$ can be expressed as
\begin{equation}
  \rbk^2 \REV{\approx} \sup_{v\in\mathbb{R}^{N_{g}^{d}}}
 \frac{v^{T} Q^{T} M_{{\tt b}} Q v}{v^{T} Q^{T} K Q v},
  \label{eqn:bk2}
\end{equation}
with $M_{{\tt b}}=\sum_{l=1}^{d} \wt{W}_{l}^{[d]}$. Similar to
Eq.~\eqref{eqn:ak3}, $\rbk^2$ can be solved as the largest eigenvalue
of
\begin{equation}
  Q^{T} M_{{\tt b}} Q v = \lambda  Q^{T} K Q v.
  \label{eqn:bk3}
\end{equation}
Eq.~\eqref{eqn:bk3} can be solved using the same iterative strategy as
for obtaining $\rak$.

\section{Numerical results}\label{sec:numerical}

In this section we test the effectiveness of the a posteriori error
estimators.  The test program is written in MATLAB, and all results are
obtained on a 2.7 GHz Intel processor with 16 GB memory.  
All numerical results are performed using the symmetric bilinear form
($\theta=1$).
The effectiveness of the upper bound and lower bound on the global domain will
be justified by comparing $\tnorm{u-u_{N}}$ and $\eta$, and
by comparing $\tnorm{u-u_{N}}$ and $\xi$, respectively.  It should be
noted that although our theory does not directly predict the effectiveness of
the estimator on each local element $\kappa$, we can measure the
local effectiveness of the upper and lower bound on each local element
$\kappa$ by defining 
\begin{equation}
  C_{\eta}(\kappa) = \frac{ \etaR+\etaF+\etaJ  }{\tnorm{u-u_{N}}_{\kappa}},
  \quad
  C_{\xi}(\kappa) = \frac{ \xi_\kappa }{\tnorm{u-u_{N}}_{\kappa}},
  \label{eqn:Ceta}
\end{equation}
where the broken energy norm $\tnorm{u-u_{N}}_{\kappa}$ is defined
according to Eq.~\eqref{eq:DGnorm}.  

The numerical results are organized as follows.  In
section~\ref{subsec:polyconst}, we apply the general approach developed
in section~\ref{sec:computation} to compute the constants $\rak, \rbk,
\dkN$ for polynomial basis functions, 
and verify that the scaling properties of the numerically computed
constants match the analytic results known in the
literature~\cite{Schwab1998}.  In section~\ref{subsec:posdef}, we illustrate the
behavior of the upper bound and the lower bound error
estimates for second order PDEs associated with positive definite
operators. We then demonstrate the results for indefinite operators in
section~\ref{subsec:indef}.  In the a posteriori error estimates of both
the upper bound and the lower bound, we make the assumption that the
non-computable number \REV{$\dku(\uN)$} can be approximated by $\dkN$ without
significant loss of effectiveness.  We justify such treatment in
section~\ref{subsec:justifydku} by directly calculating \REV{$\dku(\uN)$} using the
numerically computed reference solution.

Our test problems include both one dimensional (1D) and two dimensional
(2D) domains with periodic boundary conditions.  Our
non-polynomial basis functions are generated from the adaptive local
basis (ALB) set~\cite{LinLuYingE2012} in the DG framework.  
The ALB set was originally proposed to systematically reduce the number
of basis functions used to solve Kohn-Sham density functional theory
calculations, and in this section we demonstrate its usage to solve
second order linear PDEs.  We denote by $N$ the number of ALBs per
element.
For operators in the form of $A=-\Delta+V$ with periodic boundary
condition, the basic idea of the ALB set is to use eigenfunctions
computed local domains as basis functions corresponding to the lowest
few eigenvalues.  The eigenfunctions are associated with the same
operator $A$, but with modified boundary conditions on the local domain.
More specifically, in
a $d$-dimensional space, for each element $\kappa$, we form an
\textit{extended element} $\widetilde{\kappa}$ consisting of $\kappa$
and its $3^{d}-1$ neighboring elements in the sense of periodic boundary
condition. On $\widetilde{\kappa}$ we solve the eigenvalue problem
\begin{equation}
  -\Delta \widetilde{\varphi}_{i} + V \widetilde{\varphi}_{i} = \lambda_{i} 
  \widetilde{\varphi}_{i}.
  \label{eqn:localeig}
\end{equation}
with periodic boundary condition on $\partial \widetilde{\kappa}$.
\REV{This eigenvalue problem can be solved using standard basis set such
as finite difference, finite elements, or planewaves. Here we solve the
local eigenvalue problem~\eqref{eqn:localeig} using planewaves which
naturally satisfy periodic boundary conditions.
Since this eigenvalue problem is solved on a extended element
$\widetilde{\kappa}$ the computational cost is not large.}
The
collection of eigenfunctions (corresponding to lowest $N$ eigenvalues)
are restricted from $\widetilde{\kappa}$ to $\kappa$, i.e. 
\[
\varphi_{i}(x) = \begin{cases}
  \left[\widetilde{\varphi}_{i}\right]\vert_{\kappa}(x), & x\in \kappa;\\
  0,&\text{otherwise}.
\end{cases}
\]
After orthonormalizing $\{\varphi_{i}\}$ locally on each element
$\kappa$ and removing the linearly dependent functions, the resulting
set of orthonormal functions are called the ALB functions.  

Since periodic boundary condition is used on the global domain $\Omega$,
in all the calculations, the reference solution, which can be treated as
a numerically exact solution, is solved using a planewave basis set with
a sufficiently large number of planewaves.  The ALB basis set is also
computed using a sufficiently large number of planewaves on the extended
element $\widetilde{\kappa}$. Then a Fourier interpolation procedure
is carried out from $\widetilde{\kappa}$ to the local element
$\kappa$ on a Legendre-Gauss-Lobatto (LGL) for accurate numerical
integration.

\subsection{Estimating the constants for polynomial basis
functions}\label{subsec:polyconst}

Although the main purpose of this paper is to design a posteriori error
estimator for non-polynomial basis functions, the computational
strategies discussed in section~\ref{sec:computation} can be applied to
polynomial functions as well. Let $\kappa=[0,h]^{d}$ and
$\VN(p;\kappa)=\text{span}\{\prod_{l=1}^{d} x_{l}^{j_{l}}, j_{l}\in
\mathbb{N}, \sum_{l=1}^{d}j_{l}\le p \}$ be the space spanned by
polynomials with degree less than or equal to $p$. 
Then the asymptotic scaling of $\rak,\rbk,\dkN$ with respect to $h$ and
$p$ is known~\cite{HoustonSchotzauWihler2007}
\begin{equation}
  \rak^2 \sim \frac{h^2}{p^2}, \quad \rbk^2 \sim \frac{h}{p}, \quad
  \dkN^2 \sim
  \frac{p^2}{h}. 
  \label{eqn:standardHPscale}
\end{equation}
These results are asymptotically correct as
$p\to \infty$, and we will show that the strategy discussed in
section~\ref{sec:computation} leads to the same asymptotic result, but
the result is more accurate in the pre-asymptotic regime due to the
explicit computation of the constants.

From numerical point of view, the scaling with respect to $h$ is
naturally satisfied.  To verify this, we can simply consider a reference element
$\kappa\vert_{h=1}=[0,1]^d$ and scale the weight matrix $W^{[d]}$
and the differentiation matrix $D_{l}^{[d]}$ accordingly. The technique
is the same as that used in~\cite{Schwab1998}.

\begin{figure}[h]
  \begin{center}
    \subfloat[(a)]{\includegraphics[width=0.35\textwidth]{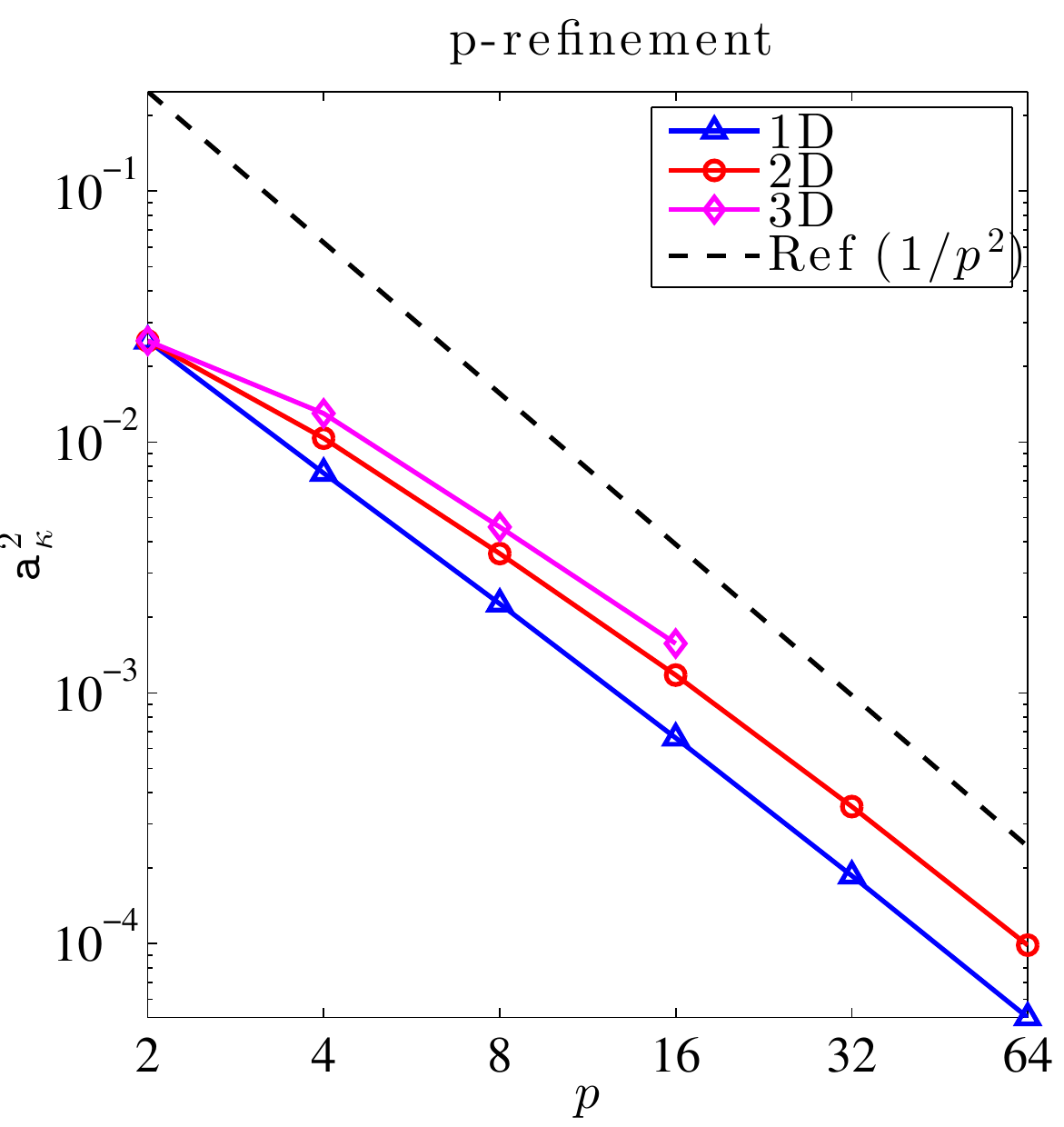}}
    \quad
    \subfloat[(b)]{\includegraphics[width=0.35\textwidth]{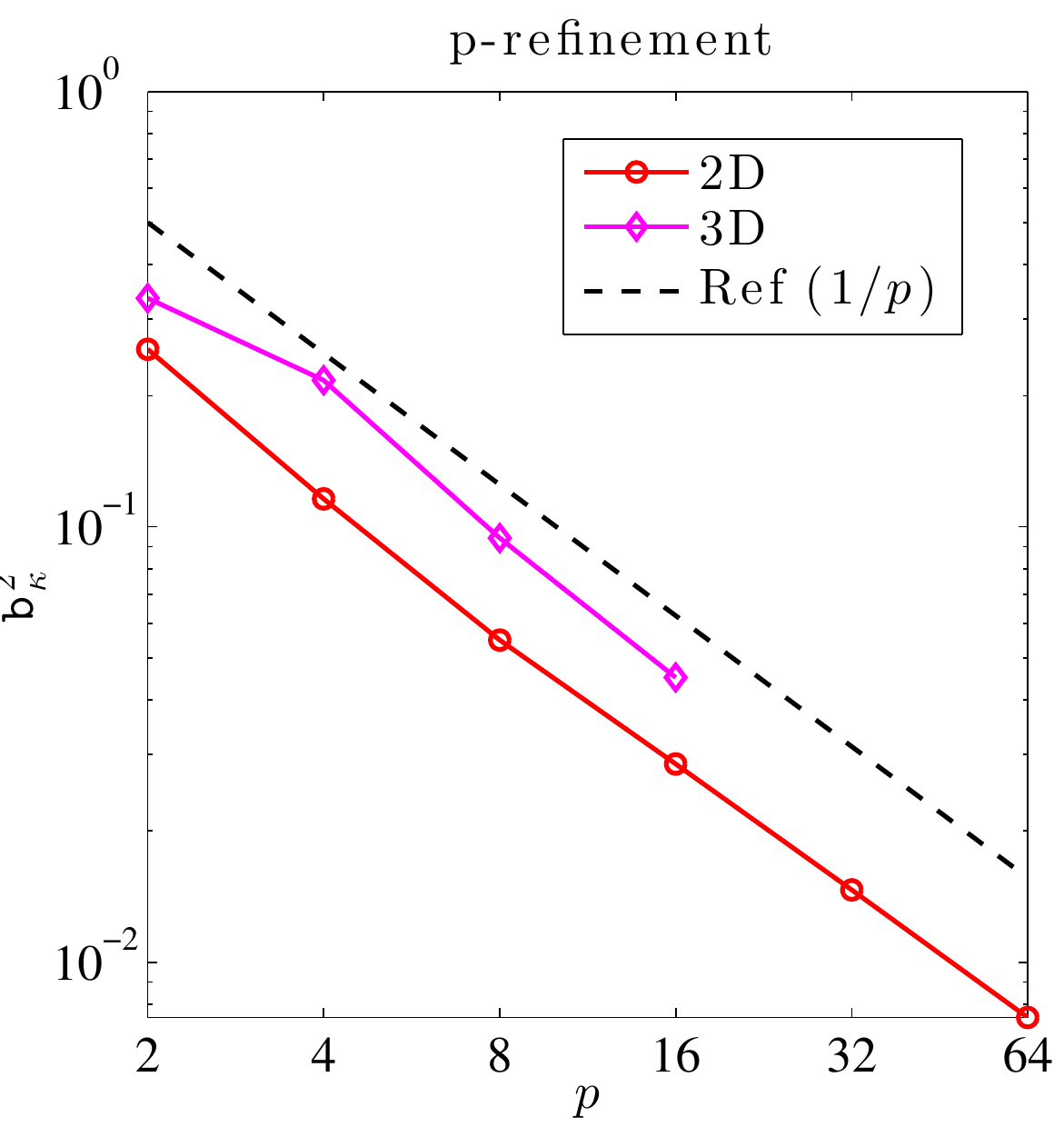}}
    \quad
    \subfloat[(c)]{\includegraphics[width=0.35\textwidth]{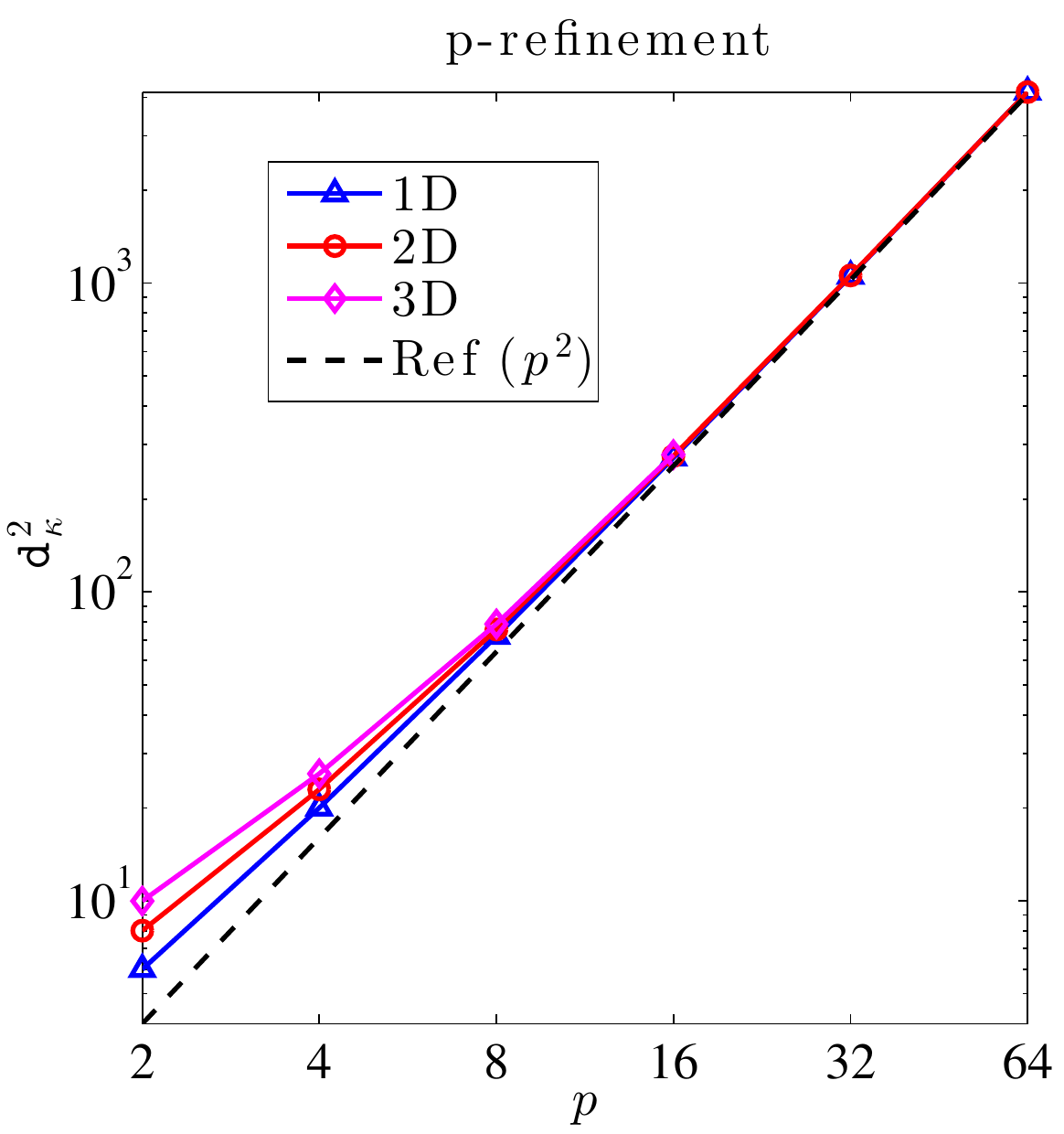}}
  \end{center}
  \caption{Numerically computed constants $\rak^2,\rbk^2,\dkN^2$ with respect to the
  polynomial degrees $p$ in 1D, 2D and 3D.}
  \label{fig:constantPrefine}
\end{figure}

We now directly verify the scaling with respect to $p$ in
Fig.~\ref{fig:constantPrefine}, using
the algorithms presented in section~\ref{sec:computation}. The LGL grid
sizes for 1D, 2D and 3D calculation are chosen to be $100$, $100\times
100$, and $50\times 50\times 50$, respectively.  The largest degree of
polynomials is $64$ for 1D and 2D, and is $16$ for the 3D case.  Note
that in the 3D case, the dimension of $\VN(p=16;\kappa)$ is already
$969$.  Fig.~\ref{fig:constantPrefine} (a) shows the behavior of
$\rak^2$, which asymptotically agrees with the $1/p^2$ scaling.
It is interesting to see that the computed $\rak^2$ can be approximated
by $C \frac{h^2}{p^2}$ where the constant $C$ is around $0.1$.
The recovery of the constant indicates that 
the numerically computed constant $\rak$ can offer a sharper estimator even for
the standard $hp$-refinement.  Similarly Fig.~\ref{fig:constantPrefine}
(b)
shows that $\rbk^2$ asymptotically scales as $1/p$ for 2D and 3D
simulation.  The 1D case is not shown in the picture, since the
numerical value of $\rbk^2$ is already as small as $10^{-20}$ for $p=2$.
This can be interpreted from Proposition~\ref{prop:rbk1D} in the
appendix.  Finally, direct computation in Fig.~\ref{fig:constantPrefine} (c) shows
that $\dkN^2$ asymptotically scales as $p^2$ for all dimensions. Again,
the computed constant $\dkN^2$ differs from the asymptotic scaling in
the pre-asymptotic regime, indicating that the numerically computed
constant should be sharper for low order polynomials ($p\le 4$).

%
%
%

\FloatBarrier

\subsection{Positive definite operators}\label{subsec:posdef}

\begin{figure}[h]
  \begin{center}
    \subfloat[(a)]{\includegraphics[width=0.35\textwidth]{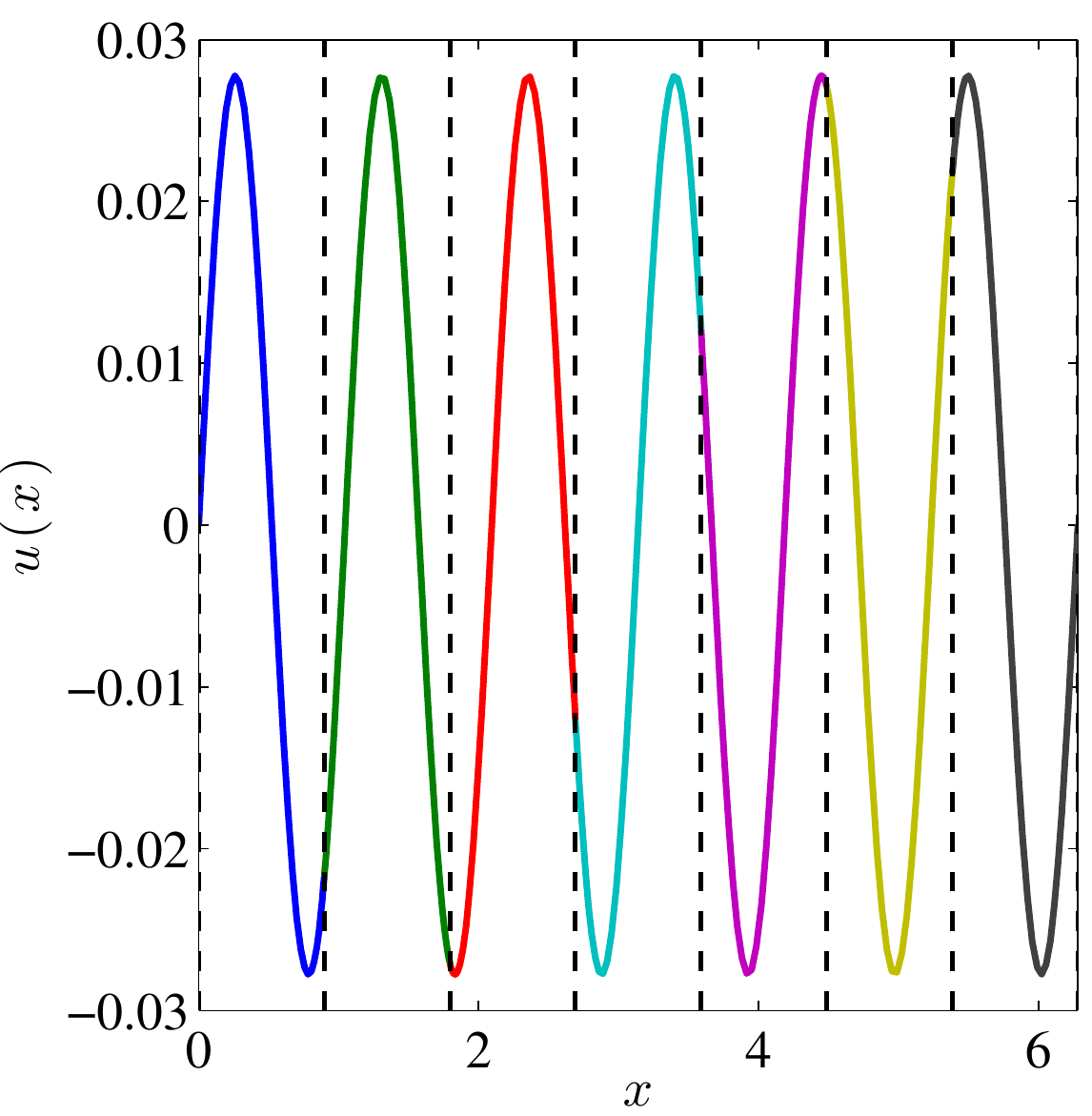}}
    \quad
    \subfloat[(b)]{\includegraphics[width=0.36\textwidth]{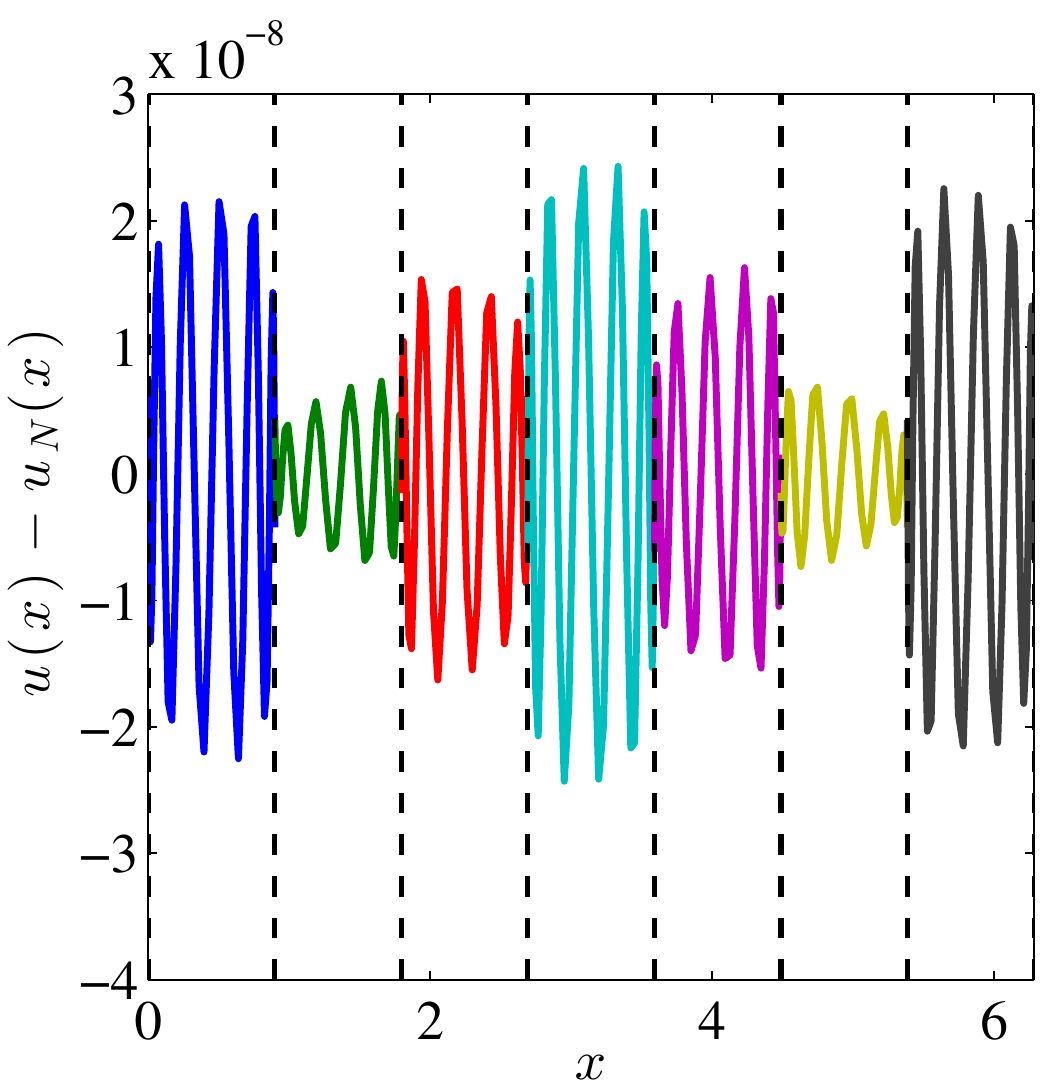}}
  \end{center}
  \caption{(a) The reference solution $u(x)$ corresponding to $V(x)=0.01$ and
  the right hand side $f(x)=\sin(6x)$. (b) Point-wise error between the
  reference solution $u(x)$ and the numerical solution $u_N(x)$ calculated
  using the ALB set with $7$ elements and $N=11$ basis functions per
  element. The domain is partitioned into $7$ elements indicated by
  black dashed lines.}
  \label{fig:uuerrPoisson1D}
\end{figure}

We first demonstrate the effectiveness of the a posteriori error
estimates for a positive definite operator on a 1D domain
\REV{$\Omega=(0,2\pi)$}, using the ALB set as non-polynomial basis functions.
Due to the periodic boundary condition, we choose $V(x)=0.01$ so that
the operator $A=-\Delta+V$ is non-singular and positive definite.  The
right hand side is chosen to be $f(x)=\sin(6x)$ which is periodic on
$\Omega$. In the ALB computation, the domain is partitioned into $7$
elements, as indicated by black dashed lines.
Fig.~\ref{fig:uuerrPoisson1D} shows solution $u$ to Eq.~\eqref{eq:Indef}
and the point-wise error $u-u_{N}$ using $N=11$ ALBs per element.  

\begin{figure}[h]
  \begin{center}
    \subfloat[(a)]{\includegraphics[width=0.37\textwidth]{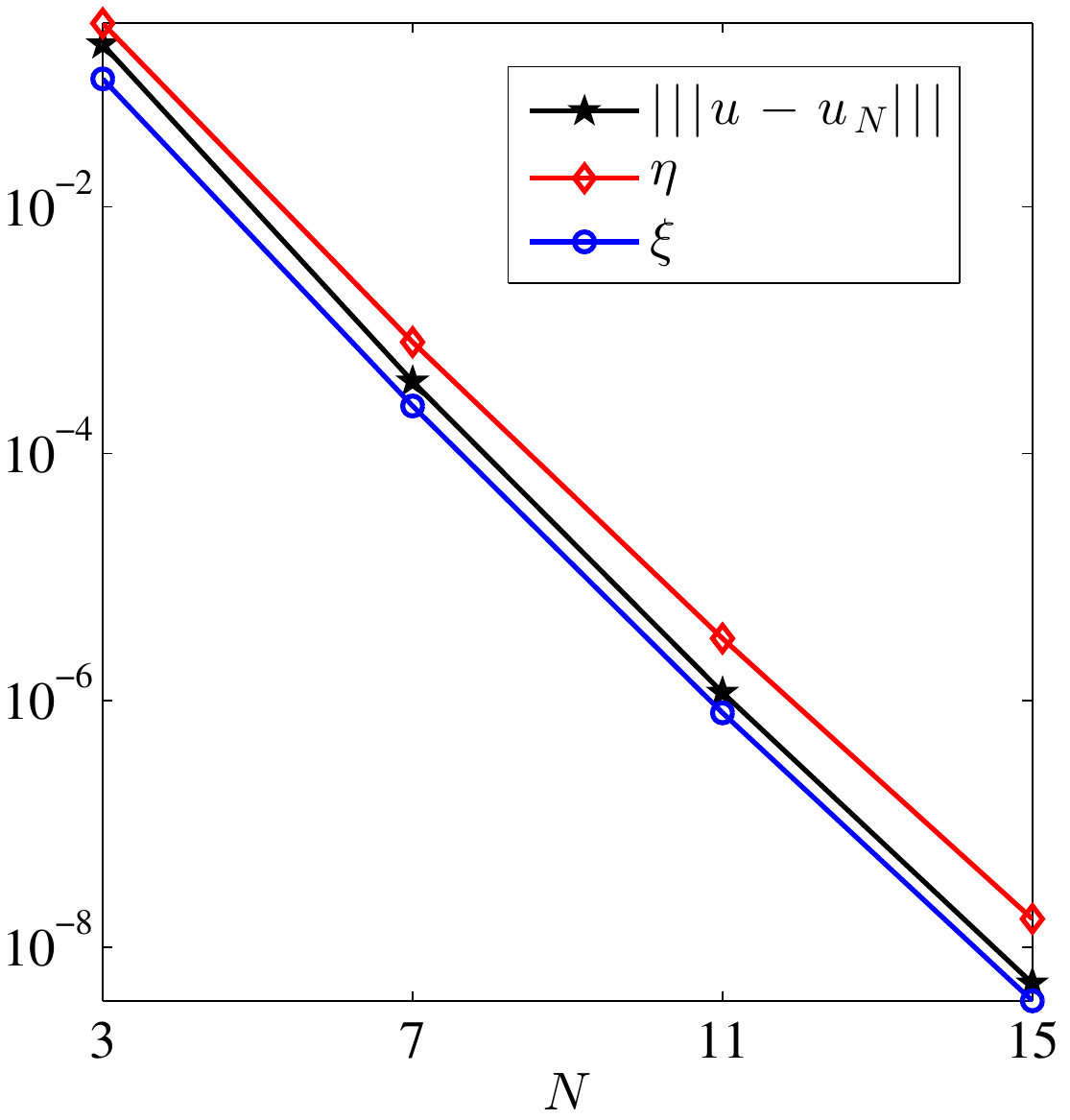}}
    \quad
    \subfloat[(b)]{\includegraphics[width=0.37\textwidth]{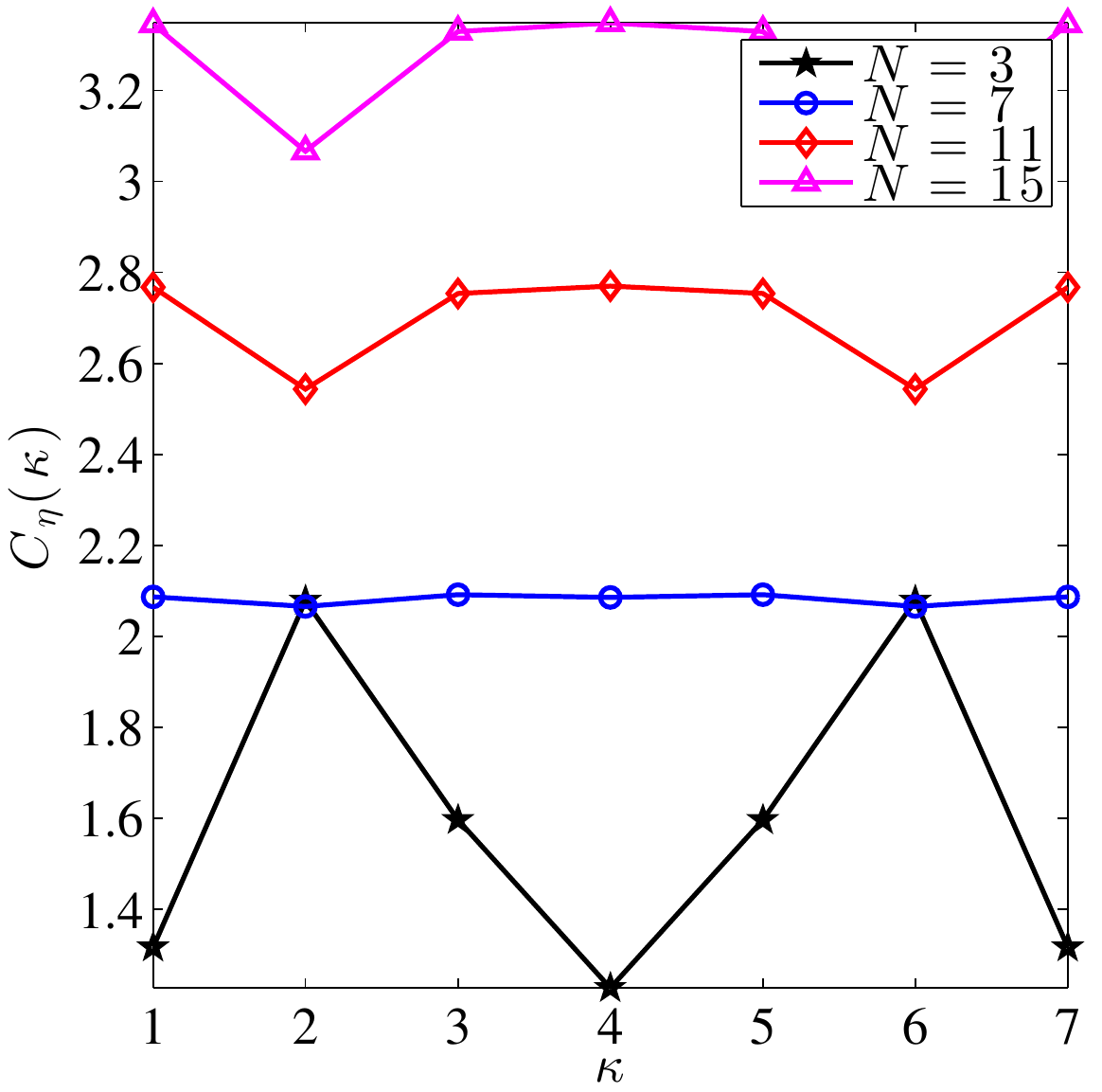}}
    \quad
    \subfloat[(c)]{\includegraphics[width=0.37\textwidth]{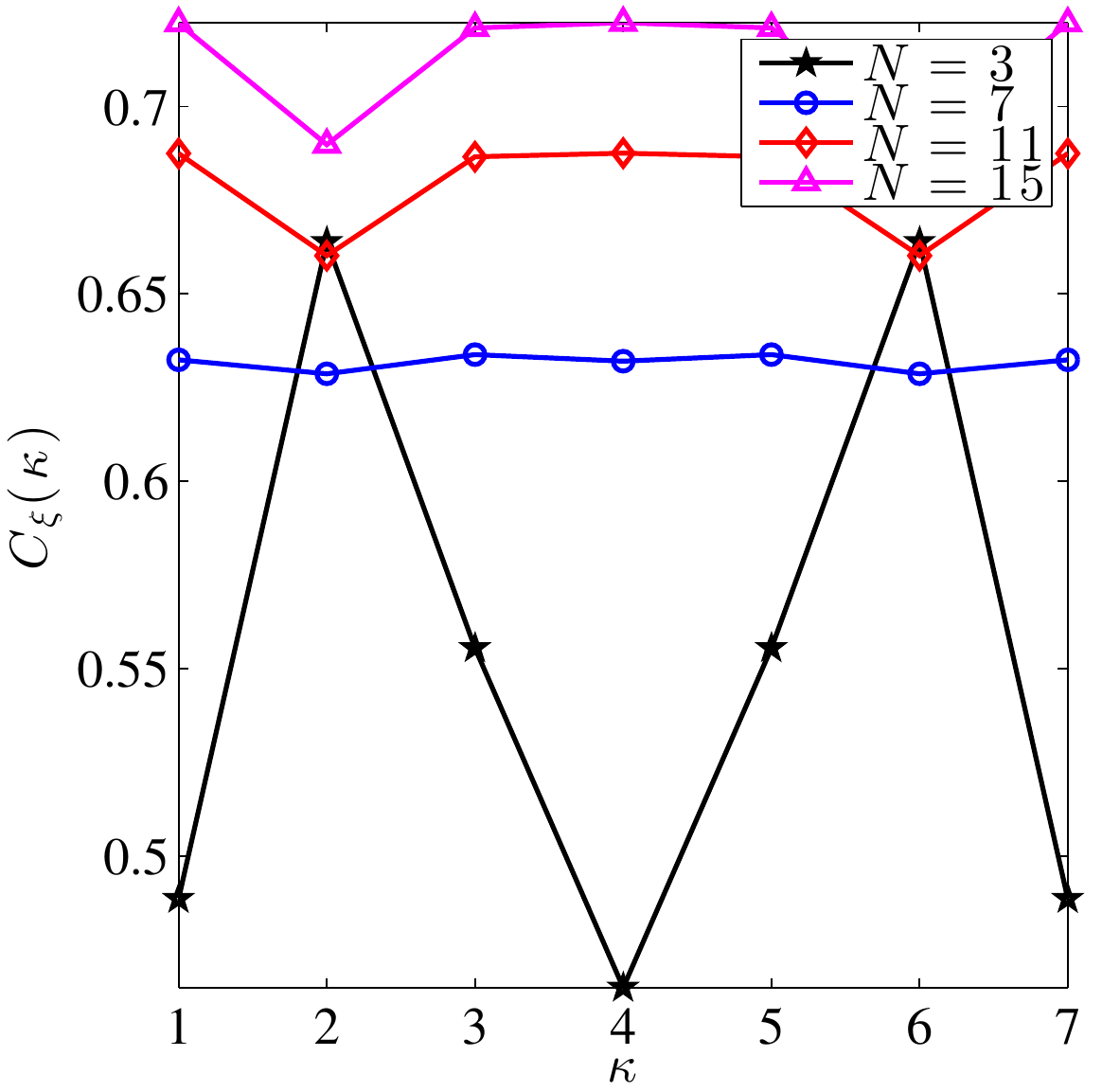}}
  \end{center}
  \caption{(a) \REV{Global error and the upper/lower bound estimator}
  for $V(x)=0.01$ and
  $f(x)=\sin(6x)$. (b) Local effectiveness of the upper bound
  characterized by $C_{\eta}(\kappa)$ for each element. (c) Local
  effectiveness of the lower bound characterized by $C_{\xi}(\kappa)$
  for each element.}
  \label{fig:estPoisson1D}
\end{figure}

Fig.~\ref{fig:estPoisson1D} (a) shows the absolute error in the energy norm, the
upper bound and lower bound estimates as the number of ALBs per element
$N$ increases from $3$ to $15$.  The relative error can be deduced by
comparing Fig.~\ref{fig:estPoisson1D} (a) and
Fig.~\ref{fig:uuerrPoisson1D} (a).
We find that the computed $\eta$ and
$\xi$ are indeed upper and lower bounds of the true error
$\tnorm{u-u_{N}}$ for all $N$ across a wide range of accuracy (from
$10^{-1}$ to $10^{-8}$).  It also
appears that the lower bound estimator $\xi$ follows the true error more
closely than the upper bound estimator $\eta$.  Fig.~\ref{fig:estPoisson1D} (b)
and (c) illustrate the local effectiveness $C_{\eta}(\kappa)$ and
$C_{\xi}(\kappa)$ for each element $\kappa$. Though not guaranteed by
our theory, we observe that $\eta_{\kappa}$ and
$\xi_{\kappa}$ are upper and lower bounds for $\tnorm{u-u_{N}}_{\kappa}$
for each element $\kappa$, respectively.  The effectiveness as
measured by $C_{\eta}(\kappa)$ and $C_{\xi}(\kappa)$ depends only
weakly on the number of adaptive local basis functions, or the accuracy
of the numerical solution.

%
%

\begin{figure}[h]
  \begin{center}
    \subfloat[(a)]{\includegraphics[width=0.4\textwidth]{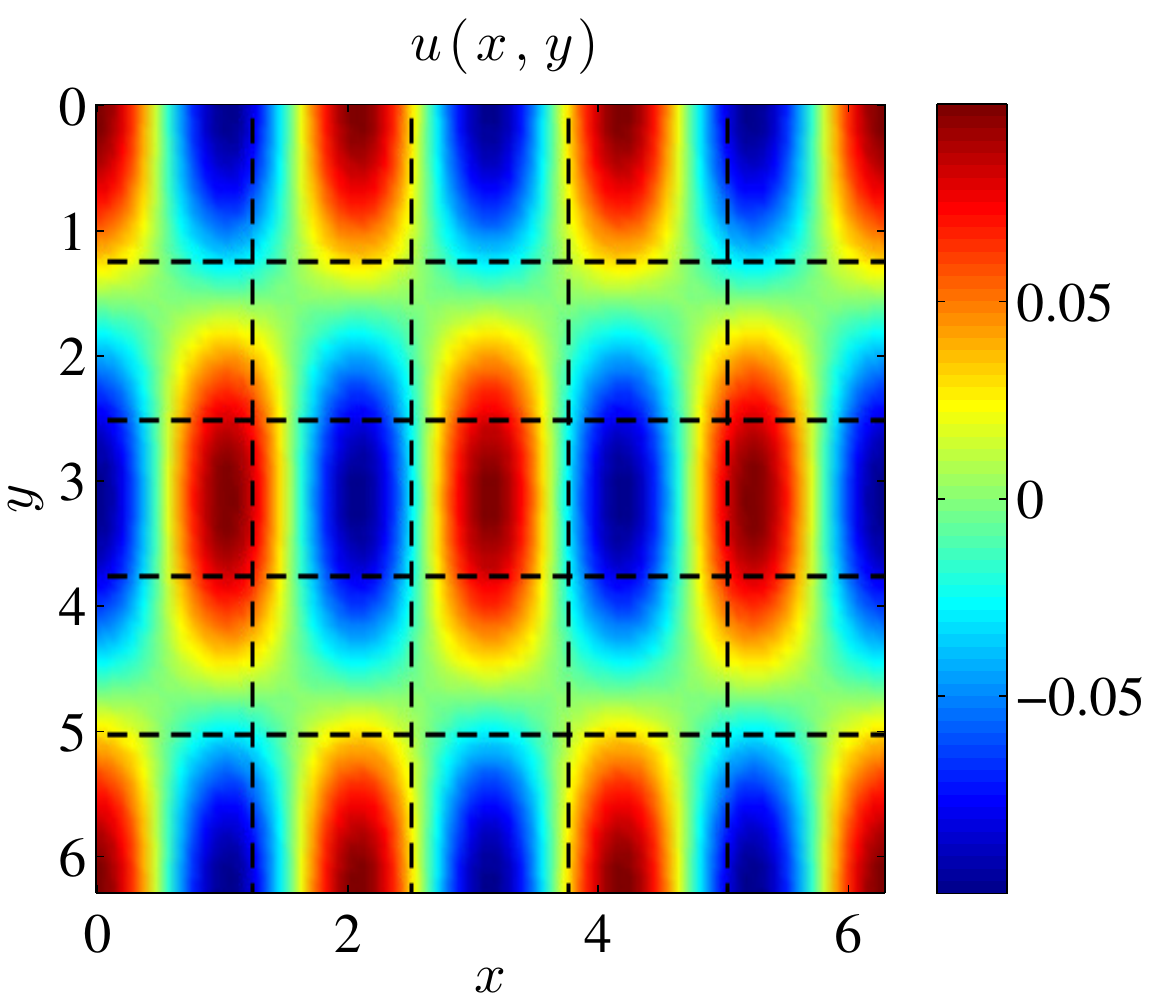}}
    \quad
    \subfloat[(b)]{\includegraphics[width=0.4\textwidth]{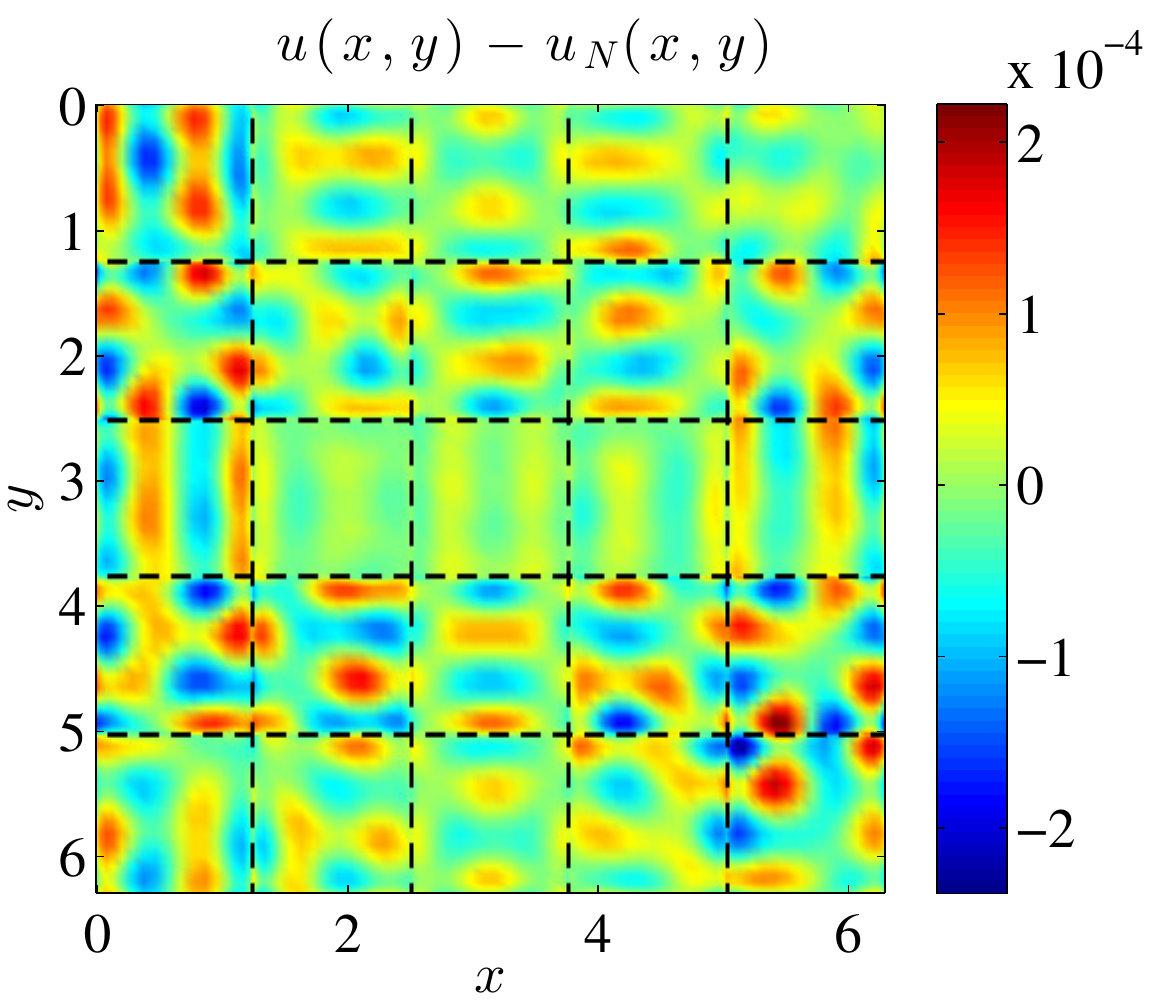}}
  \end{center}
  \caption{(a) The reference solution $u(x,y)$ corresponding to $V(x,y)=0.01$ and
  $f(x,y)=\cos(3x) \cos(y)$. (b) Point-wise error between the reference solution
  $u(x,y)$ and the numerical solution $u_N(x,y)$ calculated using the ALB
  set with $5\times 5$ elements and $N=31$ basis functions per element.}
  \label{fig:uuerrPoisson2D}
\end{figure}

Our next example is to solve a 2D problem with \REV{$\Omega=(0,2\pi)\times
(0,2\pi)$}. Again we choose $V(x,y)=0.01$ so that $A=-\Delta+V$ is
non-singular and positive definite.  The right hand side is $f(x,y)=\cos(3x) \cos(y)$,
which satisfies the periodic boundary condition.
Fig.~\ref{fig:uuerrPoisson2D} shows the reference solution $u$ to
Eq.~\eqref{eq:Indef} and the point-wise error $u-u_{N}$ using $N=31$
ALBs per element.  In the ALB computation, the domain is partitioned
into $5\times 5$ elements, indicated by black dashed lines.

\begin{figure}[h]
  \begin{center}
    \subfloat[(a)]{\includegraphics[width=0.28\textwidth]{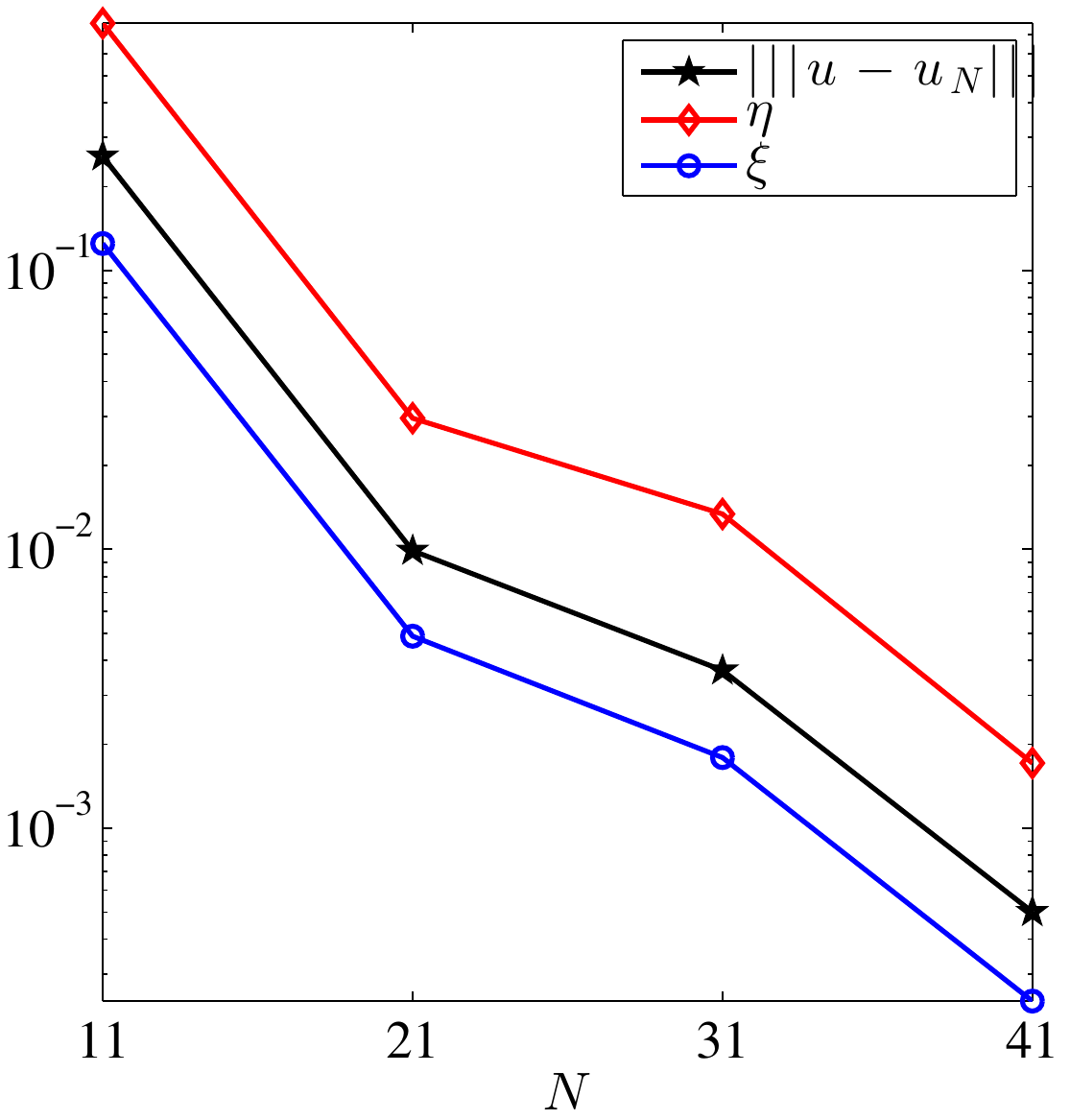}}
    \subfloat[(b)]{\includegraphics[width=0.35\textwidth]{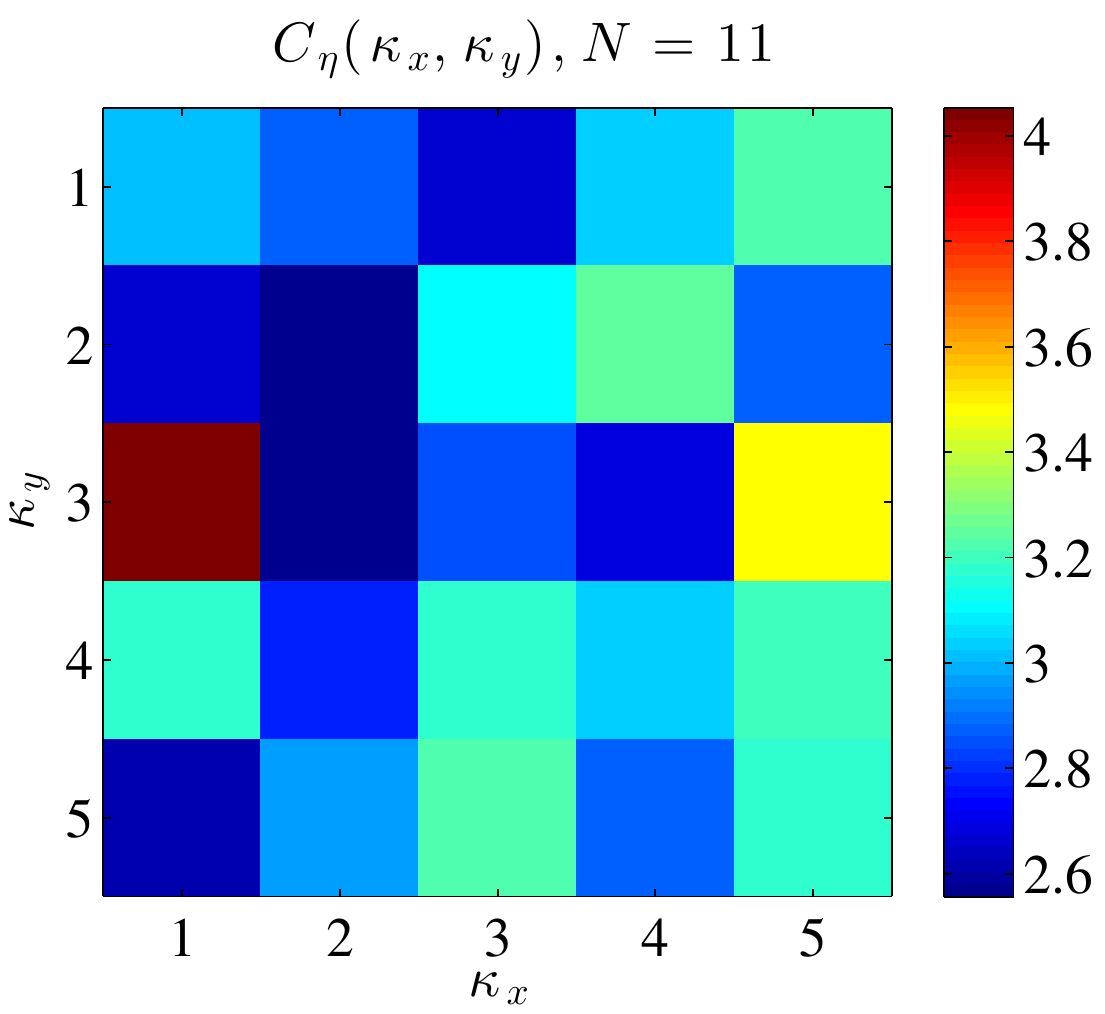}}
    \subfloat[(c)]{\includegraphics[width=0.35\textwidth]{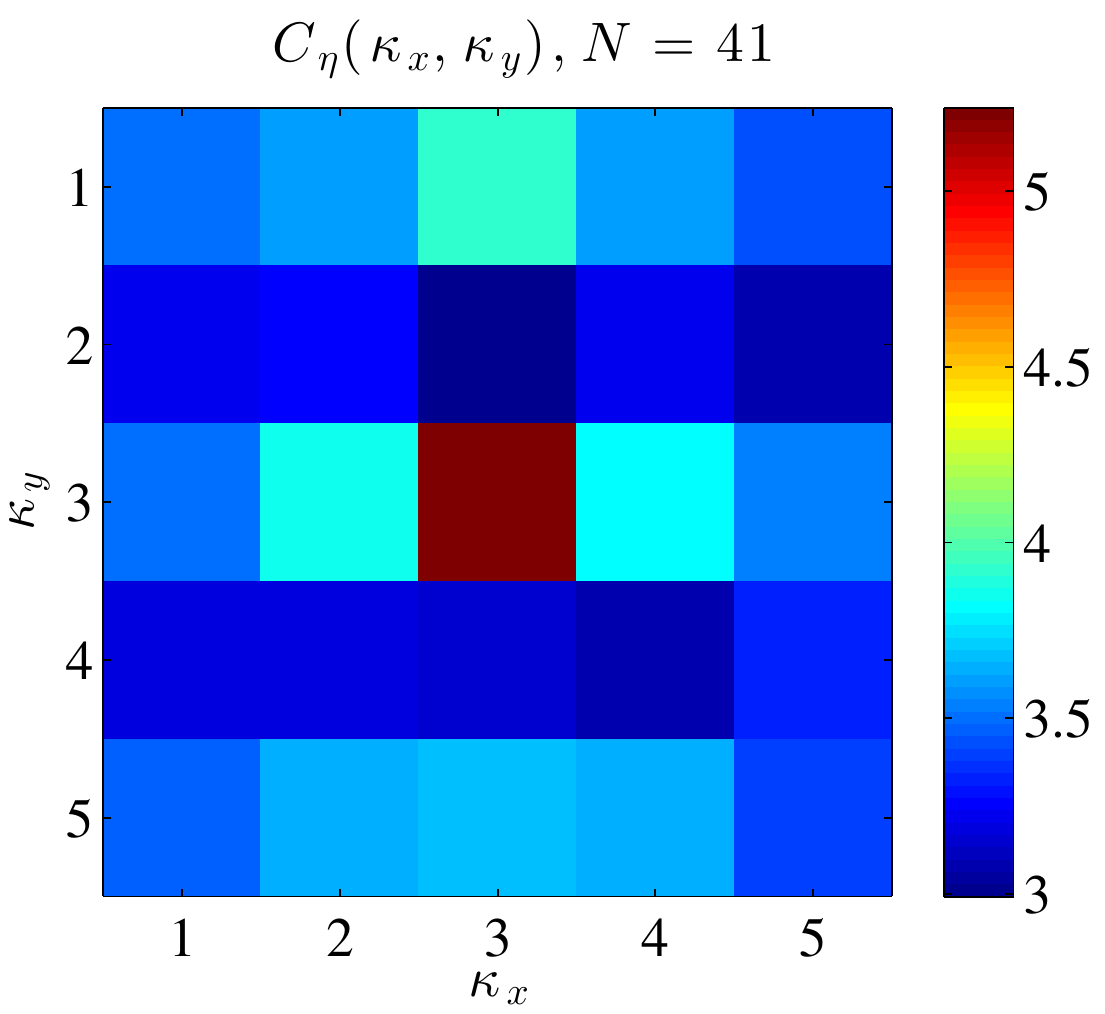}}

    \subfloat[(d)]{\includegraphics[width=0.35\textwidth]{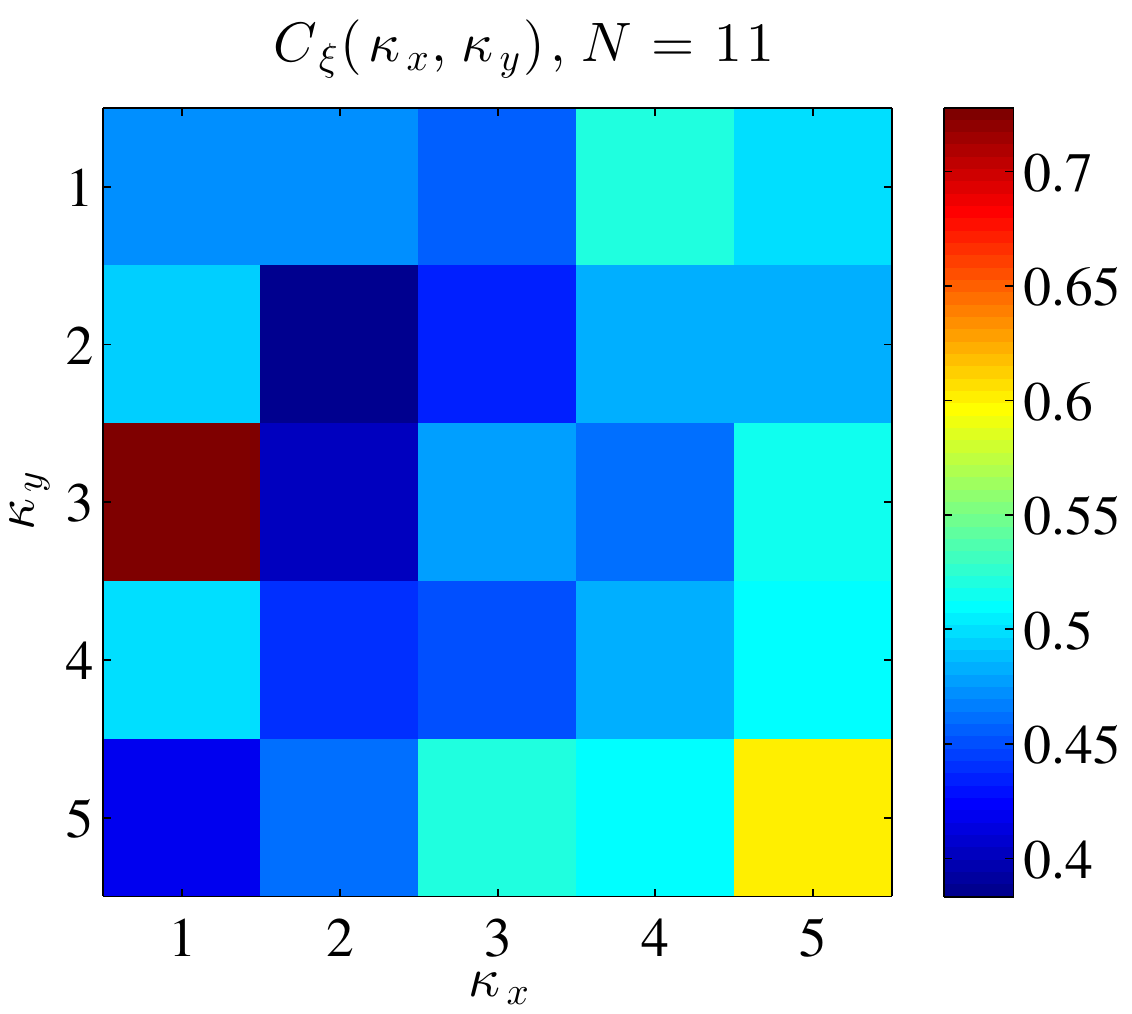}}
    \quad
    \subfloat[(e)]{\includegraphics[width=0.35\textwidth]{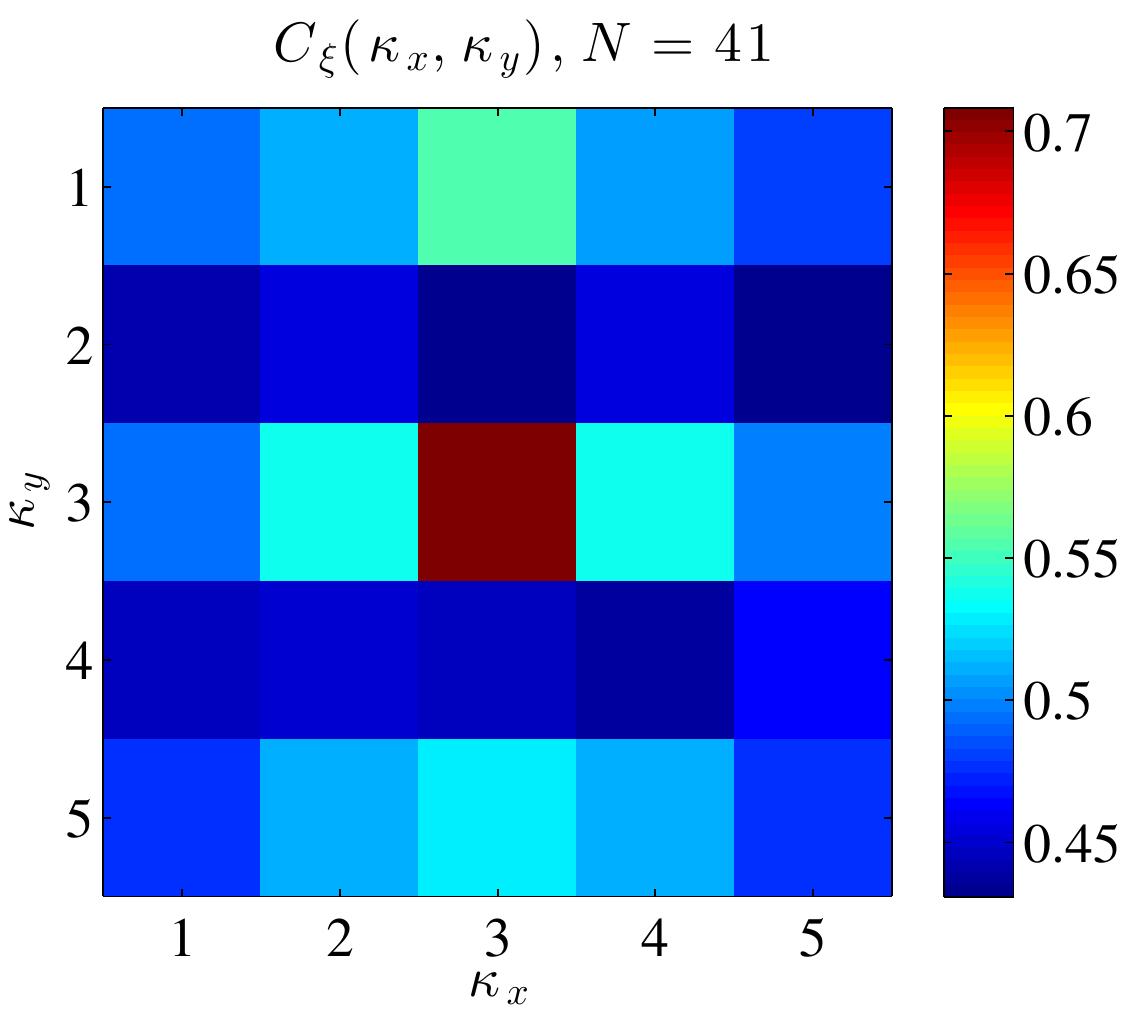}}
  \end{center}
  \caption{(a) \REV{Global error and the upper/lower bound estimator} for $V(x,y)=0.01$ and
  $f(x,y)=\cos(3x) \cos(y)$. (b) Local effectiveness of the upper bound
  characterized by $C_{\eta}$ in each element for $N=11$.
  (c) Local effectiveness of the upper bound
  characterized by $C_{\eta}$ in each element for $N=41$. 
  (d) Local effectiveness of the lower bound
  characterized by $C_{\xi}$ in each element for $N=11$.
  (e) Local effectiveness of the lower bound
  characterized by $C_{\xi}$ in each element for $N=41$. 
  }
  \label{fig:estPoisson2D}
\end{figure}

Fig.~\ref{fig:estPoisson2D} (a) shows the error in the energy norm, 
the computed upper bound and the lower bound as the number of ALBs per
element $N$ increases from $11$ to $41$.  Both the computed upper and
the lower bound estimates are effective for all calculations.
Fig.~\ref{fig:estPoisson2D} (b)-(d) illustrates the local effectiveness
of the upper and lower bound estimates for the two extreme cases
$N=11$ and $N=41$, and the estimator $\eta_{\kappa}$ and $\xi_{\kappa}$
are effective for all elements, and the effectiveness depends weakly on
the number of basis functions per element.

%

\FloatBarrier

\subsection{Indefinite operators}\label{subsec:indef}

\begin{figure}[h]
  \begin{center}
    \subfloat[(a)]{\includegraphics[width=0.35\textwidth]{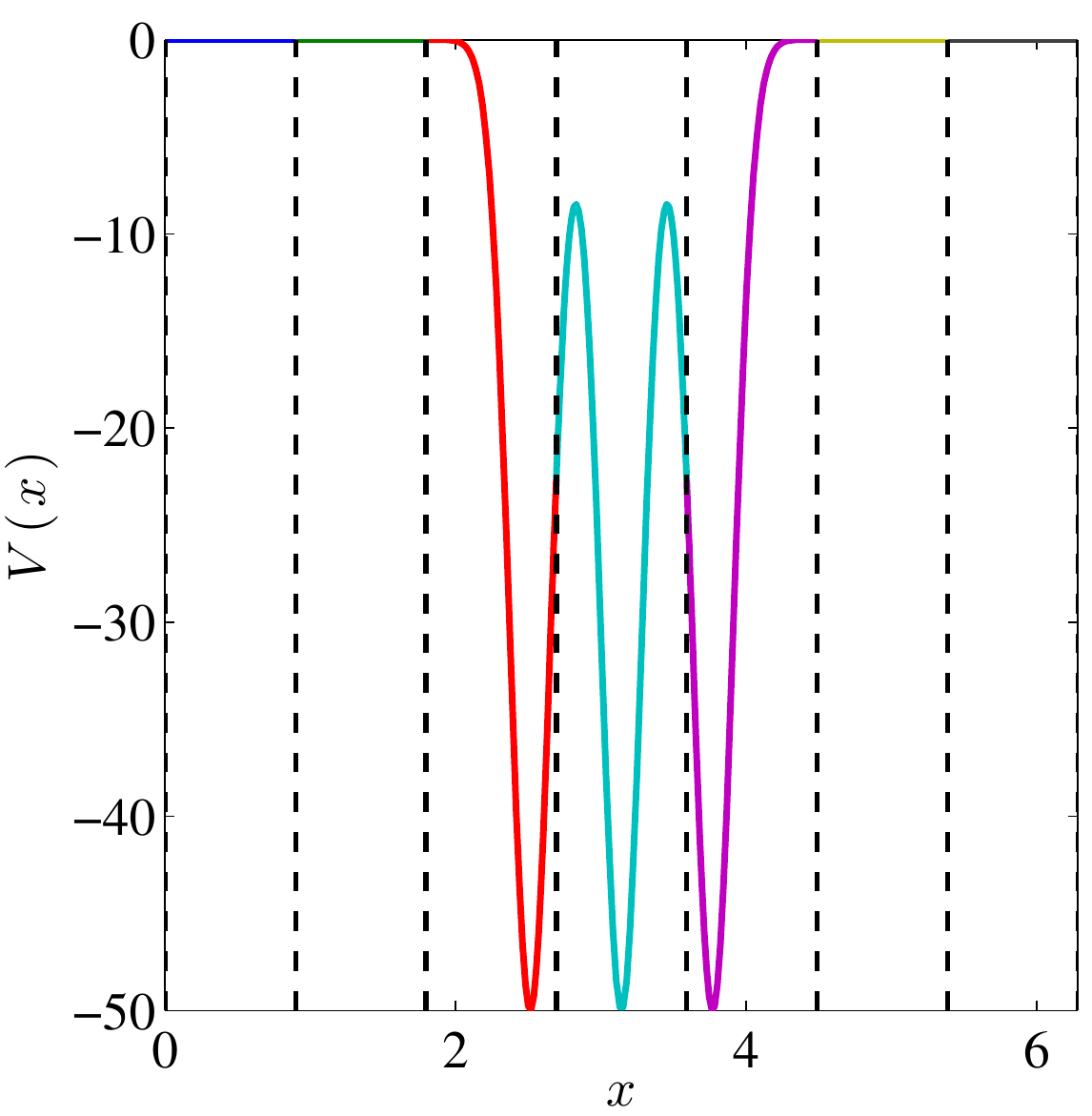}}
    \quad
    \subfloat[(b)]{\includegraphics[width=0.35\textwidth]{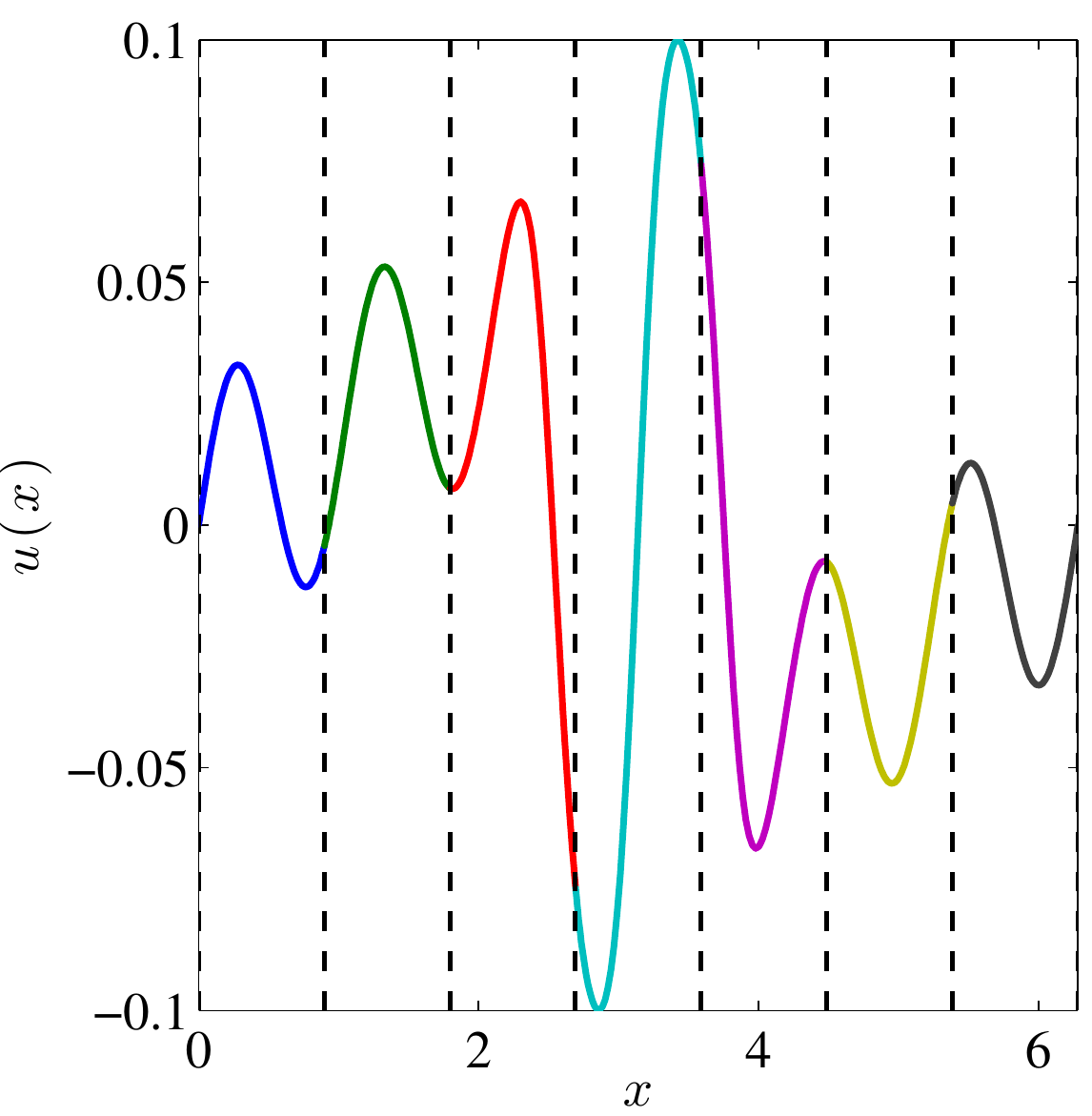}}
    \quad
    \subfloat[(c)]{\includegraphics[width=0.36\textwidth]{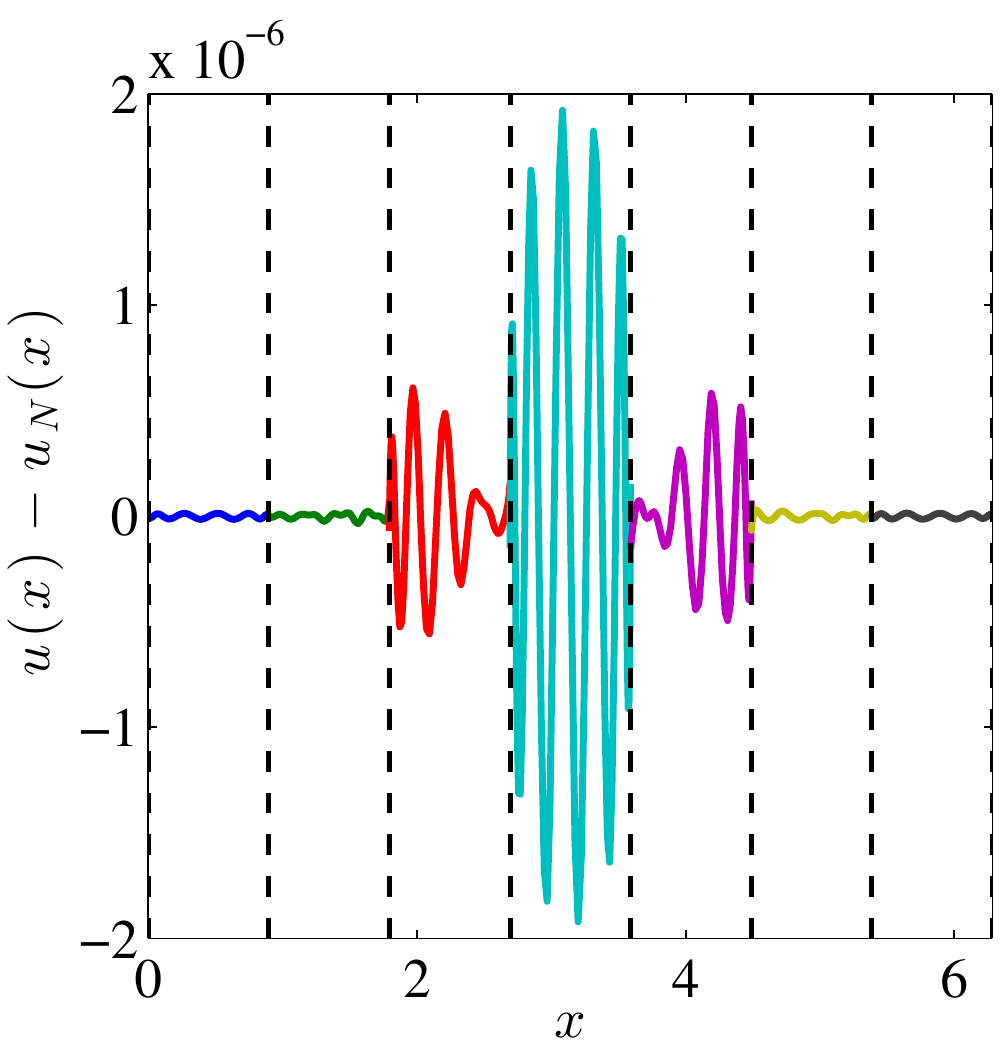}}
  \end{center}
  \caption{(a) The potential $V(x)$ given by the sum of three Gaussians
  with negative magnitude.  (b) The reference solution $u(x)$ corresponding to the
  potential $V(x)$ in (a) and the right hand side $f(x)=\sin(6x)$. (c)
  Point-wise error between the reference solution $u(x)$ and the numerical
  solution $u_N(x)$ calculated using the ALB set with $7$ elements and
  $N=11$ basis functions per element.}
  \label{fig:uuerrIndef1D}
\end{figure}

We now demonstrate the effectiveness of the upper and lower bound
estimates for indefinite operators.  We start from a 1D example on a
domain \REV{$\Omega=(0,2\pi)$} with periodic boundary conditions. The
potential function $V(x)$ is given by the sum
of three Gaussians with negative magnitude, as shown in
Fig.~\ref{fig:uuerrIndef1D} (a).  The operator $A=-\Delta+V$ has
$3$ negative eigenvalues and is indefinite. The right hand side is
$f(x)=\sin(6x)$. The domain is
partitioned into $7$ elements for the ALB calculation.
Fig.~\ref{fig:uuerrIndef1D} (b) shows the reference solution $u$ to
Eq.~\eqref{eq:Indef},
and Fig.~\ref{fig:uuerrIndef1D} (c) shows the point-wise error $u-u_{N}$
using $N=11$ ALBs per element. 

\begin{figure}[h]
  \begin{center}
    \includegraphics[width=0.33\textwidth]{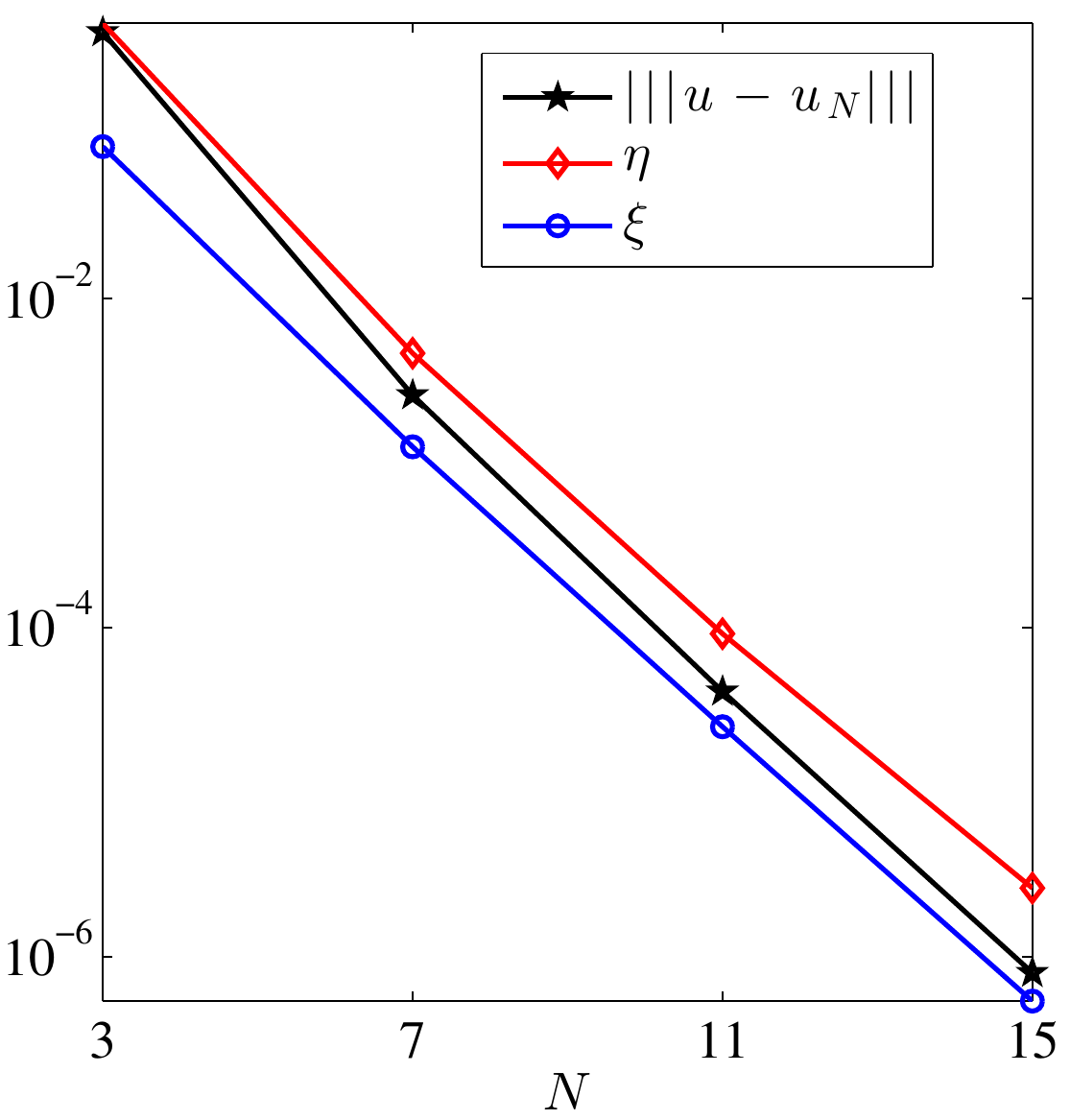}
    \includegraphics[width=0.34\textwidth]{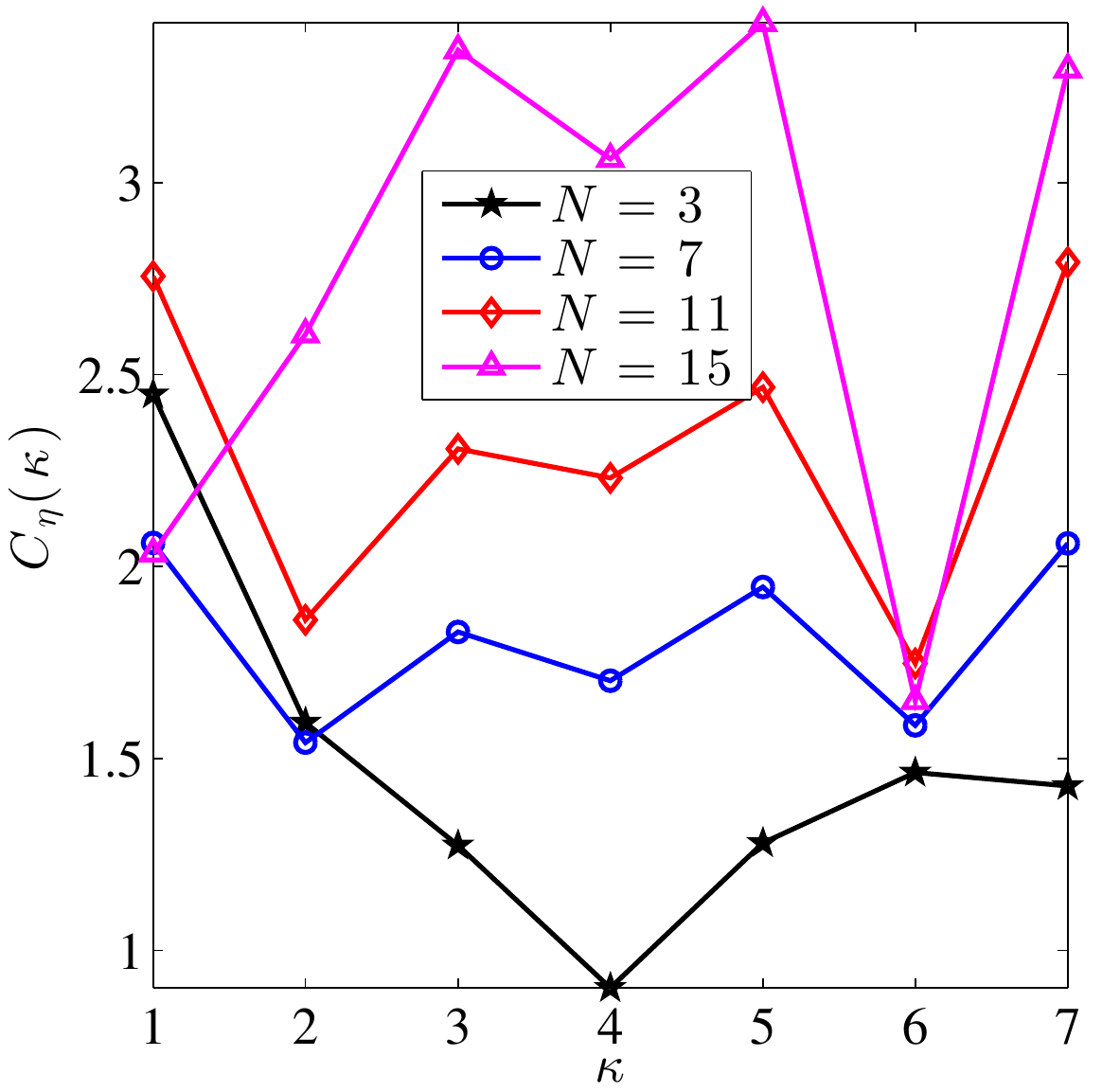}
    \includegraphics[width=0.34\textwidth]{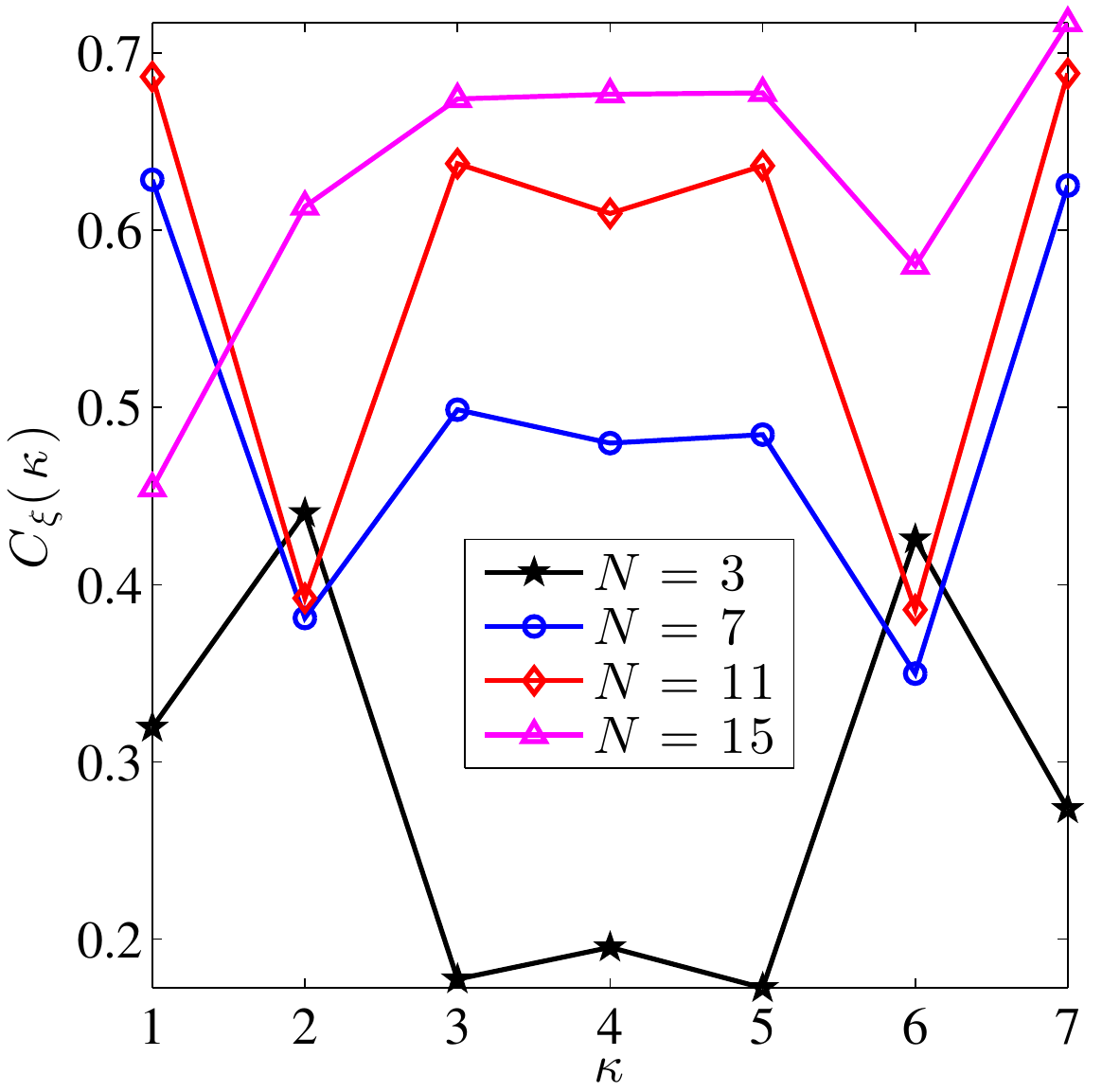}
  \end{center}
  \caption{(a) \REV{Global error and the upper/lower bound estimator} for $V(x)$ given in
  Fig.~\ref{fig:uuerrIndef1D} (a) and
  $f(x)=\sin(6x)$. (b) Local effectiveness of the upper bound
  characterized by $C_{\eta}$ in each element.  (c) Local effectiveness
  of the lower bound characterized by $C_{\eta}$ in each element. 
  }
  \label{fig:estIndef1D}
\end{figure}

Fig.~\ref{fig:estIndef1D} (a) shows the error in the energy norm, the
computed upper and lower bound estimates as the number of ALBs per
element $N$ increases from $3$ to $15$.  
Similar to
Fig.~\ref{fig:estPoisson1D}, the computed $\eta$ and
$\xi$ are upper and lower bounds for the true error $\tnorm{u-u_{N}}$
for all $N$ across a wide range of accuracy. Furthermore, the computed
$\xi$ is always a lower bound of $\tnorm{u-u_{N}}$ from $N=3$ to
$N=15$. This is guaranteed by the property of the lower bound in
Proposition~\ref{prop:GlobLowerBound}.

We should note that when the number of basis functions is very small
($N=3$), the accuracy is low and the ALB approximation is in its pre-asymptotic
regime. In such case, the upper bound is very close to the
true error.  In fact as indicated by Theorem~\ref{thm:IndefApost},
$\eta$ may not even be a rigorous upper bound for highly
indefinite operators with very few basis functions. 
%
%

%

\begin{figure}[h]
  \begin{center}
    \subfloat[(a)]{\includegraphics[width=0.4\textwidth]{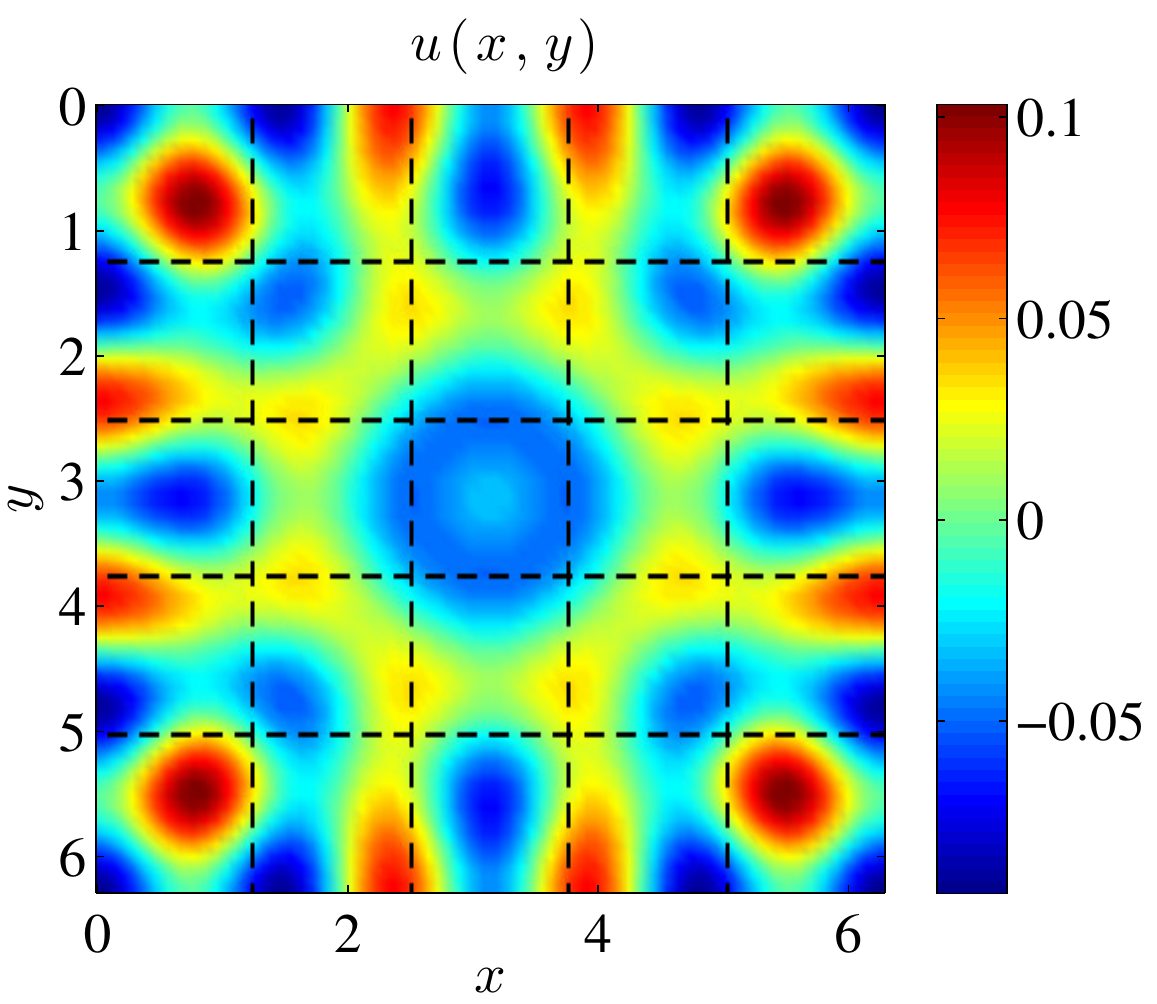}}
    \quad
    \subfloat[(b)]{\includegraphics[width=0.4\textwidth]{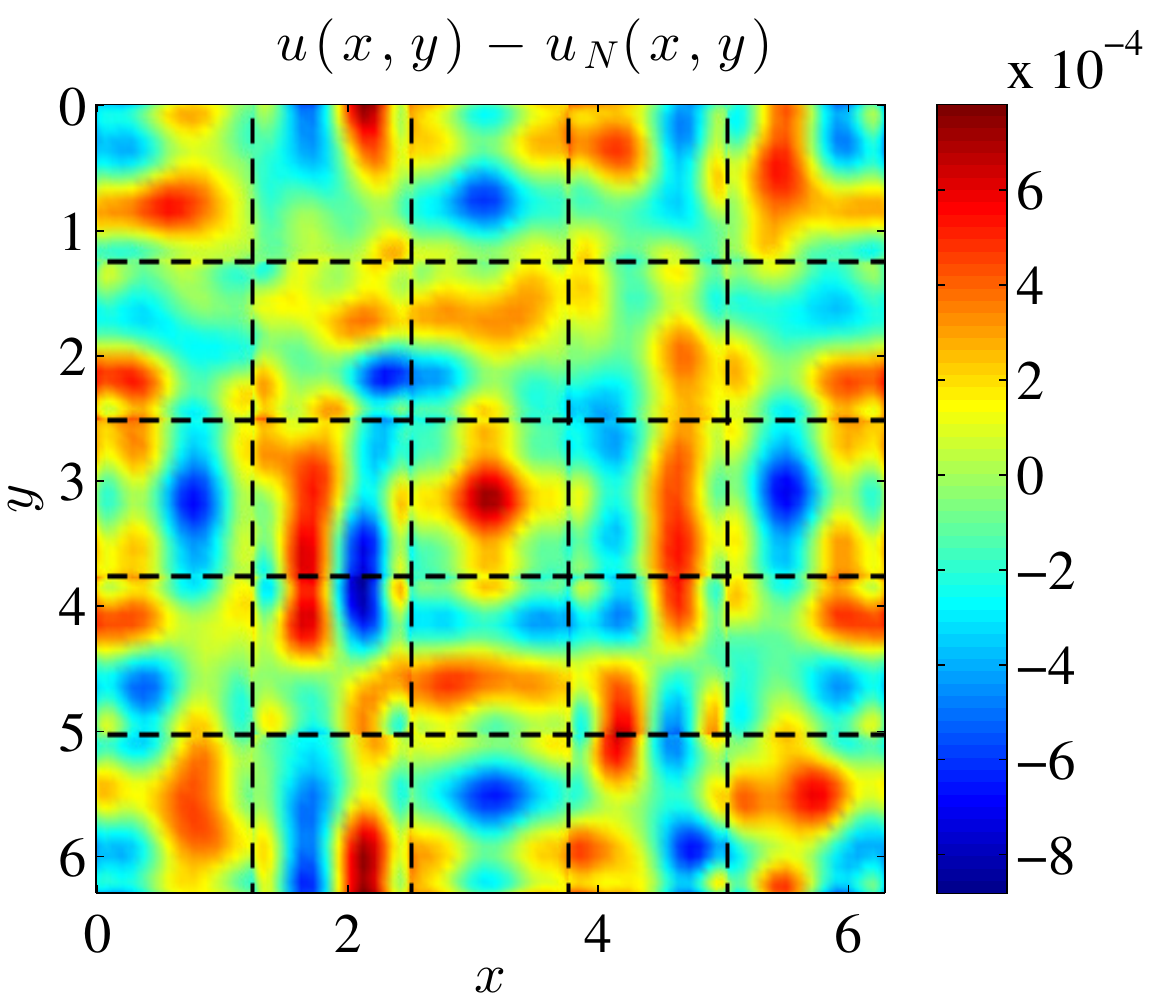}}
  \end{center}
  \caption{
  (a) The reference solution $u(x,y)$ corresponding to
  $V(x,y)=-16.5$ and $f(x,y)$ in Eq.~\eqref{eqn:rhsHelm2D}, which is a
  Gaussian localized at the center of $\Omega$. (b)
  Point-wise error between the reference solution $u(x,y)$ and the numerical
  solution $u_N(x,y)$ calculated using the ALB set with $5\times 5$
  elements and $N=31$ basis functions per element.}
  \label{fig:uuerrHelm2D}
\end{figure}

\begin{figure}[h]
  \begin{center}
    \subfloat[(a)]{\includegraphics[width=0.4\textwidth]{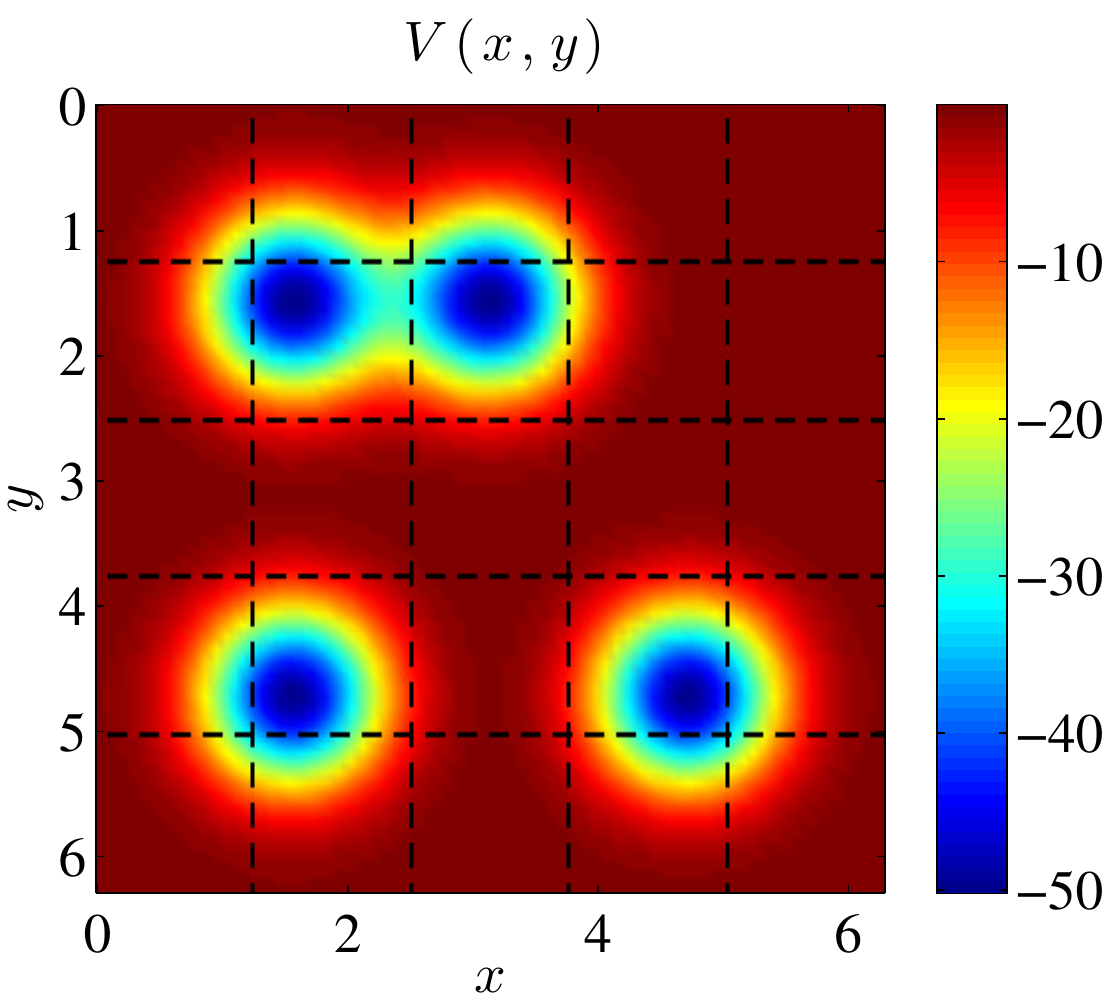}}
    \quad
    \subfloat[(b)]{\includegraphics[width=0.4\textwidth]{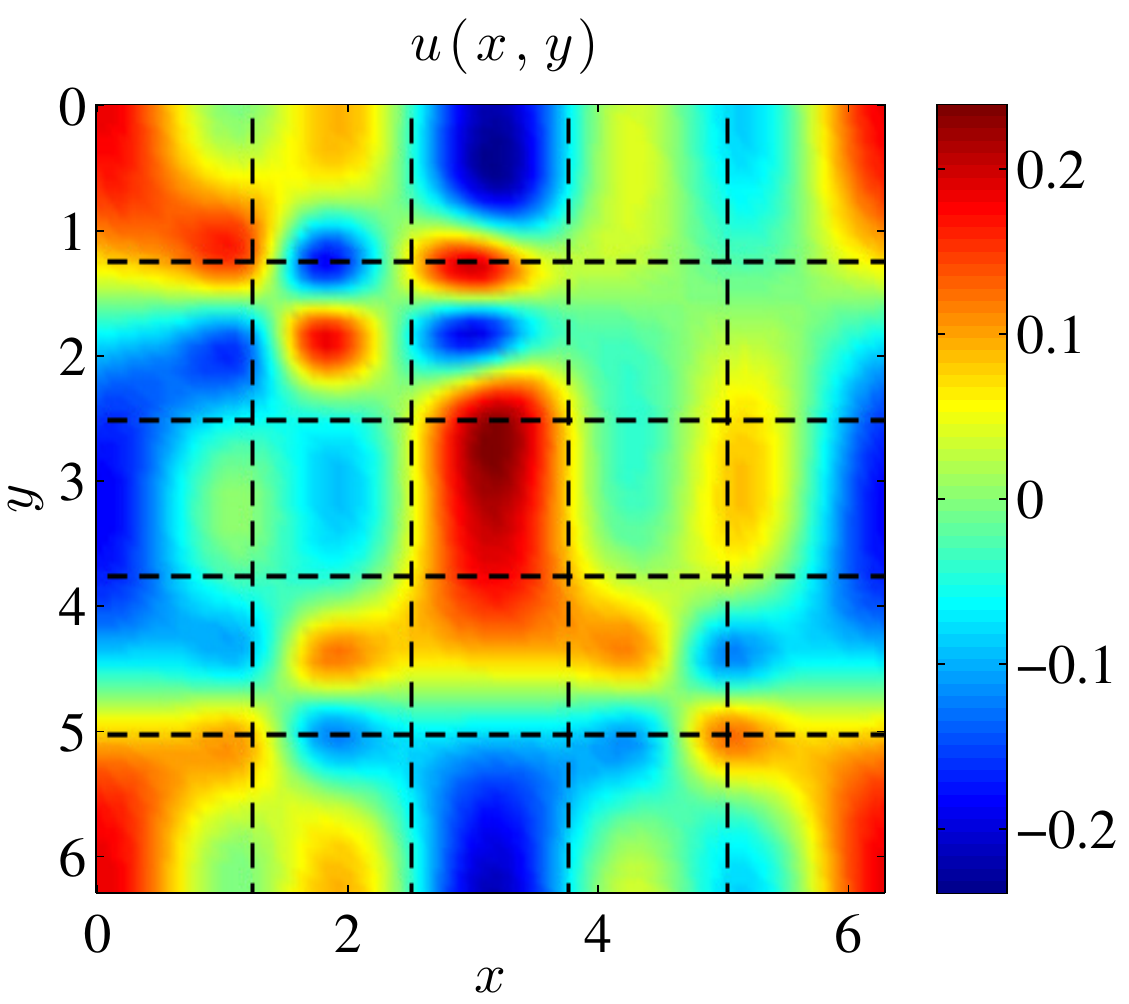}}
    \quad
    \subfloat[(c)]{\includegraphics[width=0.4\textwidth]{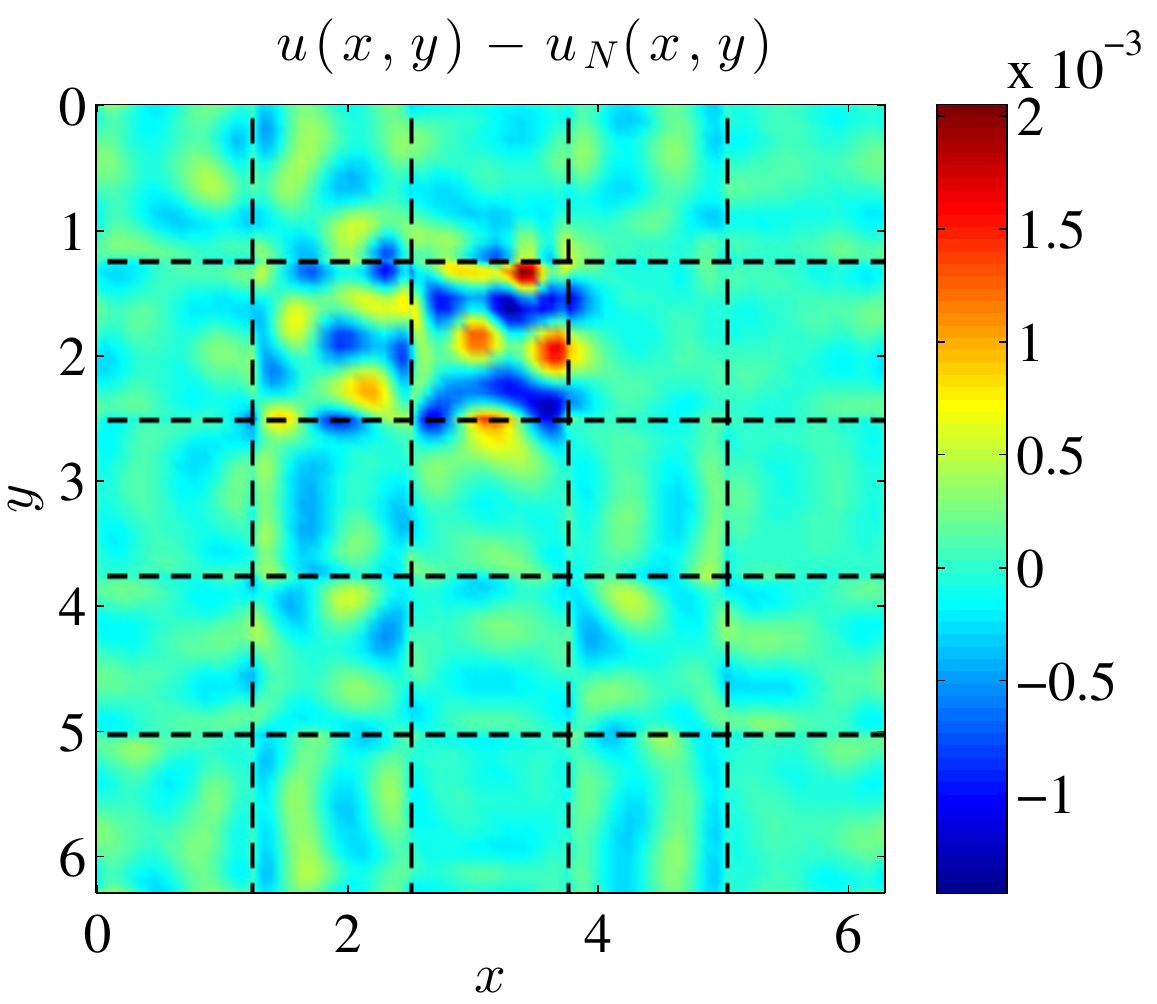}}
  \end{center}
  \caption{(a) The potential $V(x,y)$ four Gaussians with negative
  magnitude. (b) Solution $u(x,y)$ corresponding to $V(x,y)$ given in
  (a) and $f(x,y)=\cos(3x) \cos(y)$. (c) Point-wise error between the reference solution
  $u(x,y)$ and the numerical solution $u_N(x,y)$ calculated using the ALB
  set with $5\times 5$ elements and $N=31$ basis functions per element.}
  \label{fig:uuerrIndef2D}
\end{figure}

Our final examples are two indefinite problems on a 2D domain 
\REV{$\Omega=(0,2\pi)\times (0,2\pi)$}.  The first problem is a homogeneous Helmholtz
equation 
with $V(x,y)=-16.5$ and the operator $A=-\Delta+V$ has $49$ negative
eigenvalues. The right hand side is
\begin{equation}
  f(x,y)=\exp(-2(x-\pi)^2-2(y-\pi)^2),
  \label{eqn:rhsHelm2D}
\end{equation}
which is a Gaussian located at the center of $\Omega$.  
The second problem is that $V$ is given by the sum of four Gaussians
with negative magnitude, as illustrated in Fig.~\ref{fig:uuerrIndef2D}
(a). The operator $A=-\Delta+V$ has $26$ negative eigenvalues. The right
hand side is chosen to be $f(x,y)=\cos(3x) \cos(y)$ satisfying the
periodic boundary condition.  For the first problem,
Fig.~\ref{fig:uuerrHelm2D} (b) shows the reference solution $u$ to
Eq.~\eqref{eq:Indef} and Fig.~\ref{fig:uuerrHelm2D} (c) shows the
point-wise error $u-u_{N}$ using $N=31$ ALBs per element.   In the ALB
computation, the domain is partitioned into $5\times 5$ elements,
indicated by black dashed lines. Similarly for the second problem,
Fig.~\ref{fig:uuerrIndef2D} shows solution $u$ to Eq.~\eqref{eq:Indef}
and the point-wise error $u-u_{N}$ using $N=31$ ALBs per element. 



\begin{figure}[h]
  \begin{center}
    \subfloat[(a)]{\includegraphics[width=0.28\textwidth]{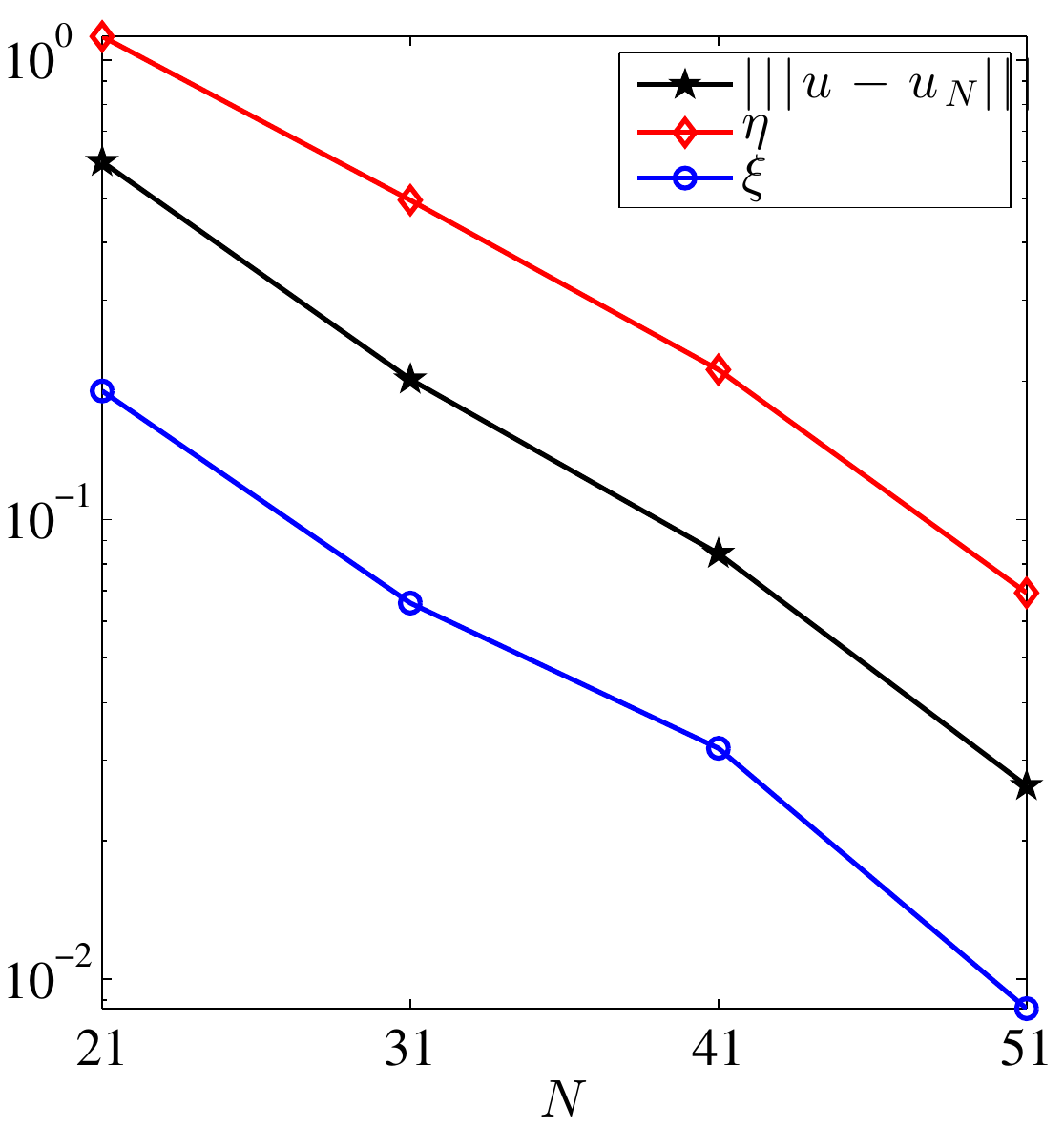}}
    \subfloat[(b)]{\includegraphics[width=0.35\textwidth]{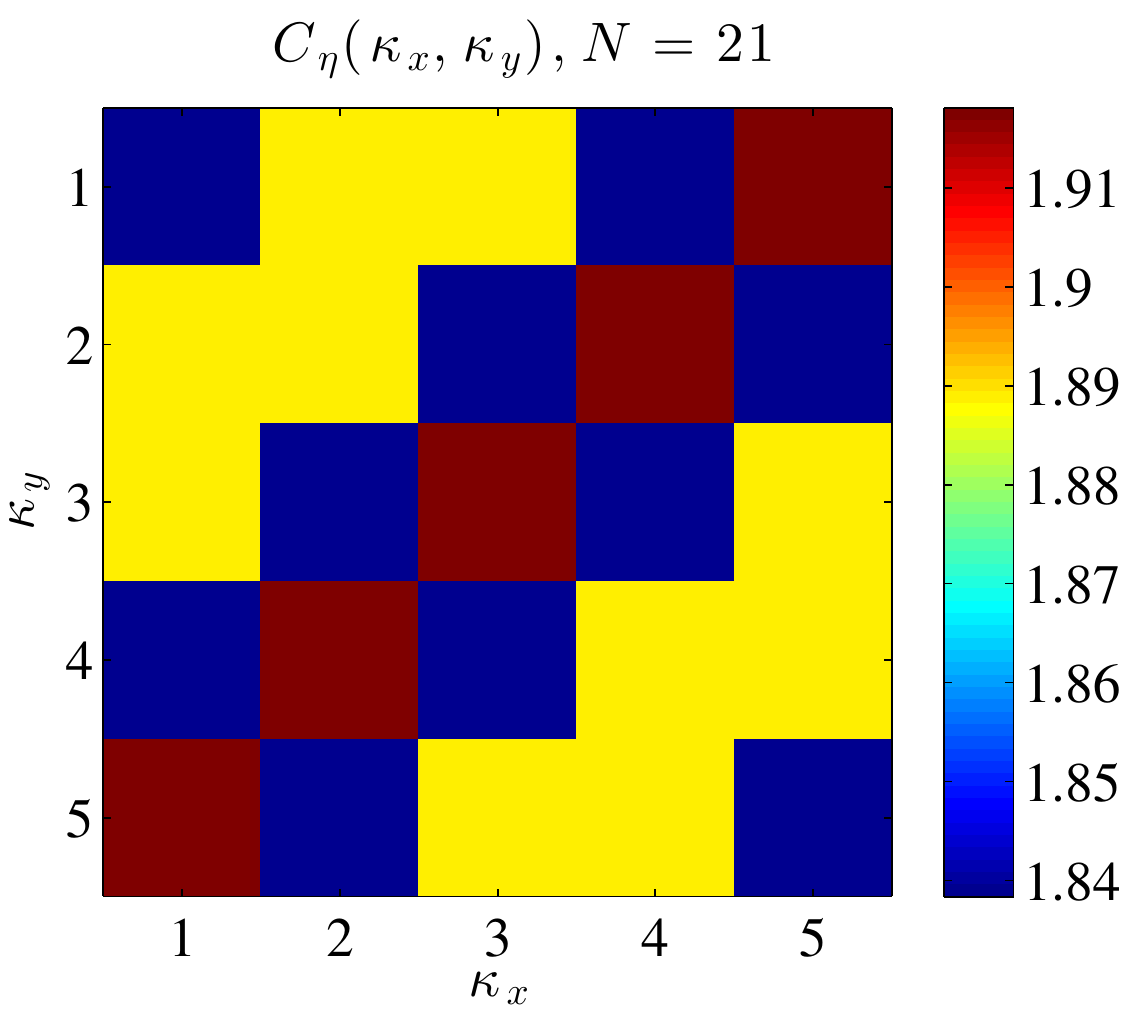}}
    \subfloat[(c)]{\includegraphics[width=0.35\textwidth]{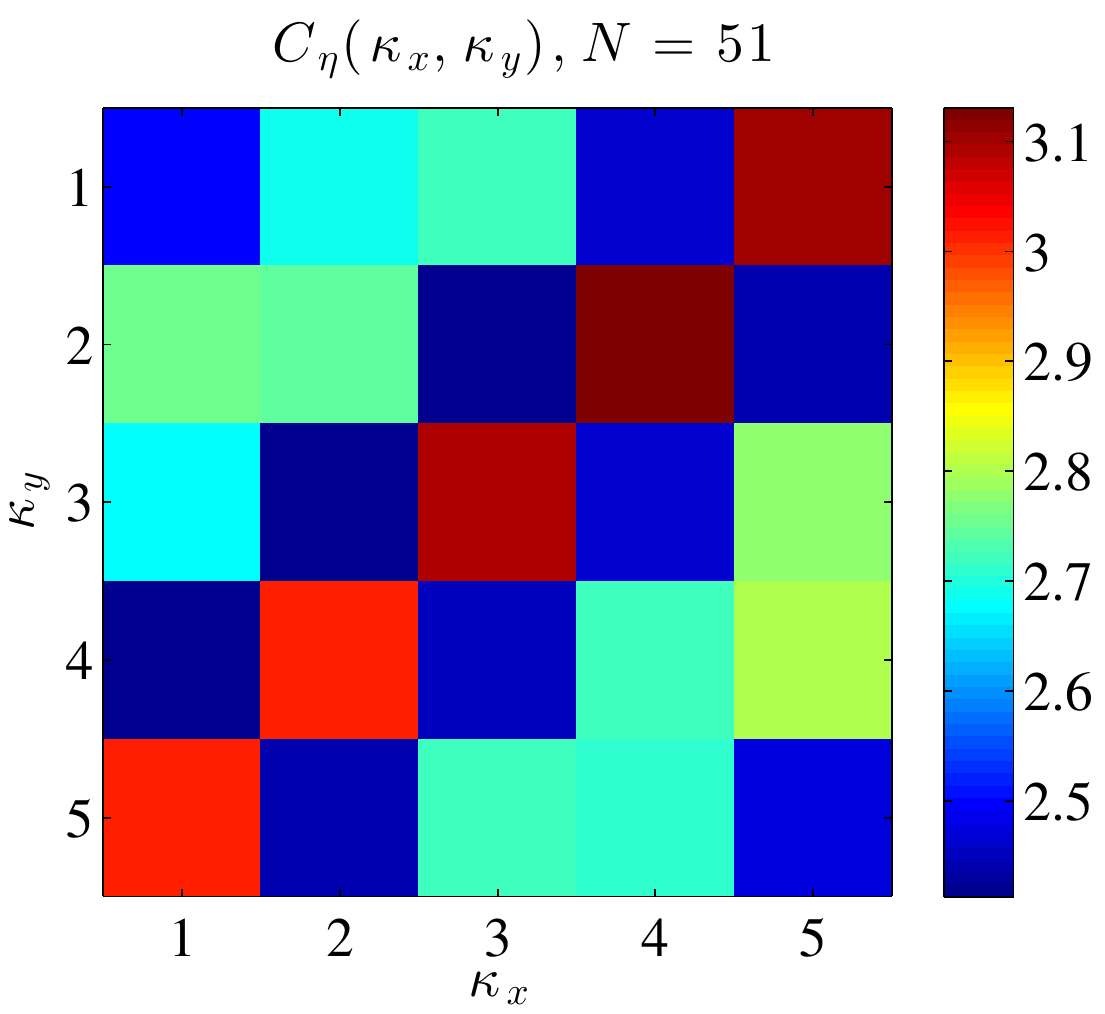}}

    \subfloat[(d)]{\includegraphics[width=0.35\textwidth]{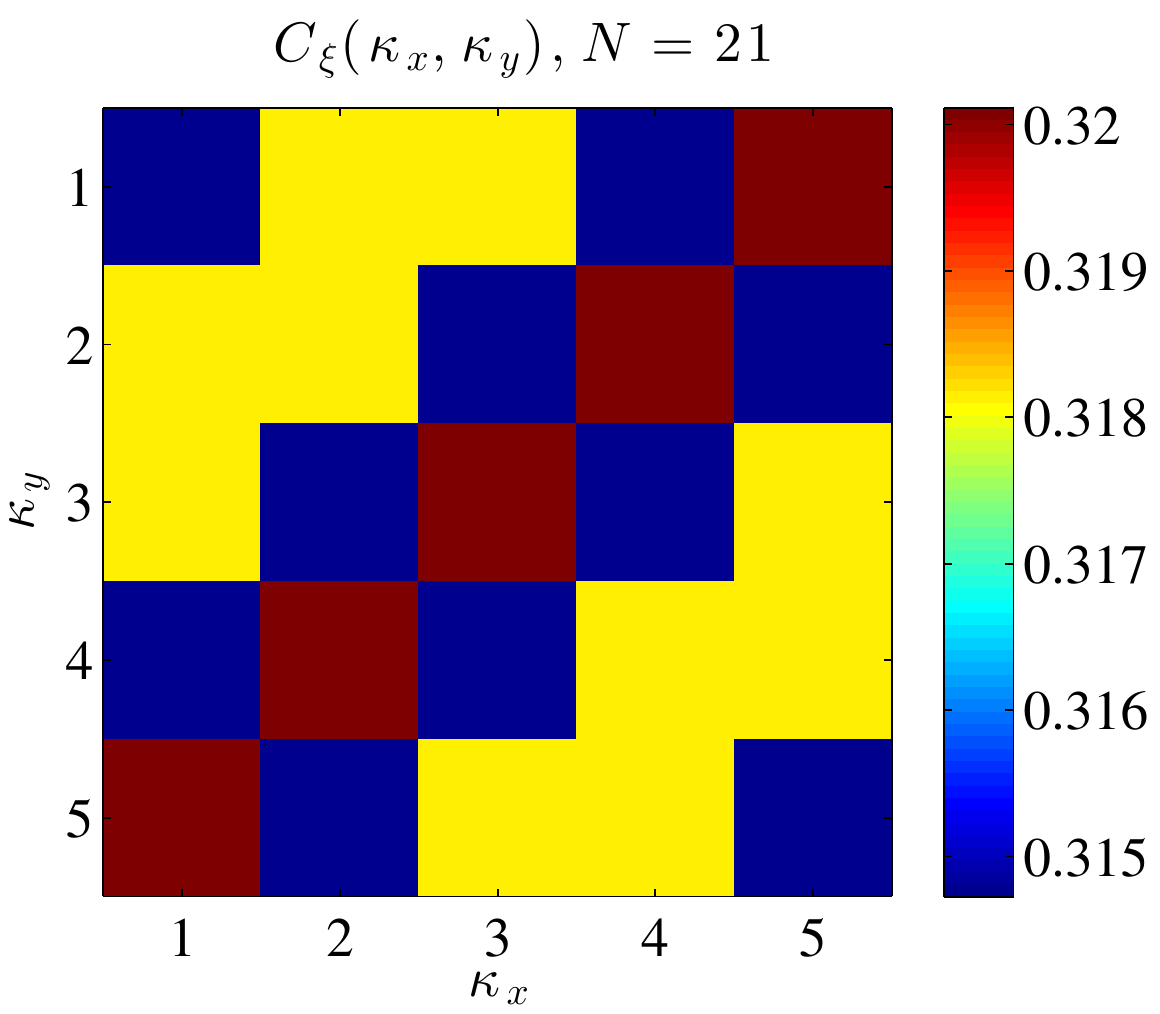}}
    \quad
    \subfloat[(e)]{\includegraphics[width=0.35\textwidth]{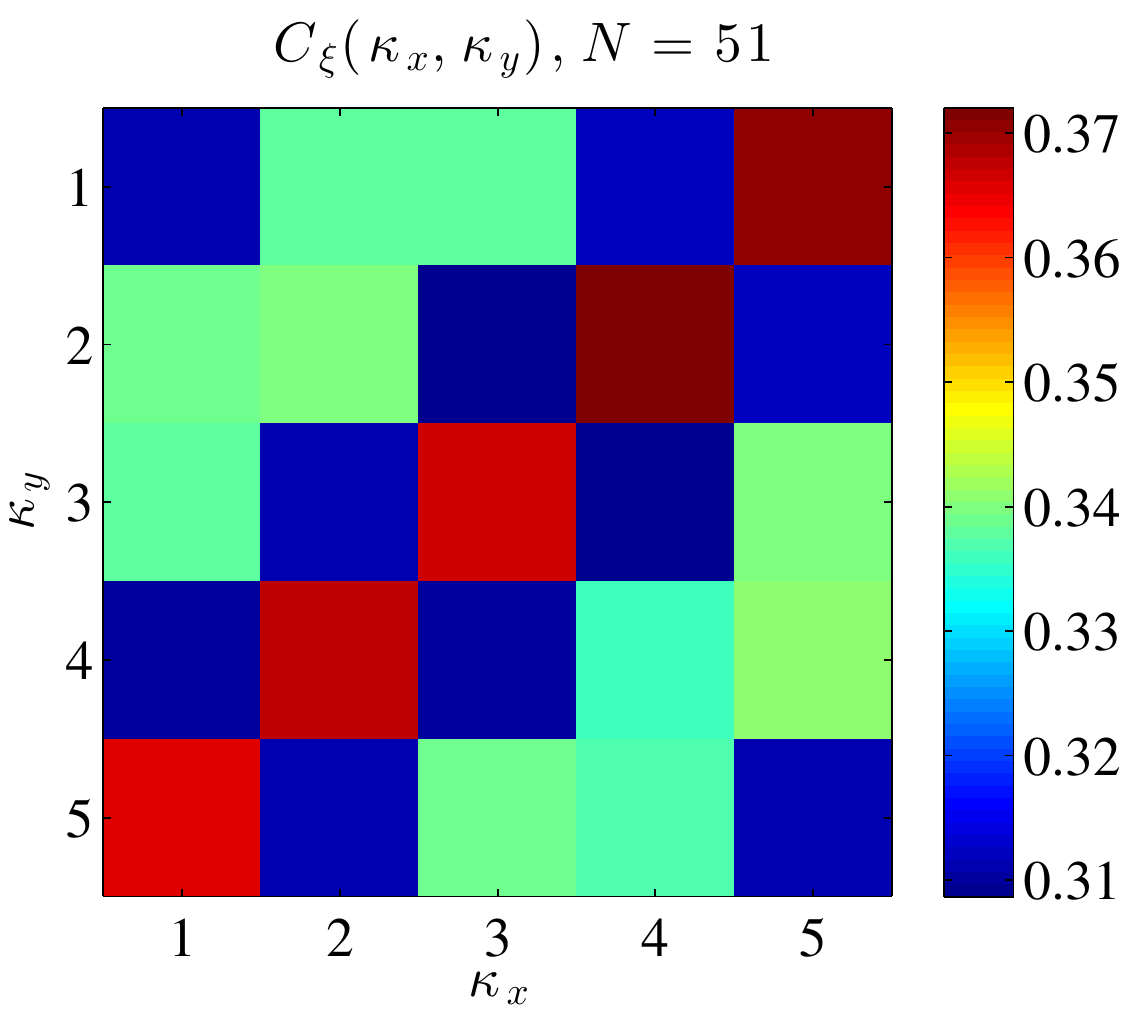}}
  \end{center}
  \caption{(a) \REV{Global error and the upper/lower bound estimator} for $V(x,y)=-16.5$ and
  $f(x,y)$ in Eq.~\eqref{eqn:rhsHelm2D}, which is a
  Gaussian localized at the center of $\Omega$. (b) Local effectiveness of the upper bound
  characterized by $C_{\eta}$ in each element for $N=21$.
  (c) Local effectiveness of the upper bound
  characterized by $C_{\eta}$ in each element for $N=51$. 
  (d) Local effectiveness of the lower bound
  characterized by $C_{\xi}$ in each element for $N=21$.
  (e) Local effectiveness of the lower bound
  characterized by $C_{\xi}$ in each element for $N=51$. 
 }
  \label{fig:estHelm2D}
\end{figure}

\begin{figure}[h]
  \begin{center}
    \subfloat[(a)]{\includegraphics[width=0.28\textwidth]{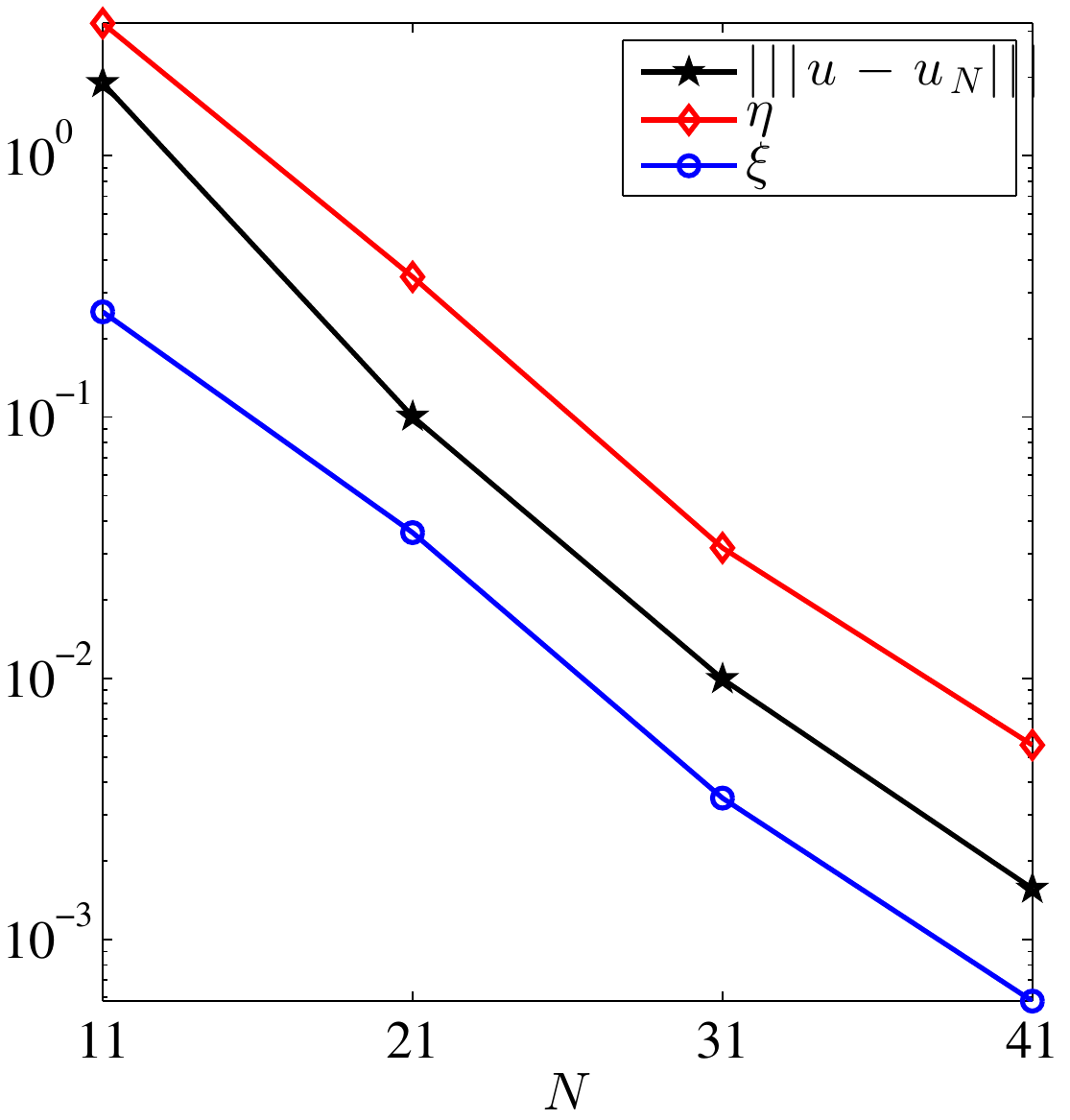}}
    \subfloat[(b)]{\includegraphics[width=0.34\textwidth]{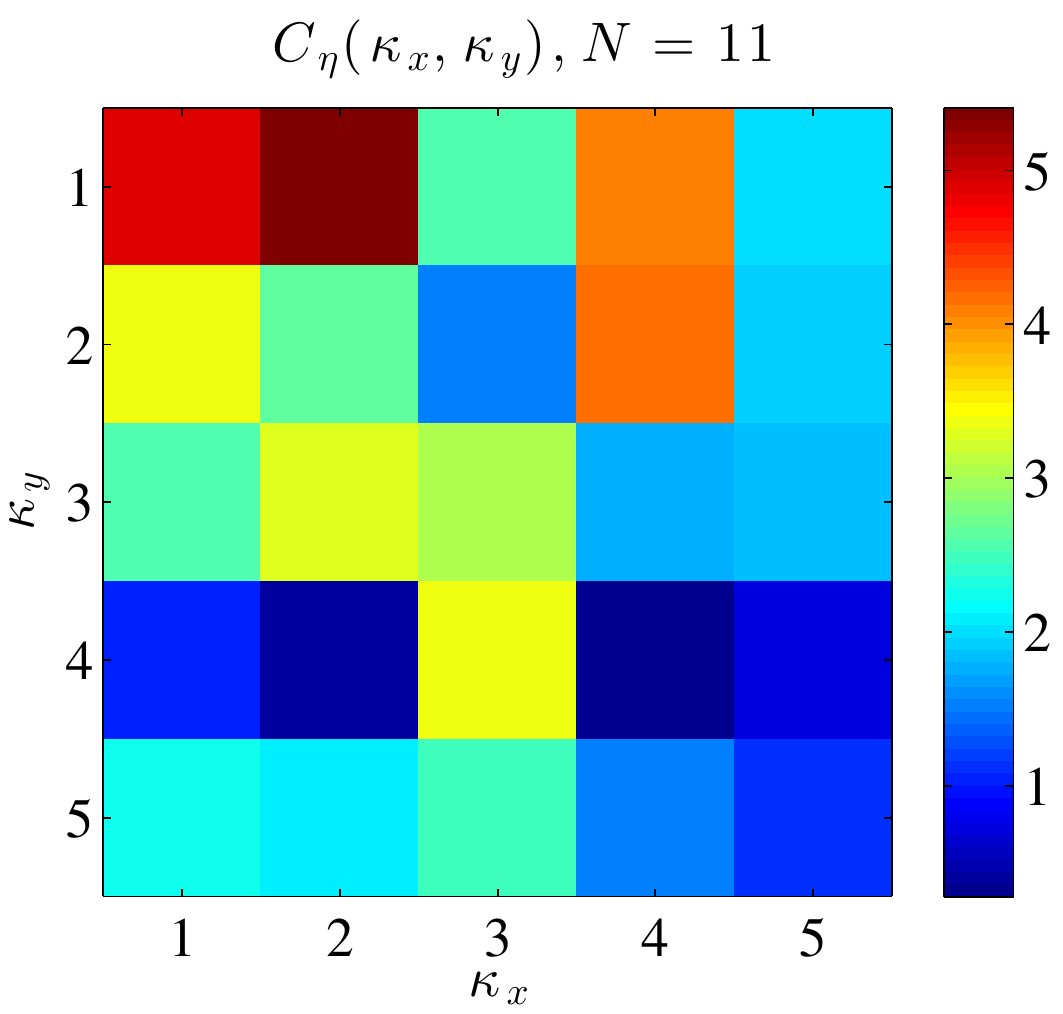}}
    \subfloat[(c)]{\includegraphics[width=0.35\textwidth]{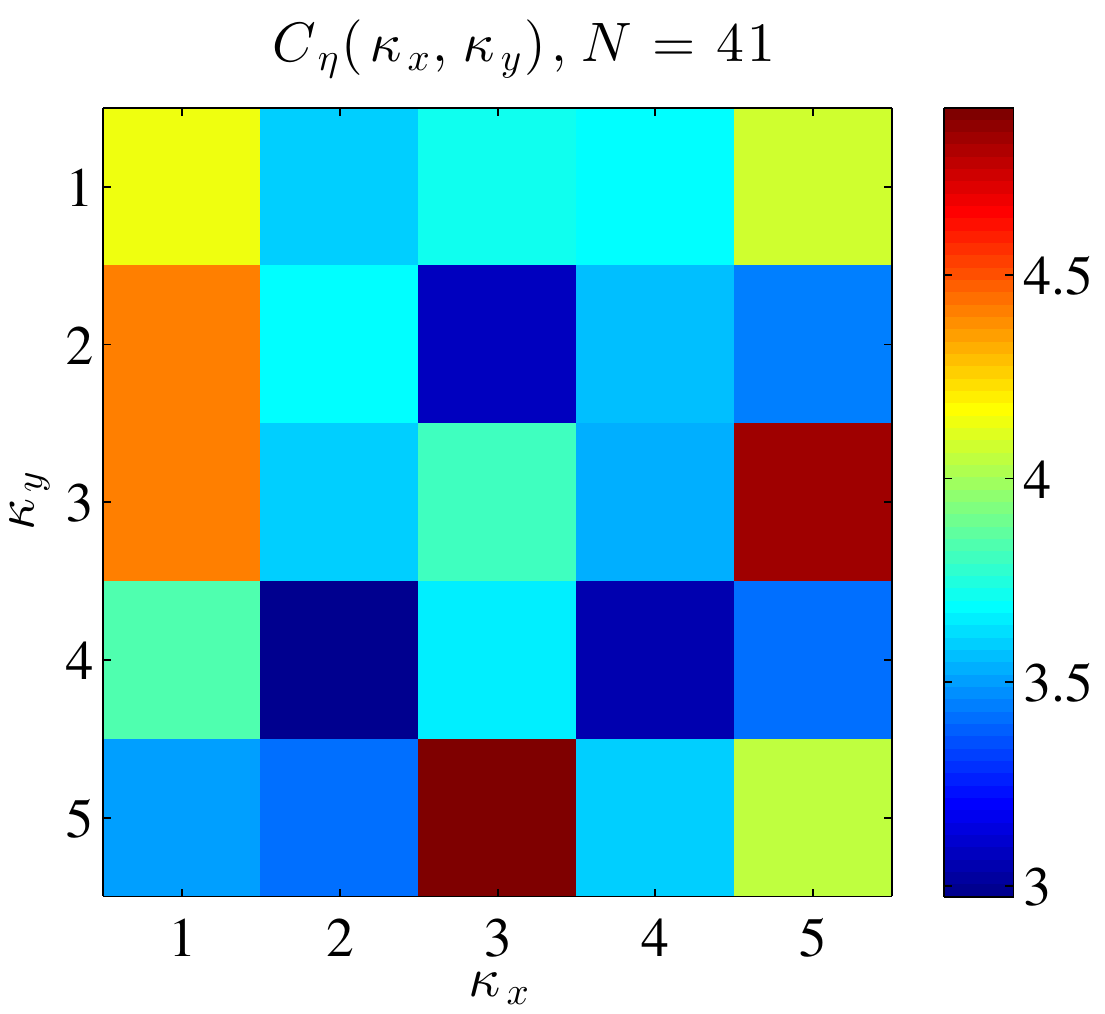}}
    \quad
    \subfloat[(d)]{\includegraphics[width=0.35\textwidth]{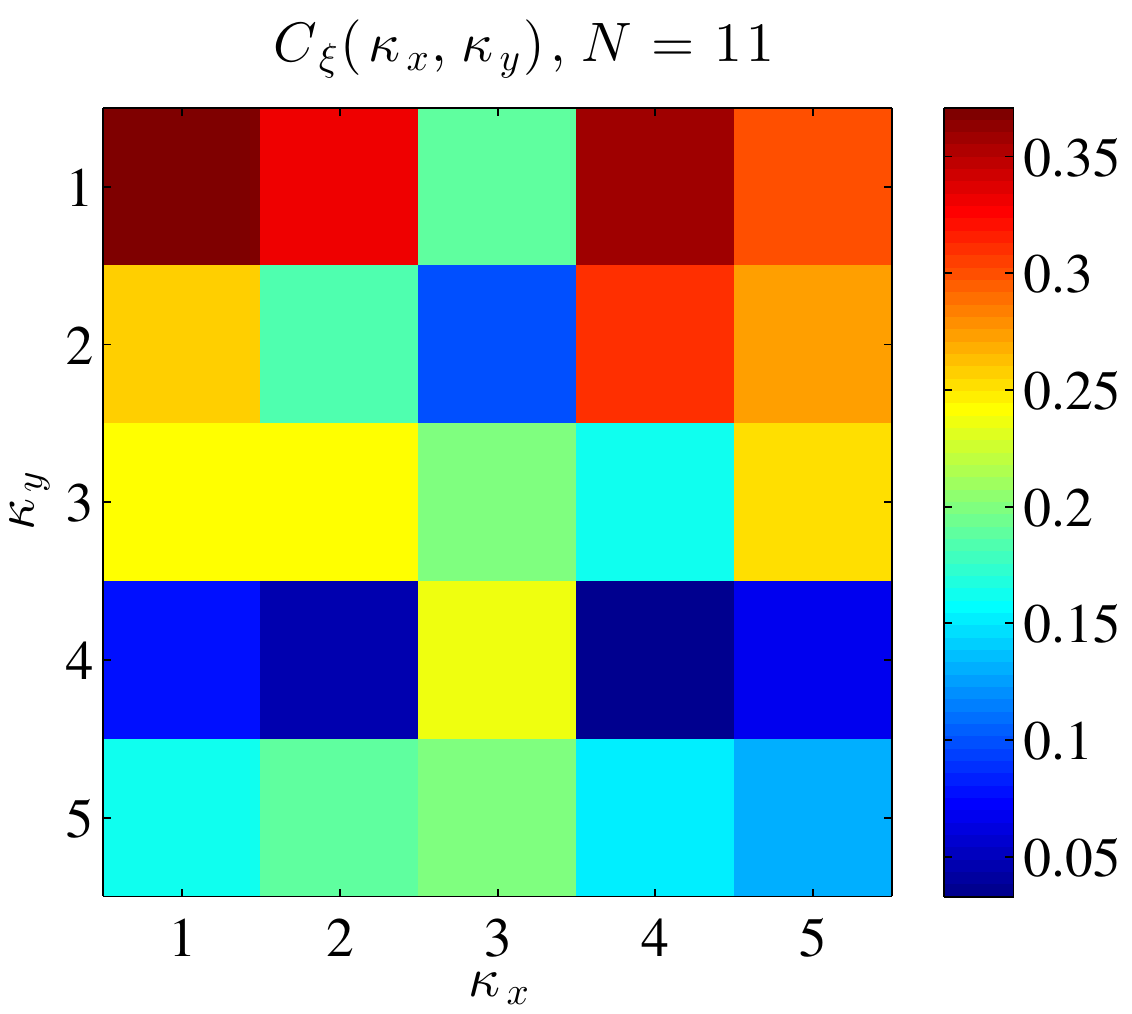}}
    \quad
    \subfloat[(e)]{\includegraphics[width=0.35\textwidth]{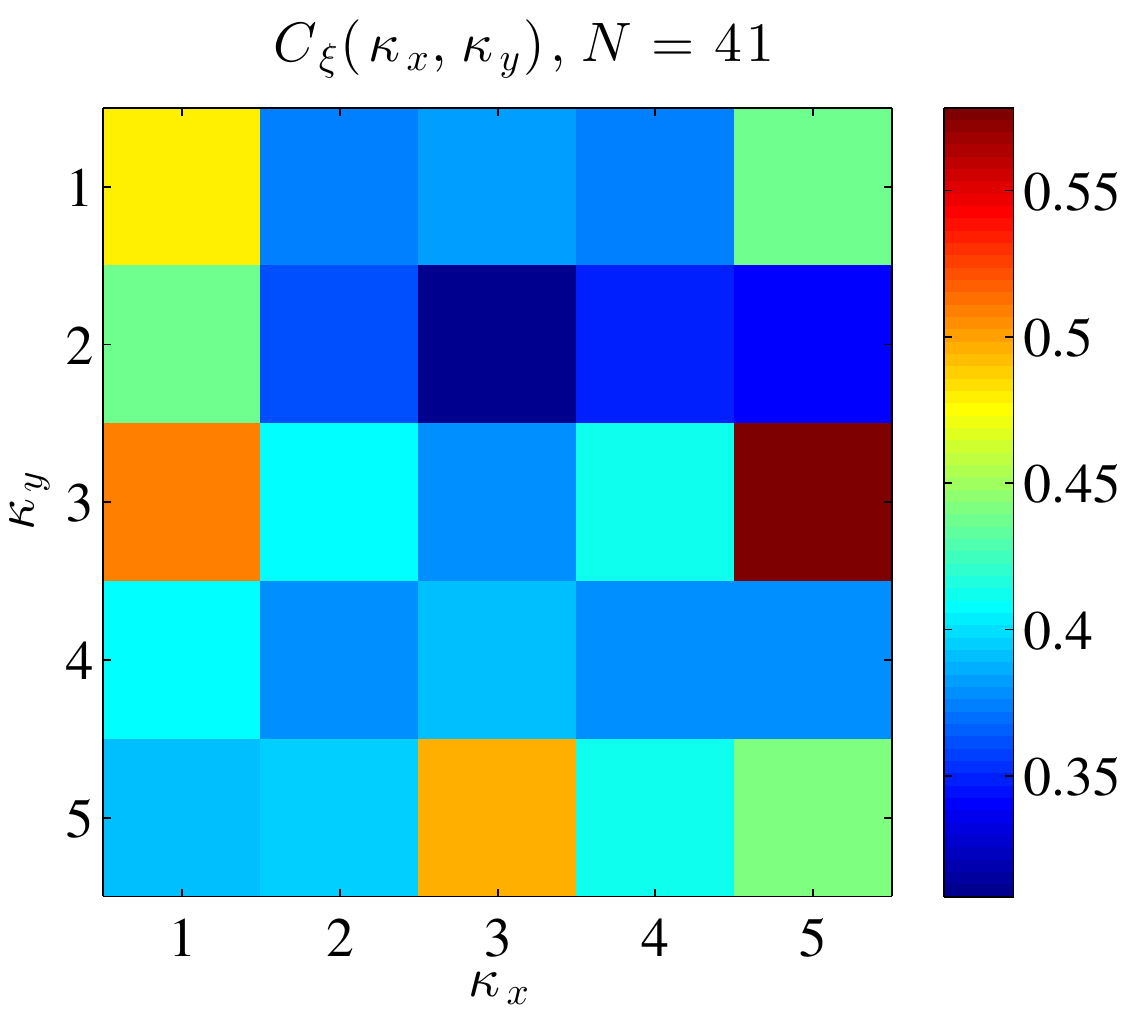}}
  \end{center}
  \caption{(a) \REV{Global error and the upper/lower bound estimator} for $V(x,y)$ given
  in Fig.~\ref{fig:uuerrIndef2D} (a) and
  $f(x,y)=\cos(3x) \cos(y)$. (b) Local effectiveness of the upper bound
  characterized by $C_{\eta}$ in each element for $N=11$.
  (c) Local effectiveness of the upper bound
  characterized by $C_{\eta}$ in each element for $N=41$. 
  (d) Local effectiveness of the lower bound
  characterized by $C_{\xi}$ in each element for $N=11$.
  (e) Local effectiveness of the lower bound
  characterized by $C_{\xi}$ in each element for $N=41$. 
 }
  \label{fig:estIndef2D}
\end{figure}

Fig.~\ref{fig:estHelm2D} (a)-(e) illustrates the global and local
effectiveness of the upper and lower bound estimates for the Helmholtz
problem, as the number of ALBs per element $N$ increases from $21$ to
$51$. Compared to the positive definite case in
Fig.~\ref{fig:estPoisson2D}, the true error is larger using a comparable
number of basis functions, reflecting that the Helmholtz equation is
more difficult to solve.  Nonetheless, $\eta$ and $\xi$ provide
effective bounds for the true error in all cases.  
Similar results can be found for the
indefinite example with negative Gaussian potentials in
Fig.~\ref{fig:estIndef2D} (a)-(e).  
In all calculations, the computed lower
bound estimator remains a lower bound for the true error. In particular, the
estimators still hold quite tightly in the pre-asymptotic regime
$(N=11)$ where the ALB approximation is crude and has large numerical
error.



%
%

\FloatBarrier

\subsection{Justification of the treatment of \REV{$\dku(\uN)$}}
\label{subsec:justifydku}

In the numerical computation of the upper and lower bound estimates, \REV{we approximated the non-computable constant
$\dku(\uN)$} by the computable constant $\dkN$.  
Below we provide numerical justification of such
approximation by direct computation of \REV{$\dku(\uN)$} via the reference
solution. We compare with $\dkN$ and $\rbk\gk$ since these three terms
appear together in $\etaJ$ in Eq.~\eqref{eq:Estim3}. 

\begin{figure}[h]
  \begin{center}
    \subfloat[(a)]{\includegraphics[width=0.37\textwidth]{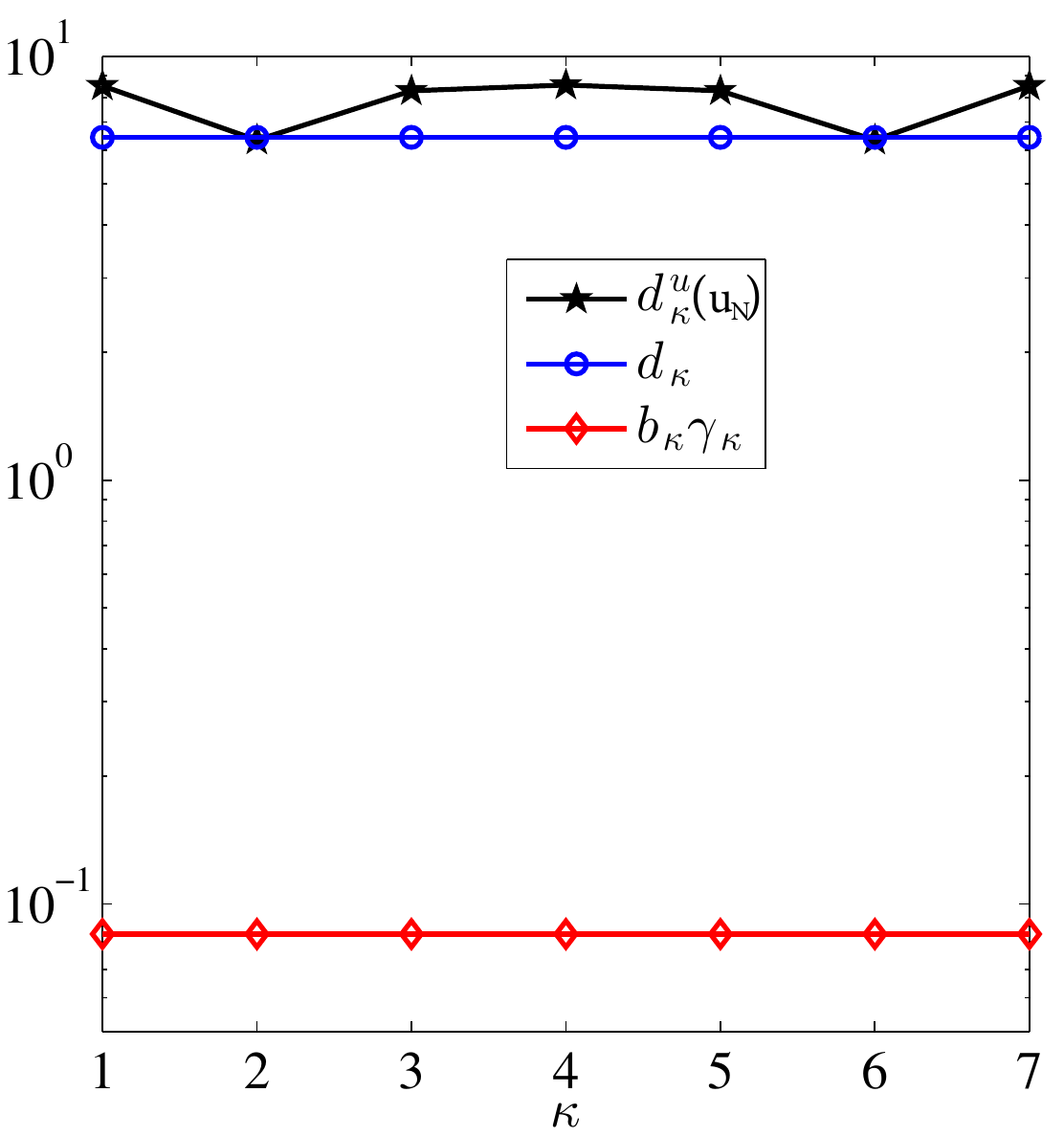}}
    \quad
    \subfloat[(b)]{\includegraphics[width=0.37\textwidth]{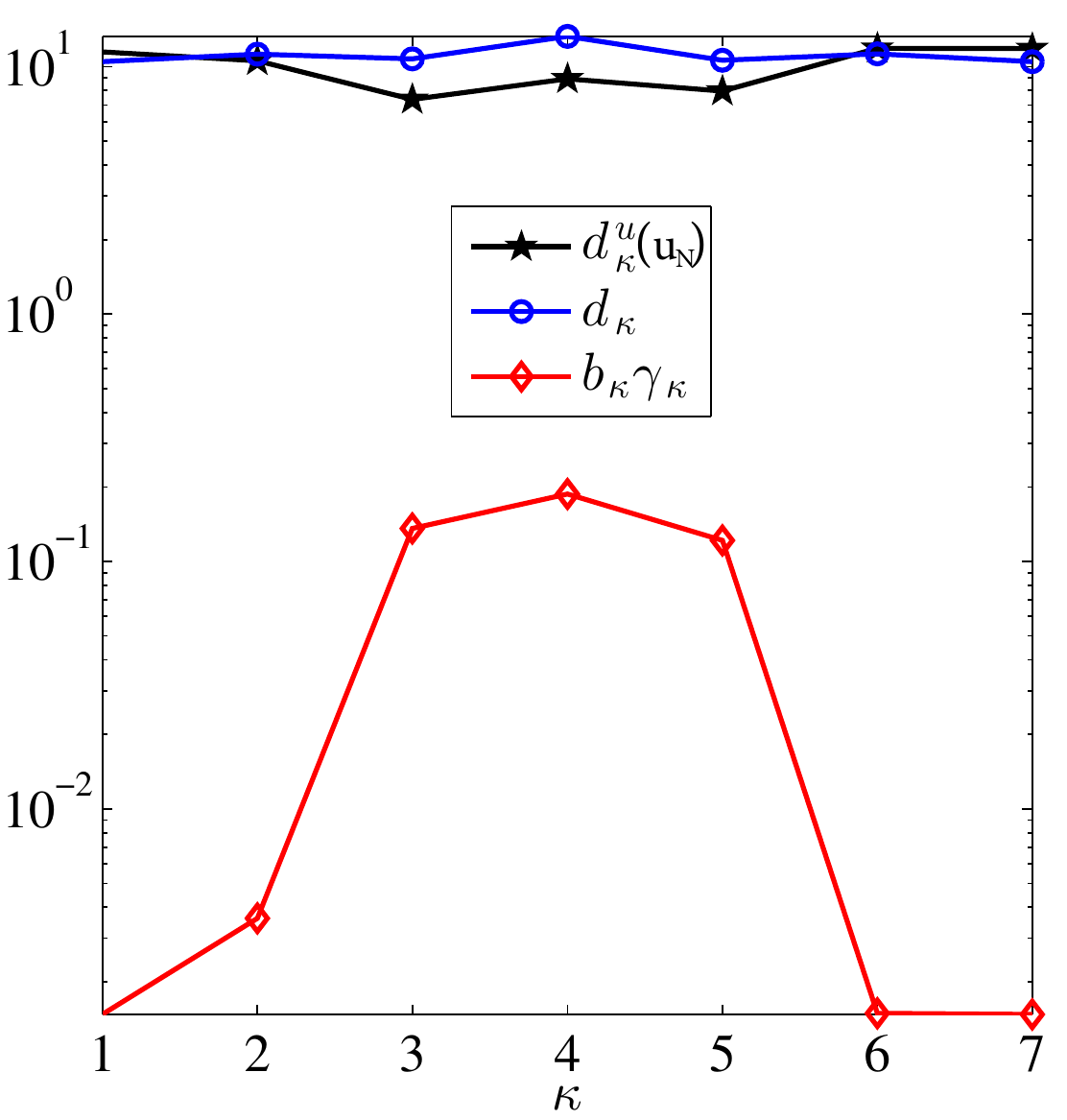}}
  \end{center}
  \caption{Comparison of \REV{$\dku(\uN)$}, $\dkN$ and $\rbk \gk$ for (a) the positive
  definite case with $V(x)=0.01$ with $N=7$. (b) the indefinite case with $V(x)$
  given in Fig.~\ref{fig:uuerrIndef1D} (a) with $N=7$.}
  \label{fig:dkuComp1D}
\end{figure}

Fig.~\ref{fig:dkuComp1D} (a) and (b) compare \REV{$\dku(\uN)$}, $\dkN$ and $\rbk
\gk$ for the positive definite and the indefinite 1D examples, 
respectively. We observe that the magnitude of \REV{$\dku(\uN)$} is comparable
to that of $\dkN$.  $\rbk\gk$ is much smaller compared to \REV{$\dku(\uN)$} and
$\dkN$. This is a direct consequence of Proposition~\ref{prop:rbk1D},
which states that $\rbk$ is in general very small for 1D systems.

\begin{figure}[h]
  \begin{center}
    \subfloat[(a)$V=0.01$]{\includegraphics[width=0.30\textwidth]{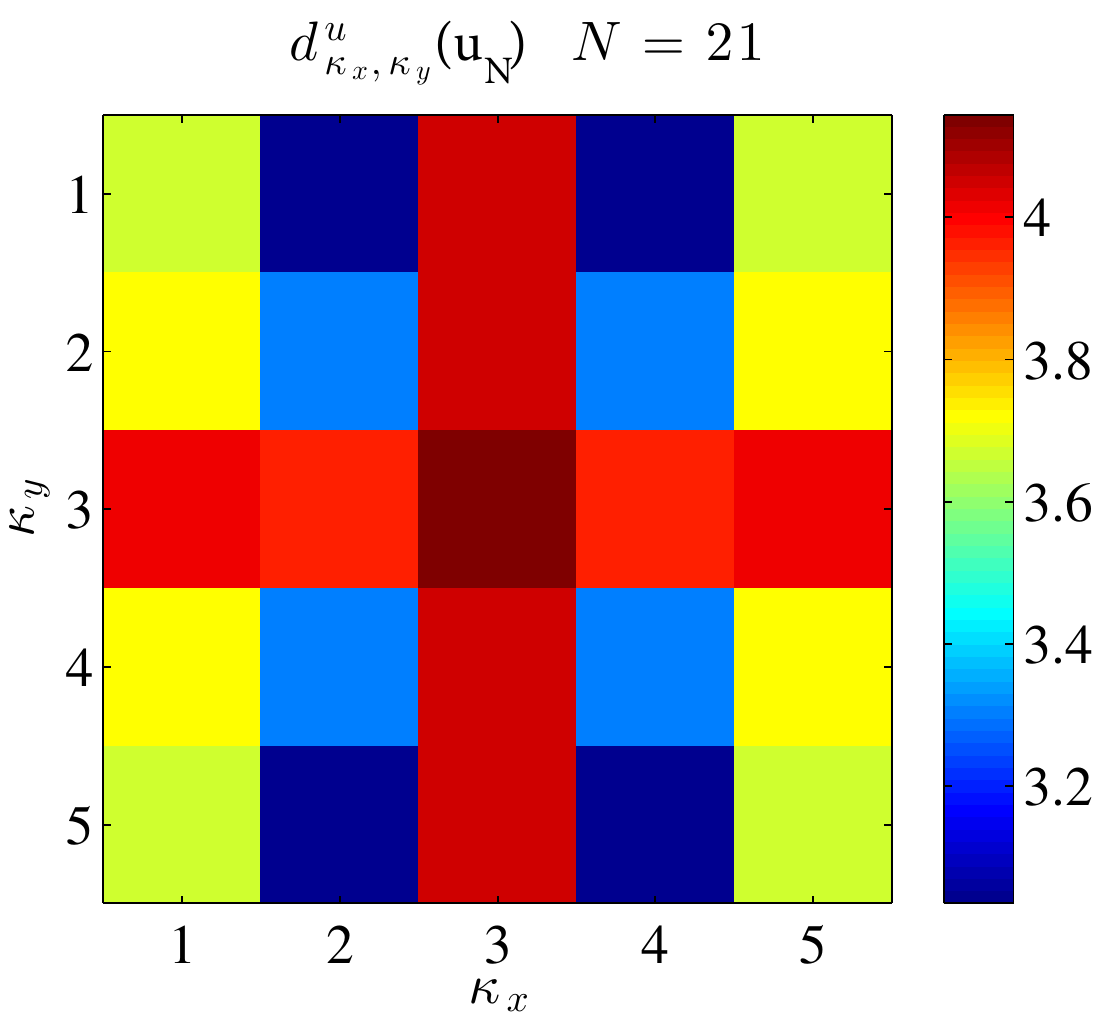}}
    \subfloat[(b)$V=0.01$]{\includegraphics[width=0.32\textwidth]{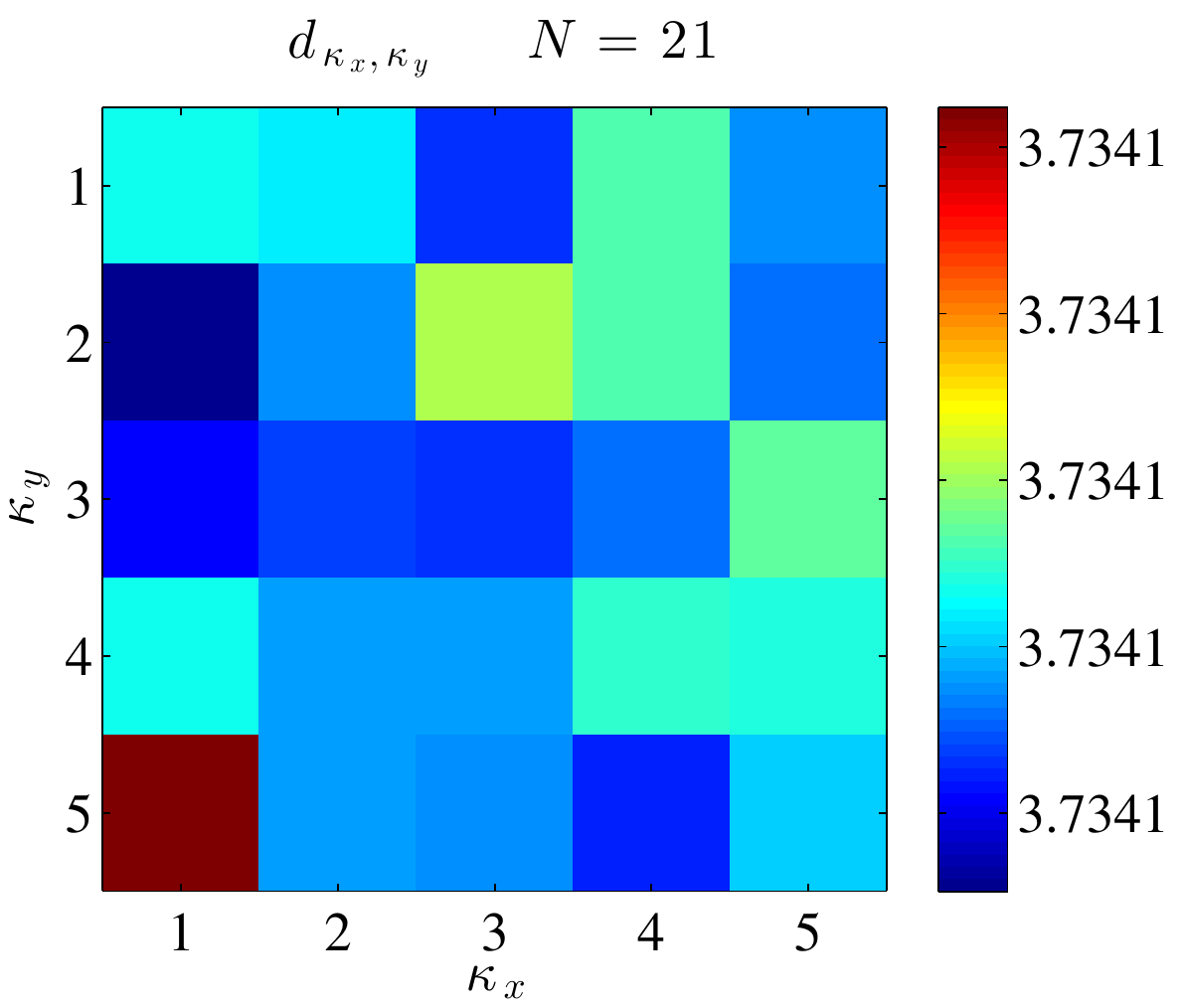}}
    \subfloat[(c)$V=0.01$]{\includegraphics[width=0.32\textwidth]{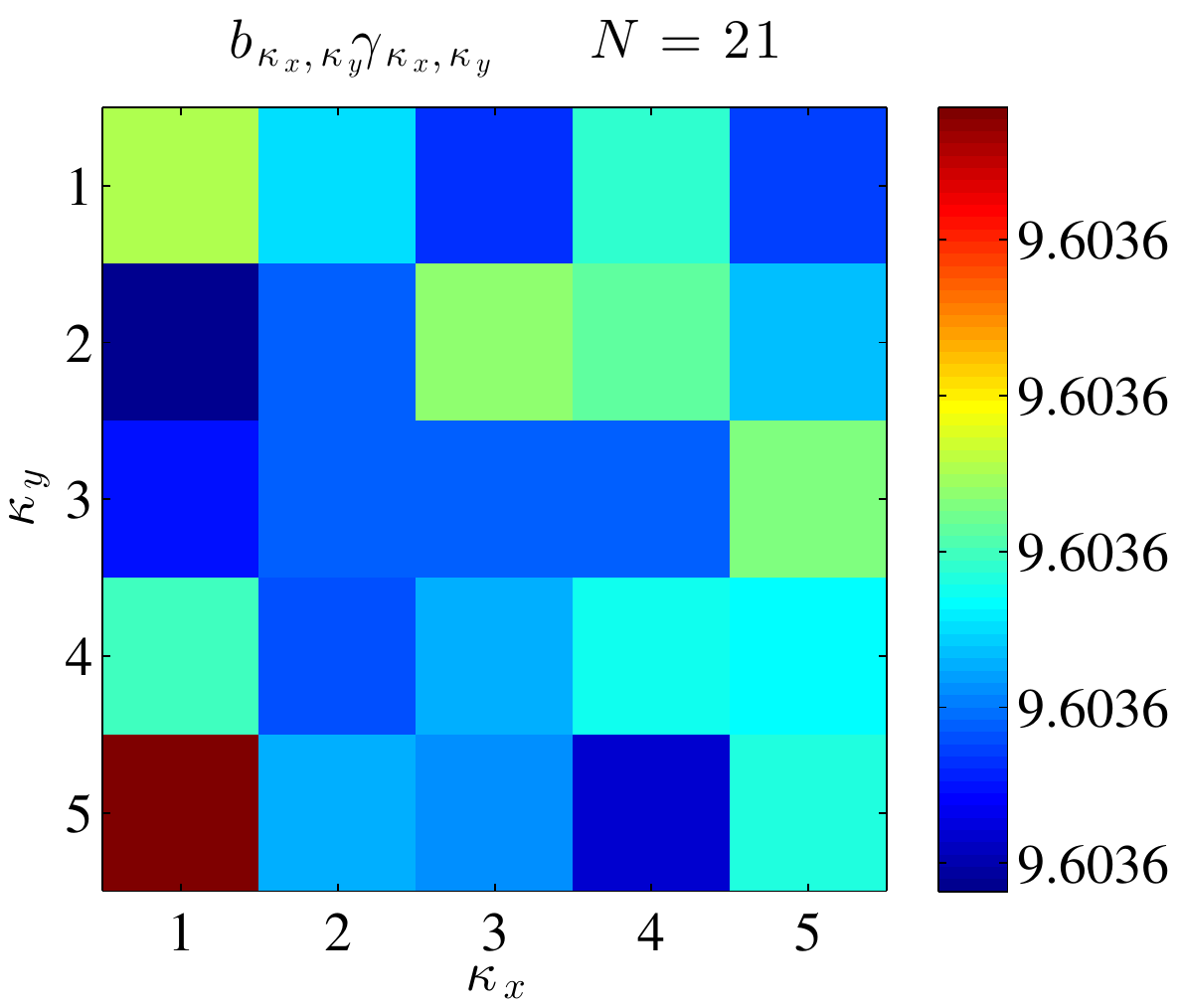}}

    \subfloat[(d)$V=-16.5$]{\includegraphics[width=0.31\textwidth]{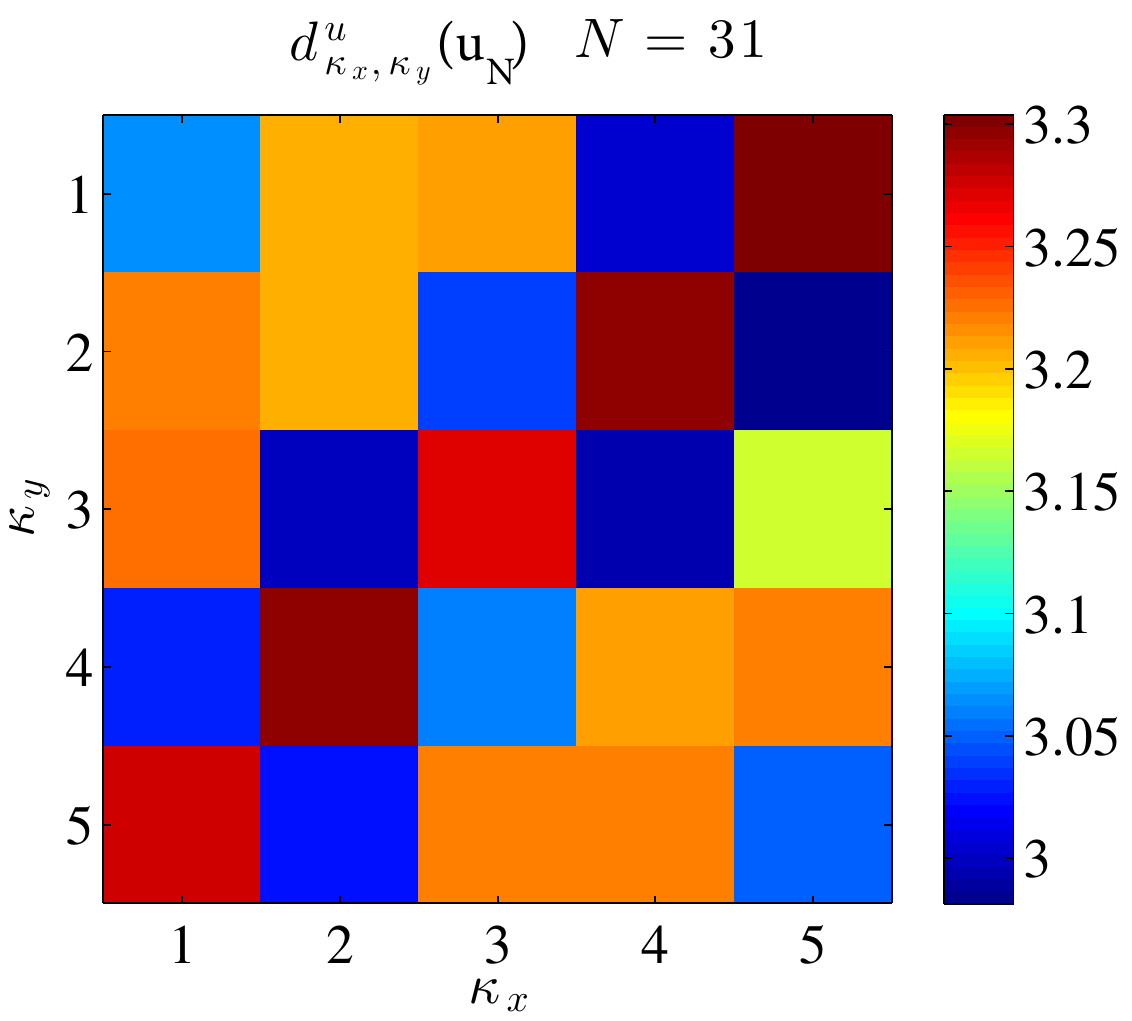}}
    \subfloat[(e)$V=-16.5$]{\includegraphics[width=0.33\textwidth]{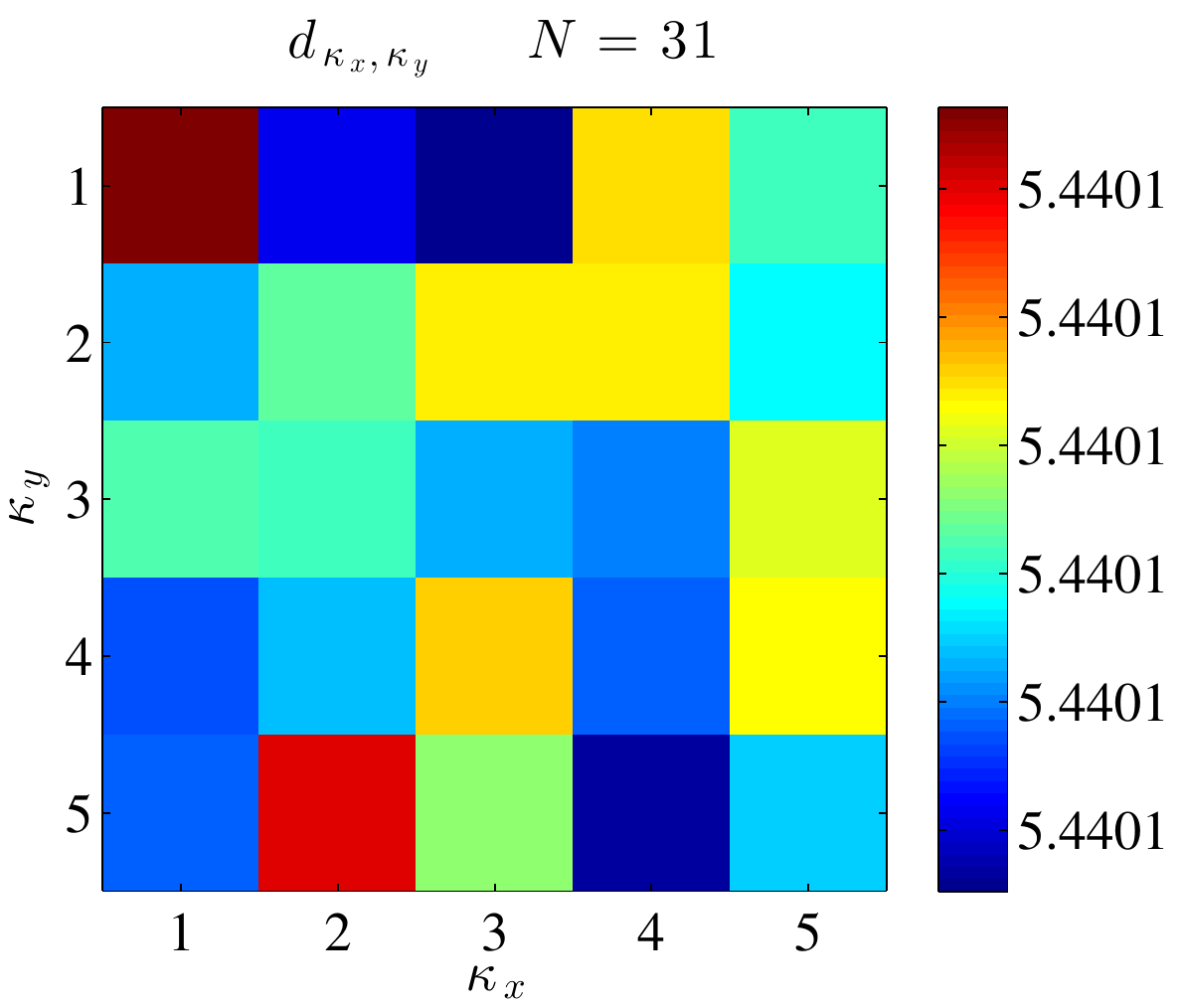}}
    \subfloat[(f)$V=-16.5$]{\includegraphics[width=0.32\textwidth]{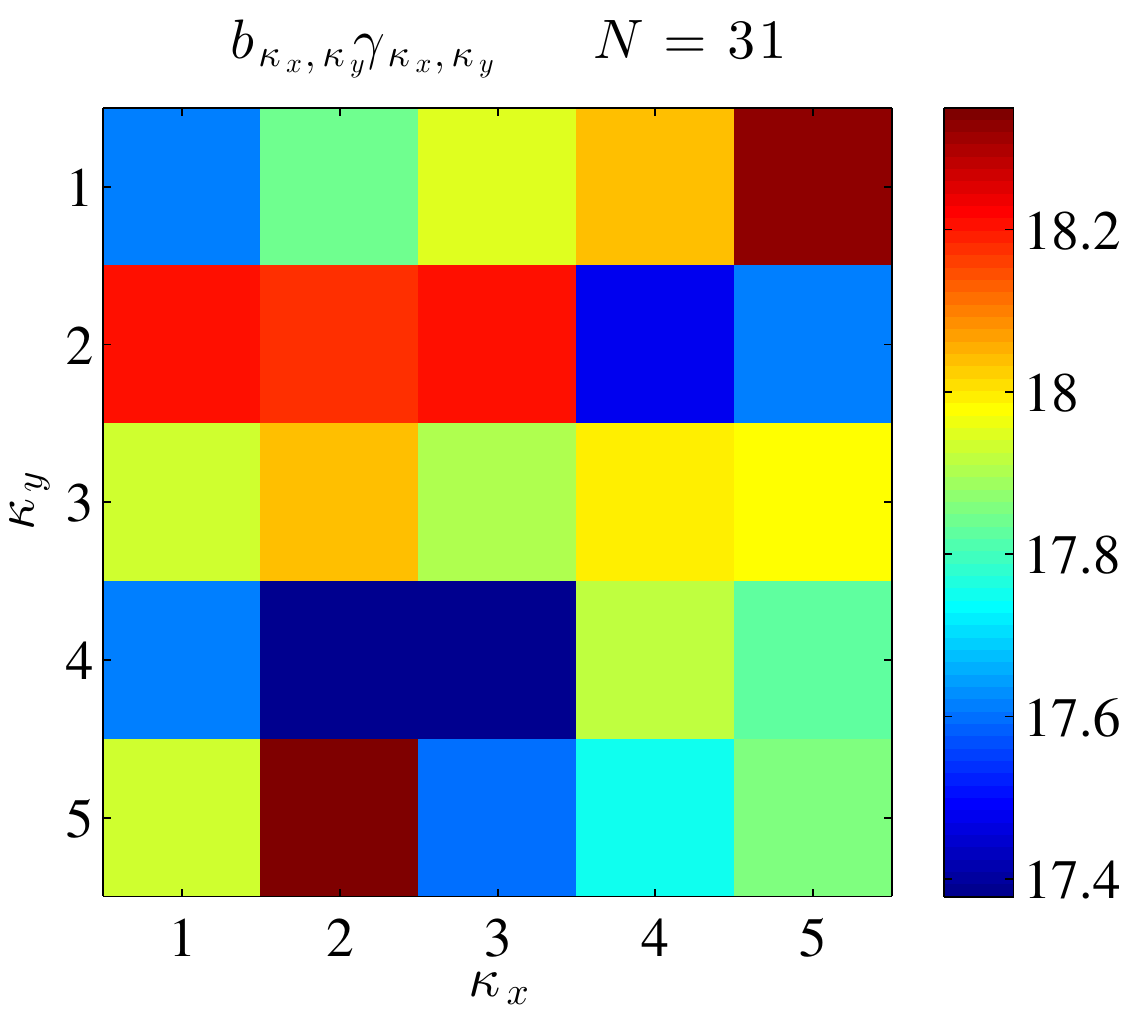}}

    \subfloat[(g)$V$ Gaussian]{\includegraphics[width=0.32\textwidth]{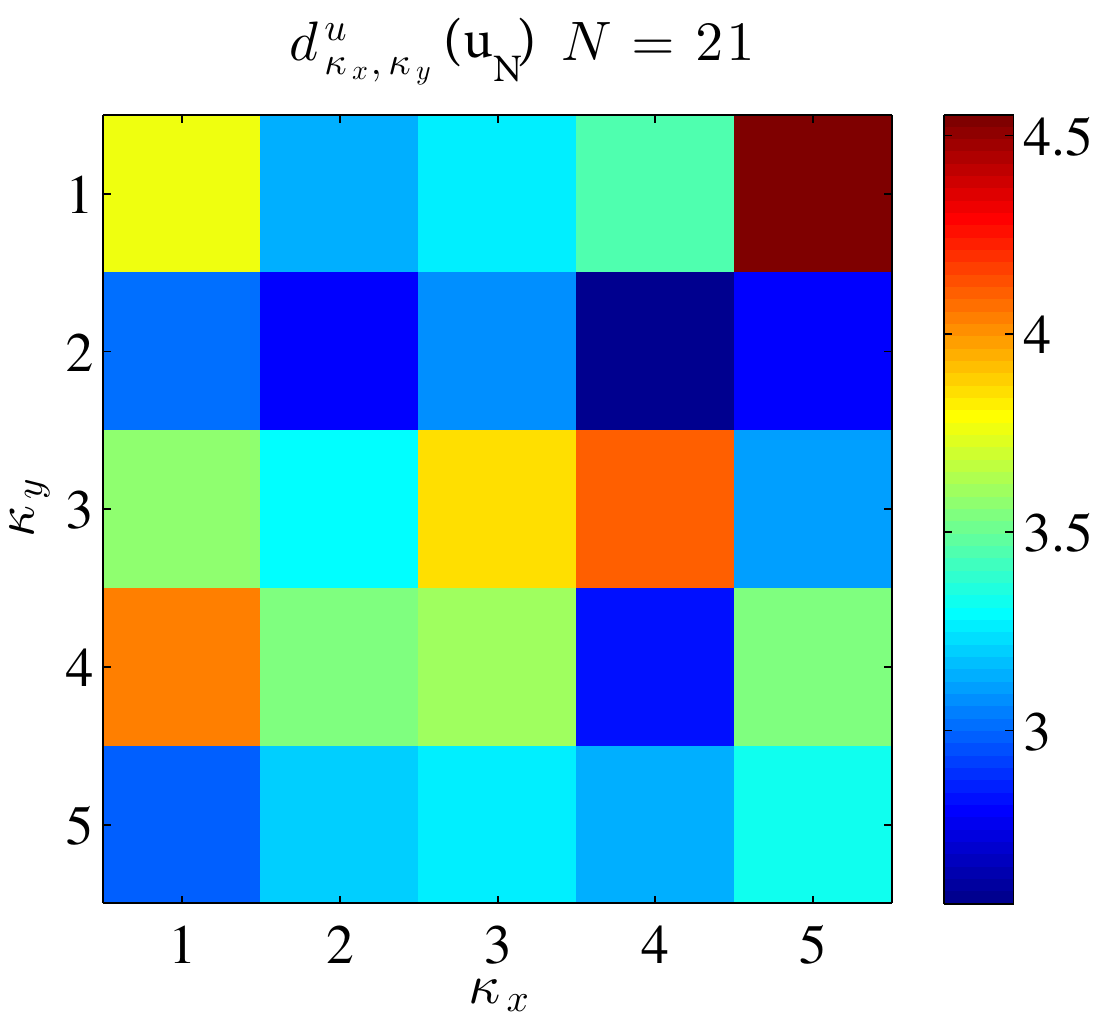}}
    \subfloat[(h)$V$ Gaussian]{\includegraphics[width=0.32\textwidth]{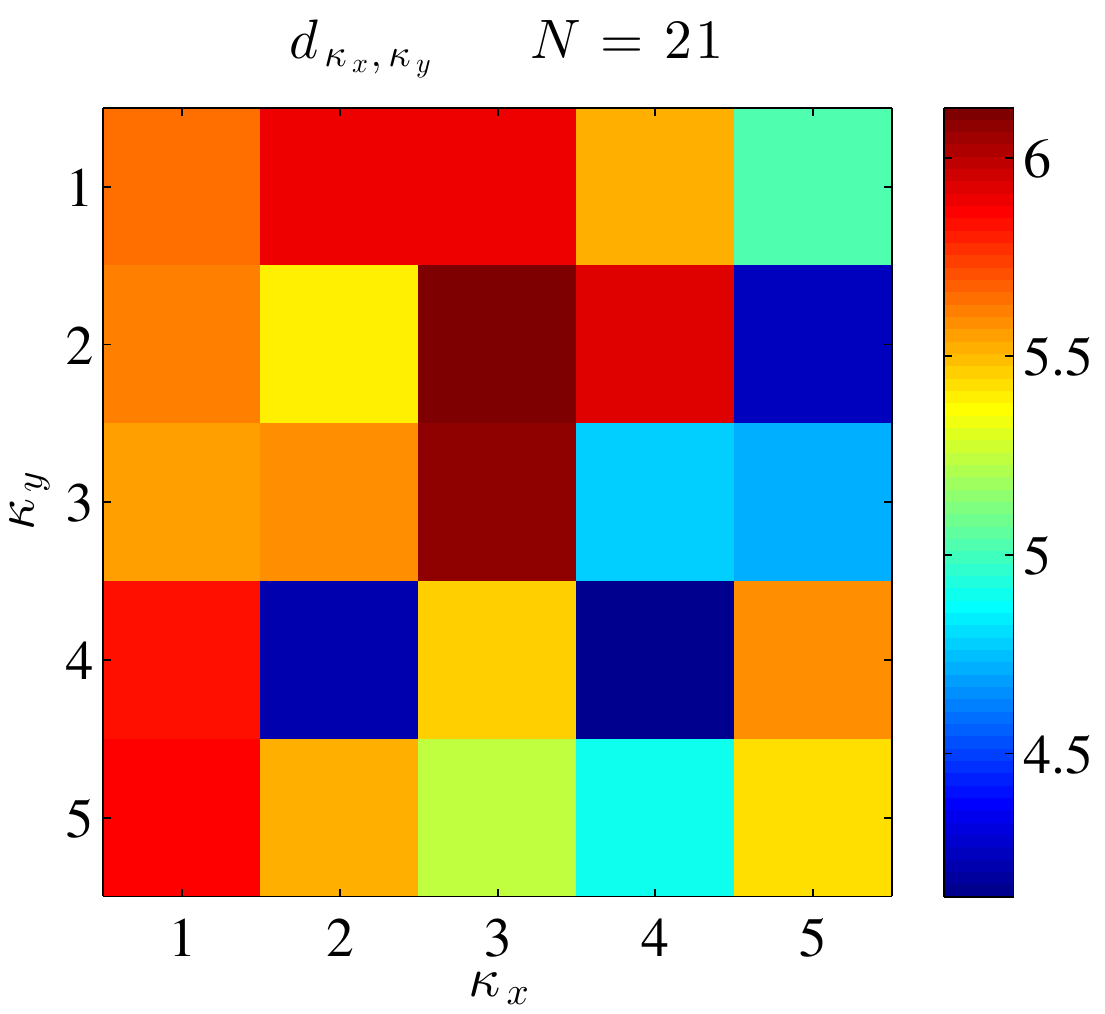}}
    \subfloat[(i)$V$ Gaussian]{\includegraphics[width=0.32\textwidth]{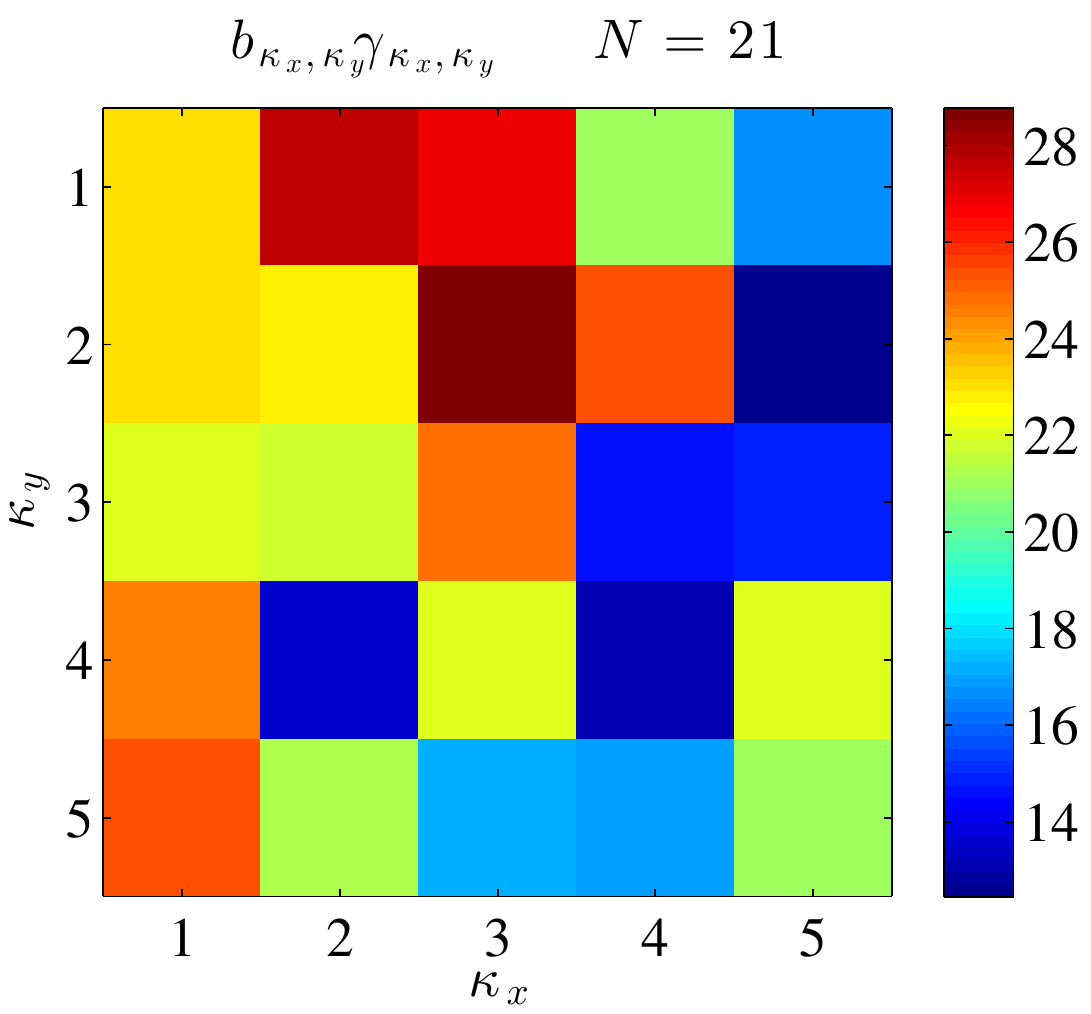}}
  \end{center}
  \caption{Comparison of \REV{$\dku(\uN)$}, $\dkN$ and $\rbk \gk$ for 2D test
  problems for (a-c) the positive definite case $V=0.01$ (d-f) the
  indefinite case $V=-16.5$ (g-i) the indefinite case with $V$ given by
  the sum of negative Gaussians in Fig.~\ref{fig:uuerrIndef2D} (a). }
  \label{fig:dkuComp2D}
\end{figure}

Fig.~\ref{fig:dkuComp2D} compare \REV{$\dku(\uN)$}, $\dkN$ and $\rbk \gk$ for the
positive definite case $V=0.01$, the indefinite case $V=-16.5$, and the
indefinite case with $V$ given by the sum of negative Gaussians in
Fig.~\ref{fig:uuerrIndef2D} (a). In all cases, the magnitude of \REV{$\dku(\uN)$}
is comparable to that of $\dkN$.  Furthermore, both \REV{$\dku(\uN)$} and
$\dkN$ are much smaller compared to $\bk\gk$.  Therefore the
effectiveness of the estimator remains unchanged even if \REV{$\dku(\uN)$} is
neglected.  We expect similar results can be observed for systems of
higher dimensionality.

Finally we provide a second justification by comparing the total
contribution of the jump term in the upper bound estimator 
\[
\eta_{J}^{2}=\sum_{\kappa} \eta_{J,\kappa}^2,
\]
and the total contribution of the jump term in the energy norm
\[
E_{J}=\sum_{\kappa} \tfrac{\gk}{2}\|\jmp{u_N}\|_{\partial\kappa}^2.
\]
This is given in Table~\ref{tab:jumpCompare}. It shows that the
approximation $\REV{\dku(\uN)}\approx \dkN$ does not lead to underestimation of the
jump term, which is consistent with the observation in
Fig.~\ref{fig:dkuComp1D} and~\ref{fig:dkuComp2D}.

\begin{table}
  \centering
  \begin{tabular}{c|c|c|c}
    \hline
    Problem & $N$ & $E_{J}$ & $\eta_{J}^{2}$\\
    \hline
    1D $V=0.01$ &  $7$ & $2.0179\times 10^{-8}$ & $2.0182\times
    10^{-8}$\\
    2D $V=0.01$ & $21$ & $1.2030\times 10^{-5}$ &
    $9.1593\times 10^{-5}$ \\
    1D Gaussian & $11$ &  $6.4687\times 10^{-11}$ & $6.4697\times
10^{-11}$\\
    2D $V=-16.5$ & $31$ & $4.7352\times 10^{-3}$ & $5.6649\times 10^{-2}$
    \\
    2D Gaussian & $21$ & $1.6226\times 10^{-3}$ & $2.8348\times
10^{-2}$\\
    \hline
  \end{tabular}
  \caption{Comparison of the total contribution of the jump term  in the
  estimator $\eta_{J}^{2}$, and the total contribution of the jump term
  in the energy error $E_{J}$.}
  \label{tab:jumpCompare}
\end{table}


%
%

\FloatBarrier

\section{Conclusion and future work}
\label{sec:Conclusion}

We present the first systematic work for deriving a posteriori error
estimates for general non-polynomial basis functions in a interior
penalty discontinuous Galerkin (DG) formulation for solving second order
linear PDEs.   The estimates not only serve to quantify the error
sharply for a given computation, but also can lead an adaptive algorithm
to refine the elements non-uniformly by adding (or even
removing/coarsening) basis functions to certain elements.  This allows a
best approximation for a given number of degrees of freedom in order to
reduce the computing time even when relatively few degrees of freedom
are employed.
A non-uniform distribution of the number of local basis functions is in
this case mandatory to develop powerful solvers, in particular when
inhomogeneous data of the PDE is involved.  It turns out that the
standard polynomial $hp$ DG-method may benefit from this analysis as it
involves numerically computed constants.

Our analysis requires the exact solution to lie in $H^2(\kappa)$ for
each element $\kappa$ which may seem limiting when dealing with
a posteriori estimates for Poisson's equation as a uniform
refinement leads to optimal convergence rates in the asymptotic limit.
We remark that despite the above asymptotic reasoning there are
numerous cases where the a posteriori analysis for regular
functions is still interesting, for example if the PDE involves a strong
small-scale character (but still being smooth) either due to strongly
oscillating material coefficients or a wave-like character of the
underlying PDE (Helmholtz equation for instance).  Or, if the data of
the PDE and thus the solution as well has an inhomogeneous character so
that a uniform refinement involves too many degrees of freedom. In this
case, combining the estimates with an adaptive algorithm as outlined
above will result in an optimal balance of degrees of freedom per
element. 

Our framework for developing explicitly computable constants for a
posteriori error estimates are not limited to second order PDEs, nor
it is necessarily limited to discontinuous Galerkin framework. In a
forthcoming publication we will demonstrate the method for eigenvalue
problems.  It is also possible to generalize the method to multiscale
methods and reduced basis methods.


\section*{Acknowledgments}

This work was partially supported by the Scientific Discovery through Advanced
Computing (SciDAC) program, and by the Center for Applied Mathematics for Energy
Research Applications (CAMERA) funded by U.S. Department of Energy, Office of
Science, Advanced Scientific Computing Research and Basic Energy Sciences
(L. L.). L. L. would like to thank the hospitality of the Jacques-Louis
Lions Laboratory (LJLL) during his visit. We sincerely thank Yvon Maday
for thoughtful suggestions and critical reading of the paper.

\section*{Appendix}

\begin{prop}
  \label{prop:rbk1D}
  Let $\kappa=[a,b]$ be a 1D element and
  $\VN(p;\kappa)=\mathrm{span}\{x^{j}, j\le p \}$ be the function space
  spanned by polynomials with degree less than or equal to $p$. Then
  $\forall p\ge 2$, $\rbk = 0$.
\end{prop}
\begin{proof}
  Define $c=(a+b)/2$. For any $v\in H^1(\kappa)$, $v\perp\VN(p;\kappa)$
  with $p\ge 2$, we have
  \[
    (v,1)_{\star,\kappa} = 0, \quad (v,(x-c))_{\star,\kappa} = 0, \quad
    (v,(x-c)^2)_{\star,\kappa} = 0.
  \]
  Using the definition of the inner product $(\cdot,
  \cdot)_{\star,\kappa}$ 
  \[
  \int_{a}^{b} v(x) \,dx = 0, \quad \int_{a}^{b} v'(x) \,dx = 0, 
  \quad \int_{a}^{b} v'(x) (x-c) d\,x = 0.
  \]
  With integration by parts, we have $v(a) = v(b) = 0$. Therefore
  $\|v\|_{\partial\kappa} = 0$. Using the definition of $\rbk$ we obtain
  $\rbk=0$.
\end{proof}


\begin{thebibliography}{10}

\bibitem{Ainsworth:2012kv}
{\sc M.~Ainsworth and R.~Rankin}, {\em {Technical note: A note on the selection
  of the penalty parameter for discontinuous Galerkin finite element schemes}},
  Numer. Methods Partial Differential Equations, 28 (2012), pp.~1099--1104.

\bibitem{AmaraDjellouliFarhat2009}
{\sc M.~Amara, R.~Djellouli, and C.~Farhat}, {\em Convergence analysis of a
  discontinuous {G}alerkin method with plane waves and {L}agrange multipliers
  for the solution of {H}elmholtz problems}, SIAM J. Numer. Anal., 47 (2009),
  pp.~1038--1066.

\bibitem{Arnold1982}
{\sc D.~N. Arnold}, {\em An interior penalty finite element method with
  discontinuous elements}, SIAM J. Numer. Anal., 19 (1982), pp.~742 -- 760.

\bibitem{BabuskaZlamal:73}
{\sc I.~Babu\v{s}ka and M.~Zl{\'a}mal}, {\em Nonconforming elements in the
  finite element method with penalty}, SIAM J. Numer. Anal., 10 (1973), pp.~863
  -- 875.

\bibitem{Baker77}
{\sc G.~A. Baker}, {\em Finite element methods for elliptic equations using nonconforming elements},
Math. Comp., 31 (1977), pp.~45--59.

\bibitem{BaumannOden99}
{\sc C.~E. Baumann and J.~T. Oden}, {\em A discontinuous {$hp$} finite element method for convection-diffusion problems}, Comput. Methods Appl. Mech. Engrg., 175 (1999), p.~311--341.

\bibitem{EEngquist2003}
{\sc W.~E and B.~Engquist}, {\em The heterognous multiscale methods}, Comm.
  Math. Sci., 1 (2003), pp.~87--132.

\bibitem{DD76}
{\sc J.~Douglas, Jr. and T.~Dupont}, {\em Interior penalty procedures for elliptic and parabolic {G}alerkin
  methods}, Computing methods in applied sciences (Second Internat.
  Sympos., Versailles, 1975), 58 (1976), pp.~207--216. 
\bibitem{Epshteyn:2007vz}
{\sc Y.~Epshteyn and B.~Rivi\`ere}, {\em {Estimation of penalty parameters for
  symmetric interior penalty Galerkin methods}}, J. Comput. Appl. Math., 206
  (2007), pp.~843--872.

\bibitem{FrischPopleBinkley1984}
{\sc M.~J. Frisch, J.~A. Pople, and J.~S. Binkley}, {\em Self-consistent
  molecular orbital methods 25. supplementary functions for gaussian basis
  sets}, J. Chem. Phys., 80 (1984), pp.~3265--3269.

\bibitem{GianiHall2012}
{\sc S.~Giani and E.~J.~C. Hall}, {\em An a posteriori error estimator for
  \textit{hp}-adaptive discontinuous {G}alerkin methods for elliptic eigenvalue
  problems}, Math. Mod. Meth. Appl. Sci., 22 (2012), pp.~1250030--1250064.

\bibitem{HiptmairMoiolaPerugia2011}
{\sc R.~Hiptmair, A.~Moiola, and I.~Perugia}, {\em Plane wave discontinuous
  {G}alerkin methods for the 2{D} {H}elmholtz equation: analysis of the
  p-version}, SIAM J. Numer. Anal., 49 (2011), pp.~264--284.

\bibitem{HouWu1997}
{\sc T.~Y. Hou and X.-H. Wu}, {\em A multiscale finite element method for
  elliptic problems in composite materials and porous media}, J. Comput. Phys.,
  134 (1997), pp.~169--189.

\bibitem{HoustonSchotzauWihler2007}
{\sc P.~Houston, D.~Sch{\"o}tzau, and T.~P. Wihler}, {\em Energy norm a
  posteriori error estimation of hp-adaptive discontinuous {G}alerkin methods
  for elliptic problems}, Math. Mod. Meth. Appl. Sci., 17 (2007), pp.~33--62.

\bibitem{Houston:2002uo}
{\sc P.~Houston, C.~Schwab, and E.~S{\"u}li}, {\em {Discontinuous $hp$-finite
  element methods for advection-diffusion-reaction problems}}, SIAM J. Numer.
  Anal., 39 (2002), pp.~2133--2163.

\bibitem{Junquera:01}
{\sc J.~Junquera, O.~Paz, D.~Sanchez-Portal, and E.~Artacho}, {\em Numerical
  atomic orbitals for linear-scaling calculations}, Phys. Rev. B, 64 (2001),
  pp.~235111--235119.

\bibitem{KarakashianPascal2003}
{\sc O.~A. Karakashian and F.~Pascal}, {\em A posteriori error estimates for a
  discontinuous {G}alerkin approximation of second-order elliptic problems},
  SIAM J. Numer. Anal., 41 (2003), pp.~2374--2399.

\bibitem{KayeLinYang2014}
{\sc J.~Kaye, L.~Lin, and C.~Yang}, {\em A posteriori error estimator for
  adaptive local basis functions to solve {Kohn-Sham} density functional
  theory}, Commun. Math. Sci, in press,  (2015).

\bibitem{Knyazev2001}
{\sc A.~V. Knyazev}, {\em Toward the optimal preconditioned eigensolver:
  Locally optimal block preconditioned conjugate gradient method}, SIAM J. Sci.
  Comp., 23 (2001), p.~517.

\bibitem{LinLuYingE2012}
{\sc L.~Lin, J.~Lu, L.~Ying, and W.~E}, {\em {Adaptive local basis set for
  Kohn-Sham density functional theory in a discontinuous Galerkin framework I:
  Total energy calculation}}, J. Comput. Phys., 231 (2012), pp.~2140--2154.

\bibitem{Nitsche71}
{\sc J.~Nitsche}, {\em\"{U}ber ein {V}ariationsprinzip zur {L}\"osung von
  {D}irichlet-{P}roblemen bei {V}erwendung von {T}eilr\"aumen, die keinen
  {R}andbedingungen unterworfen sind},
	Abh. Math. Sem. Univ. Hamburg, 36 (1971), pp.~9--15.
  
\bibitem{Ohlberger2005}
{\sc M.~Ohlberger}, {\em A posteriori error estimates for the heterogeneous
  multiscale finite element method for elliptic homogenization problems},
  Multiscale Model. Simul., 4 (2005), pp.~88--114.

\bibitem{SchotzauZhu2009}
{\sc D.~Sch{\"o}tzau and L.~Zhu}, {\em A robust a-posteriori error estimator
  for discontinuous {G}alerkin methods for convection--diffusion equations},
  Appl. Numer. Math., 59 (2009), pp.~2236--2255.

\bibitem{Schwab1998}
{\sc C.~Schwab}, {\em $p$-and $hp$-Finite Element Methods}, Oxford Univ. Pr.,
  New York, 1998.

\bibitem{StammWihler2010}
{\sc B.~Stamm and T.~Wihler}, {\em {hp}-{O}ptimal discontinuous {G}alerkin
  methods for linear elliptic problems}, Math. Comp., 79 (2010),
  pp.~2117--2133.

\bibitem{TezaurFarhat2006}
{\sc R.~Tezaur and C.~Farhat}, {\em Three-dimensional discontinuous {G}alerkin
  elements with plane waves and {L}agrange multipliers for the solution of
  mid-frequency {H}elmholtz problems}, Int. J. Numer. Meth. Eng., 66 (2006),
  pp.~796--815.

\bibitem{Wheeler78}
{\sc M.~F. Wheeler}, {\em An elliptic collocation-finite element method with interior
  penalties.}, SIAM J. Numer. Anal., 15 (1978), p.~152--161.

\end{thebibliography}
\end{document}